\documentclass[11pt]{article}%
\usepackage{graphicx}
\usepackage{subcaption}
\usepackage{makeidx}
\usepackage{amssymb}
\usepackage[all]{xy}
\usepackage{float}
\usepackage[latin1]{inputenc}
\usepackage{amsfonts}
\usepackage{amsmath}
\usepackage{amssymb}
\usepackage{url}
\usepackage{hyperref}
\usepackage{graphicx}%
\usepackage{mathabx}
\usepackage{wrapfig}

\usepackage{amsthm}
\usepackage{url}
\usepackage{hyperref}
\providecommand{\U}[1]{\protect\rule{.1in}{.1in}}

\usepackage{xcolor}

\newtheorem{prop}{Proposition}[section]
\newtheorem{cor}[prop]{Corollary}
\newtheorem{defi}[prop]{Definition}

\newtheorem{lem}[prop]{Lemma}

\newtheorem{theo}[prop]{Theorem}

\newcommand{\EE}{\mathbb{E}}

\newcommand{\LL}{\mathbb{L}}

\newcommand{\RR}{\mathbb{R}}

\newcommand{\Aa}{ {\cal A }}

\newcommand{\Ca}{ {\cal C }}

\newcommand{\La}{ {\cal L }}

\newcommand{\Ea}{ {\cal E }}
\newcommand{\Sa}{ {\cal S }}
\newcommand{\Ra}{ {\cal R }}

\newcommand{\Fa}{ {\cal F }}

\newcommand{\Qa}{ {\cal Q }}

\newcommand{\Xa}{ {\cal X }}
\newcommand{\Ma}{ {\cal M }}

\newcommand{\Pa}{ {\cal P }}
\newcommand{\Za}{ {\cal Z }}
\newcommand{\Ya}{ {\cal Y }}
\newcommand{\Wa}{ {\cal W }}

\newcommand{\point}{\mbox{\LARGE .}}

\newcommand{\cqfd}{\hfill\blbx \\}
\def\blbx{\hbox{\vrule height 5pt width 5pt depth 0pt}\medskip}

\def \RR{\mathbb{R}}
\def \SS{\mathbb{S}}
\def \EE{\mathbb{E}}

\def \LL{\mathbb{L}}

\def \WW{\mathbb{W}}

\def\ricc{\mbox{\rm Ricc}}

\begin{document}

  \title{On the stability and the uniform propagation of chaos properties of Ensemble Kalman-Bucy filters }
  \author{P. Del Moral, J. Tugaut}


\maketitle

\begin{abstract}
The Ensemble Kalman filter is a sophisticated and powerful 
data assimilation method  for filtering high dimensional problems arising in fluid mechanics and
geophysical sciences. This Monte Carlo method can be interpreted as a mean-field 
McKean-Vlasov type particle interpretation of the Kalman-Bucy diffusions. In contrast to more conventional particle filters and
nonlinear Markov processes these models are designed in terms of a diffusion process with a diffusion matrix that depends on particle covariance matrices. 

Besides some recent advances on the stability of nonlinear Langevin type diffusions with drift interactions,  the long-time behavior of 
models with interacting diffusion matrices and conditional distribution interaction functions has never been discussed in the literature. One of the main contributions of the article is to initiate the study of this new class of models The article presents a series of new functional inequalities to quantify the stability of these nonlinear diffusion processes. 

In the same vein, despite some recent contributions on the convergence
of the Ensemble Kalman filter when the number of sample tends to infinity very little is known on stability and the long-time behaviour
of these mean-field interacting type particle filters. 
The second contribution of this article is to provide uniform propagation of chaos properties as well as
 $\LL_n$-mean error 
estimates w.r.t. to the time horizon. Our regularity condition is also shown to be sufficient and necessary for the uniform convergence of the Ensemble Kalman filter.

The stochastic analysis developed in this article is based on an original combination of functional inequalities and Foster-Lyapunov techniques
with coupling, martingale techniques, random matrices and spectral analysis theory. \\



\end{abstract}


\section{Introduction}
\subsection{The Ensemble Kalman filter}
The Ensemble Kalman filter ({\em abbreviated EnKF}) has been introduced by G. Evensen 
 in the seminal article~\cite{evensen-intro} published  in 1994. In the last two decades the EnKF has became one of the main numerical technique
for solving high dimensional forecasting and data assimilation problems, particularly in ocean and atmosphere
sciences~\cite{allen,lisa,majda,kalnay,ott}, weather forecasting~\cite{anderson-jl,anderson-jl-2,burgers,houte},  environmental and ecological statistics~\cite{eknes,johns},
as well as in oil reservoir simulations~\cite{evensen-reservoir,nydal,seiler,skj,weng}, and many others. We also refer the reader to~\cite{bittanti,einicke}
for recent reviews on Riccati equations, estimation and linear filtering techniques.

The mathematical foundations and the convergence of the EnKF is more recent. It has started in 2011 with the independent pioneering works of 
F. Le Gland, V. Monbet and V.D. Tran~\cite{legland}, and the one by J. Mandel, L. Cobb, J. D. Beezley~\cite{mandel}.
These articles provide $\LL_n$-mean error estimates for discrete time EnKF and show that the EnKF converges towards the Kalman filter as the number of samples tends to infinity. In a more recent study by X. T. Tong, A. J. Majda and D. Kelly
the authors analyze the long-time behaviour and the ergodicity of discrete generation EnKF using Foster-Lyapunov techniques ensuring that the filter is asymptotically stable w.r.t. any erroneous initial condition~\cite{tong}. 
These important properties ensure that the EnKF has a single invariant measure 
and initialization errors of the EnKF will not dissipate w.r.t. the time parameter.

Beside the importance of these properties, the only ergodicity of the particle process does not give any information of the convergence and the accuracy of 
the EnKF towards the optimal filter as the number of samples tends to infinity.

One of the main objective of this article is to analyze this convergence and quantify the fluctuation of errors on large-time horizon. We provide 
 uniform $\LL_n$-mean error estimates  w.r.t. the time parameter for the sample mean as well as for the 
 sample covariance matrices. Incidentally, the stochastic analysis we have developed also allows to quantify the 
 stability properties of the Kalman-Bucy filter and the corresponding matrix valued Riccati equations. These estimates are deduced from the stability properties of a nonlinear diffusion interpretation of the Kalman-Bucy equations.

To better connect this work with existing literature on nonlinear Markov processes and particle methods
we emphasize that the EnKF can be seen as a mean-field particle interpretation 
 of a nonlinear McKean-Vlasov type diffusion. These probabilistic models were introduced in the end of the 60s
 by H. P. McKean~\cite{mckean}. For a detailed discussion on these models and their application domains we refer the reader to the lecture notes of A. S. Sznitman~\cite{sznitman}, the ones by S. M\'el\'eard~\cite{meleard2}, and the more research monograph~\cite{dm-crc2013}. 
 
 The refined convergence as well as the long-time behaviour of nonlinear diffusion processes is still an active research area. When the interaction function only enters in the drift part of the diffusion
 several results including uniform estimates w.r.t. the time horizon are available~\cite{bolley,cattiaux,tugaut-1,tugaut-2,malrieu-1}. Most of these works are based on power and sophisticated 
 coupling methods, nonlinear semigroup analysis, as well as Gamma-two type techniques and optimal transport theory. 
 
In our context, the Kalman-Bucy filter and the Riccati equation represents the evolution equations of the mean and 
 the covariance matrices of the random states of a nonlinear diffusion process. We shall call this process the Kalman-Bucy diffusion. The diffusion part of this class of processes depend on the covariance matrix of its random states. The recent techniques developed in nonlinear Markov processes theory are not suited to analyze the stability of these complex nonlinear processes. To the best of our knowledge the long-time behaviour of such nonlinear diffusions with
covariances matrices depending on the distribution of the random states remains an open and important research question. 

In the present article we initiate the study of the stability of this class of nonlinear diffusion models. We present a series of functional inequalities to quantify the stability of these nonlinear diffusion processes. We  also analyze the exponential stability of these processes w.r.t. Wasserstein distances and relative entropy inequalities. The stability properties of the Kalman-Bucy filter are deduced by a direct application of Jensen type inequalities.

At the level of the particle population model, the EnKF also belongs to the class of mean-field type particle filters. 
The stochastic analysis of particle filters and related diffusions Monte Carlo schemes is rather well understood, see for instance~\cite{dm-springer2004,dm-crc2013} and the references therein.
Nevertheless the EnKF strongly differs from particle filters or Sequential Monte Carlo methods currently used in nonlinear filtering theory, Bayesian inference and computational physics.
 Roughly speaking, the EnKF is designed to approximate the Kalman filter (as well as the extended Kalman filter) for high dimensional problems. In the reverse angle, particle filters are designed to estimate the nonlinear filtering equation, and to sample sequentially according to the flow of conditional distributions. In continuous time settings, the EnKF is an interacting diffusion while particle filters are interacting jump particle systems. As a result none of the techniques developed in particle filtering theory applies to analyze the fluctuations of the EnKF uniformly w.r.t. the time horizon.
 It is clearly not the scope of this article to compare in full details these two particle filtering methods. For a more thorough discussion on particle filtering techniques we  refer the reader to~\cite{dm-springer2004,dm-crc2013}, and the references therein.

 In the same vein, the stochastic analysis developed so far in the literature on more general classes of mean-field particle methods 
 cannot be used to analyze the uniform convergence of particle approximating schemes involving interacting covariances matrices.
As mentioned above, the EnKF belongs to this class of nonlinear diffusions with a mean-field particle interpretation based on interacting covariance matrices of multi-dimensional particles. 
To the best of our knowledge the uniform propagation of chaos estimates developed in the present article seems to be the first result of this type 
for this class of nonlinear diffusions.

To derive these uniform estimates we develop a novel stochastic fluctuation analysis which combines 
Foster-Lyapunov techniques with matrix valued martingale methods, as well as random matrices and spectral analysis theory. The central idea is to take advantage of the linear-Gaussian 
structure of the filtering problem to enter the stability properties of the signal process and the (nonlinear) Riccati matrix-valued equation
into the fluctuation analysis of the EnKF. We also prove that the stability property of the signal is a sufficient and necessary condition to obtain uniform propagation of chaos estimates.

\subsection{Organization of the article}
The article is organized as follows:

Section~\ref{desc-sec-intro}
is dedicated to the description of the Kalman-Bucy filter, 
 the nonlinear diffusion process interpretation of the filter, as well as the mean-field EnKF particle algorithm. In Section~\ref{statement-sec-intro} we state
the main theorems of the article. 

The first one shows that the sample mean and the random interacting
covariance matrices of the EnKF satisfy the same equation as the EnKF and the Riccati equation up to some fluctuation martingales whose angle brackets only depends on the sample covariance matrices. These diffusion equations in matrix spaces
are pivotal as they
allow to analyze the fluctuations of the EnKF using Foster-Lyapunov and martingale techniques combined with trace and spectral type inequalities.
 
The second theorem provides uniform convergence
and propagations of chaos  estimates  w.r.t. the time parameters.

Section~\ref{section-regularity} provides a detailed discussion on our regularity conditions. 
In section~\ref{sous-section-section-regularity} we analyze the stability properties and the 
catastrophic divergence issues of EnFK filters in terms of global divergence regions and ill-conditioned filtering problems.
We analyze the propagations of the fluctuations induced by the sample covariance matrices in terms of observer-type filters
and stochastic Ornstein-Ulhenbeck diffusions. 
In control theory, the terminology "observer" is often restricted to deterministic models. 

As its name indicates a stochastic observer
is a stochastic process that uses sensory history to estimate the true signal; the randomness comes from the fact that the perturbations of the
sensor are random. We design and we analyze the long time behavior of a class of 
stochastic observer driven by stochastic covariance matrices. We also discuss some pivotal semigroup contraction properties
in terms of log-norms of matrices. Several illustrations are provided in section~\ref{observability-section}. 

Section~\ref{stab-section-p} discusses 
the stability properties of Kalman-Bucy diffusions. Section~\ref{global-contraction-section} is dedicated with
uniform contraction inequalities for the nonlinear semigroups associated with the Riccati equation and
 Kalman-Bucy diffusions. Section~\ref{local-contraction} presents some local functional inequalities
 to estimate the fluctuations of the models around their steady state version w.r.t. the Wasserstein distance and the relative entropy.

The remainder of the article is mainly concerned with the proof of the main theorems presented in Section~\ref{statement-sec-intro} and Section~\ref{stab-section-p}:

Section~\ref{preliminary-sec} presents some technical preliminary results used in the further development of the article.
Section~\ref{some-uniform-moments-sec} shows that our regularity conditions that ensures the uniform convergence
of the EnKF is sharp and cannot be relaxed. The section also provides some uniform convergence estimates on the filter,
the signal states, and the Riccati equation. It also presents some semigroup estimates and related trace inequalities of current use in this study.
Section~\ref{riccati-section} is dedicated to the Riccati equation. We 
analyze the explicit solution in the one dimensional case and we present a trace type comparison lemma
to analyze multivariate models.

Section~\ref{sec-EnFK} is concerned with the stochastic analysis of the EnKF.
Section~\ref{sec-EnFK-ss1} is dedicated to the proof of the stochastic differential equations EnKF 
sample mean and the particle covariance matrices. Section~\ref{unif-moment-EnKF-ss2} is dedicated to uniform
moments estimates for the trace of the particle covariance matrices and the random states of the EnKF. These results are deduced from a technical lemma, of its own interest; combining Foster-Lyapunov with martingale techniques to control the moments of Riccati type
stochastic differential equations uniformly w.r.t. the time horizon. 

Section~\ref{quantitative-section} is mainly concerned with the detailed proofs of the  uniform propagation of chaos theorem presented in Section~\ref{statement-sec-intro}. 

The final section, Section~\ref{conclusion}, presents a brief summary of the contributions of the article and
proposes an avenue of open research projects.

\subsection{Some basic notation and preliminary results}
This section provides with some notation and terminology used in several places in the article.
Given some random variable $Z$ with some probability measure $\mu$ and some function
$f$ on some product space $\RR^r$, we let $\mu(f)=\EE(f(Z))=\displaystyle\int~f(x)~\mu(dx)$ be the integral of $f$ w.r.t. $\mu$ or the expectation of $f(X)$.  As a rule any multivariate random variable, say $Z$,
is represented by a column vector and we use the transposition operator $Z^{\prime}$ to denote the row vector.
Given a distribution on some product space $\RR^r$ and some measurable  function $f$ from $\RR^r$ into $\RR^r$
we set $\eta(f)=\left(\eta(f_i)\right)_{1\leq i\leq r}$ the column vector with entries $\eta(f_i)$ where $f_i$ stands for the $i$-th coordinate mapping from $\RR^r$ into $\RR$. We also denote by $a_+=\max{(a,0)}$ the real part of a number $a\in\RR$.

We let $\left\Vert\point\right\Vert$ be the Euclidean norm on $\RR^{r}$, for some $r\geq 1$. We denote by $\SS_r$ the set of $(r\times r)$ symmetric matrices
with real entries, and by $\SS_r^+$ the subset of positive definite matrices. We let $\mbox{\rm Spec}(A)$ be the set of 
eigenvalues of a square matrix $A$. With a slight abuse of notation we denote 
by $Id$ the $(r\times r)$ identity matrix, for any $r\geq 1$.

We often denote by $\lambda_i(A)$, with $1\leq i\leq r$,
the non increasing sequence of eigenvalues of a symmetric $(r\times r)$-matrix $A$.
We also often denote by $\lambda_{\tiny min}(A)=\lambda_{r}(A)$ and $\lambda_{\tiny max}(A)=\lambda_{1}(A)$ the minimal and the maximal eigenvalue. We also set 
$A_{\tiny sym}:=(A+A^{\prime})/2$ for any  $(r\times r)$-square matrix $A$.
We recall that the norm $\Vert A\Vert$ and logarithmic norm $\mu(A)$ of an $(r_1\times r_1)$-square matrix $A$ 
are defined by $\Vert A\Vert:=\sup_{\|x\|=1}{\Vert A x\Vert}$ and
\begin{equation}\label{def-log-norm}
\begin{array}{rcl}
\mu(A)&:=&\inf{\left\{\alpha~:~\forall x~~\langle x,Ax\rangle\leq \alpha \left\Vert x\right\Vert^2\right\}}=\lambda_{\tiny max}\left(A_{\tiny sym}\right)\\
&=&\inf{\left\{\alpha~:~\forall t\geq 0~~\Vert \exp{(At)}\Vert_2\leq \exp{(\alpha t)}\right\}}
\end{array}
\end{equation}
The above equivalent formulations show that 
$$
\mu(A)\geq \varsigma(A):=\max{\left\{\mbox{\rm Re}(\lambda)~:~\lambda\in \mbox{\rm Spec}(A)\right\}}
$$
where $\mbox{\rm Re}(\lambda)$ stands for the real part of  the eigenvalues $\lambda$.
The parameter $\varsigma(A)$ is often called the spectral abscissa of $A$.
Also notice that $A_{\tiny sym}$ is negative semi-definite as soon as $\mu(A)<0$. 
The Frobenius matrix norm of a given $(r_1\times r_2)$ matrix $A$ is defined by
$$
\left\Vert A\right\Vert_{F}^2=\mbox{\rm tr}(A^{\prime}A)
\qquad\mbox{\rm
with the trace operator $\mbox{\rm tr}(\point)$.}
$$
If $A$ is a matrix $r\times r$, we have $\left\Vert A\right\Vert_{F}^2=\sum_{1\leq i,j\leq r}A(i,j)^2\geq \Vert A\Vert^2$.

We also need to consider the $n$-th Wasserstein distance between two probability measures
$\nu_1$ and $\nu_2$ on $\RR^r$ defined by
$$
\WW_n(\nu_1,\nu_2)=\inf{\left\{\EE\left(\Vert Z_1-Z_2\Vert^n\right)^{\frac{1}{n}}
\right\}}.
$$
The infimum in the above displayed formula is taken of all pair of random variable $(Z_1,Z_2)$
such that $\mbox{\rm Law}(Z_i)=\nu_i$, with $i=1,2$. We denote by $ \mbox{\rm Ent}\left(\nu_1~|~\nu_2\right)$
the Boltzmann-relative entropy
$$
 \mbox{\rm Ent}\left(\nu_1~|~\nu_2\right):=\int~\log\left(\frac{d\nu_1}{d\nu_2}\right)~d\nu_1\quad\mbox{if $\nu_1\ll\nu_2$, and $+\infty$ otherwise.}
$$

The state transition matrix associated with a smooth flow of $(r\times r)$-matrices $A~:~u\mapsto A_u$ is denoted by
$$
\Ea_{s,t}(A)=\exp{\left[\oint_s^t A_u~du\right]}\Longleftrightarrow \partial_t \Ea_{s,t}(A)=A_t~\Ea_{s,t}(A)\quad\mbox{\rm and}\quad
\partial_s \Ea_{s,t}(A)=-\Ea_{s,t}(A)~A_s
$$
for any $s\leq t$, with $\Ea_{s,s}=Id$, the identity matrix. Equivalently in terms of the fundamental solution matrices
$\Ea_t(A):=\Ea_{0,t}(A)$ we have
$
\Ea_{s,t}(A)=\Ea_t(A)\Ea_s(A)^{-1}
$.
Observe that for any $s\leq r\leq t$ the exponential semigroup property
$$
\exp{\left[\oint_s^t A_u~du\right]}=\exp{\left[\oint_r^t A_u~du\right]}\exp{\left[\oint_s^rA_u~du\right]}.
$$
The following technical lemma provides a pair of semigroup estimates  of the state transition matrices
associated with a sum of drif-type matrices. 
\begin{lem}[Perturbation lemma]\label{perturbation-lemma-intro}
Let $A~:~u\mapsto A_u$ and $B~:~u\mapsto B_u$ be some smooth flows of $(r\times r)$-matrices.
For any $s\leq t$ we have
  $$
 \left\Vert   \Ea_{s,t}(A+B)\right\Vert_2\leq \exp{\left(\int_s^t\mu(A_u)~du+\int_s^t~\Vert B_u\Vert_2~du\right)}
 $$
In addition, for any  matrix norm $\Vert\cdot\Vert$ we have
 $$
  \left\Vert   \Ea_{s,t}(A+B)\right\Vert\leq \alpha_A \exp{\left[-\beta_A(t-s)+\alpha_A\int_s^t \Vert B_u\Vert~du\right]}
 $$
 as soon as
  $$
 \forall 0\leq s\leq t\qquad \Vert  \Ea_{s,t}(A) \Vert\leq \alpha_A~\exp{\left(-\beta_A~(t-s)\right)}.
 $$

\end{lem}
These estimates are probably well known but we have not found a precise
reference. For the convenience of the reader, the detailed proof of this lemma is housed in the appendix, on page~\pageref{proof-lemma-intro}.
For time homogeneous matrices $A_t=A$, the state transition matrix reduces to the conventional matrix exponential  $\Ea_{s,t}(A)=e^{(t-s)A}=
\Ea_{t-s}(A)$. 

The norm of $\Ea_t(A)$ can be estimated in various ways: The first one is based on the Jordan decomposition
$T^{-1}AT=J$ decomposition of the matrix $A$ in terms of $k$ Jordan blocks associated with the eigenvalues
with multiplicities $m_i$, with $1\leq i\leq k$. In this situation, we have the Jordan type estimate
\begin{equation}\label{first-estimate-Jordan}
e^{\varsigma(A)t}\leq \displaystyle\Vert \Ea_t(A)\Vert_2\leq \kappa_{{\tiny Jor},t}(T) ~e^{\varsigma(A)t}
\end{equation}
with
$$
\kappa_{{\tiny Jor},t}(T)=\left(\vee_{0\leq j<n}\frac{t^j}{j!}\right)~ \Vert T\Vert_2\Vert T^{-1}\Vert_2
\quad \mbox{\rm and}\quad n:=\vee_{1\leq i\leq k} m_i
$$
Observe that $\kappa_{{\tiny Jor},t}(T)$ depends on the time horizon $t$ as soon as $A$ is not of full rank. In addition, whenever $A$ is close to singular, the conditioning number $\mbox{\rm cond}(T):= \Vert T\Vert_2\Vert T^{-1}\Vert_2$ 
tends to be very large.

A second strategy is based on Schur decomposition $U^{\prime}AU=D+T$ in terms of an unitary matrix $U$, with
$D=\mbox{\rm diag}(\lambda_1(A),\ldots;\lambda_r(A))$ and a strictly triangular matrix $T$ s.t. $T_{i,j}=0$ for any $i\geq j$. In this case we have the Schur type estimate
\begin{equation}\label{estimate-Schur}
\Vert \Ea_t(A)\Vert_2\leq  \kappa_{{\tiny Sch},t}(T)~e^{\varsigma(A)t}
~\quad \mbox{\rm with}\quad \kappa_{{\tiny Sch},t}(T):=\sum_{0\leq i\leq r}\frac{(\Vert T\Vert t)^i}{i!}
\end{equation}
The proof of these estimates can be found in ~\cite{vanloan-19,vanloan}. 
In both cases for any $\epsilon\in ]0,1]$ and  any $t\geq 0$ we have
\begin{equation}\label{common-estimate-Jordan-Schur}
e^{\varsigma(A)t}\leq\Vert\Ea_t(A)\Vert_2\leq \kappa(\epsilon)~~e^{(1-\epsilon)\varsigma(A)t}~
\end{equation}
for some constants $\kappa(\epsilon)$ whose values only depend on the parameters $\epsilon$. When $A$
is  asymptotically stable; that is all its eigenvalues have negative real parts, for any positive definite matrix $B$ we have
$$
e^{\varsigma(A)t}\leq\Vert\Ea_t(A)\Vert_2\leq \mbox{\rm cond}(T_B)~\exp{\left[-t/\Vert B^{-1/2} T_B~B^{-1/2}\Vert\right]}
$$
 with the positive definite matrix
$$
T_B=\int_0^{\infty}~e^{A^{\prime}t}~B~e^{At}~dt\Longleftrightarrow A^{\prime}T+TA=-B
$$
In this case, we have $-1/\Vert B^{-1/2} T~B^{-1/2}\Vert\geq 2\varsigma(A)$. The proof of these estimates can be found in ~\cite{kresemir} (theorem 13.6 and exercise 13.11). 

Recalling the norm equivalence formulae
$$\Vert A\Vert_2^2=\lambda_{\tiny max}(A^{\prime}A)\leq \mbox{\rm tr}(A^{\prime}A)=\Vert A\Vert_F^2\leq r~\Vert A\Vert_2^2$$ for any $(r\times r)$-matrix $A$, the above estimates are valid if we replace the $\LL_2$-norm by the Frobenius norm.

Most of the semigroup analysis and the contraction inequalities developed in this article are based on the logarithm norms instead of
the spectral abscissa given by the top (real part of the) eigenvalues. The reasons are twofolds: 

Firstly, as its name indicates,
 the logarithmic norm represents the logarithmic decays of semigroups w.r.t the $\LL_2$-norm (cf. (\ref{def-log-norm})). These norms facilitate 
 the stability analysis of exponential semigroups.  
 
 On the other hand, most of the matrix exponential estimates expressed in terms of spectral abscissas 
 involve numerical constants that depends on the norm of the
diagonalization matrix  and its inverse, but also on polynomial functions w.r.t. the time parameter.  When the matrix has an 
ill conditioned eigen-system these constants are generally too large to obtain an effective useful estimate.  
We refer the reader to the formulae (\ref{first-estimate-Jordan}) and (\ref{estimate-Schur}).

For a more thorough discussion on these norms and 
 their used in the stability analysis of homogeneous semigroups of the form $e^{tA}$ we refer to~\cite{Zeibur}.

We end this section with a couple of rather well-known estimates in matrix theory. For any $(r\times r)$-square matrices $(P,Q)$ by a direct application of Cauchy-Schwarz inequality we have
\begin{equation}\label{form-pq-ref}
\vert\mbox{\rm tr}(PQ)\vert\leq \Vert P\Vert_F~\Vert Q\Vert_F.
\end{equation}
For any (symmetric and) positive semi-definite  $(r\times r)$-square matrices $P$ and $Q$.
\begin{equation}\label{f30}
 \mbox{\rm tr}\left(P^2\right)\leq \left(\mbox{\rm tr}\left(P\right)\right)^2\leq r~\mbox{\rm tr}\left(P^2\right)
\quad\mbox{\rm
and}\quad
\lambda_{\tiny min}(P)~\mbox{\rm tr}\left(Q\right)\leq \mbox{\rm tr}\left(PQ\right)\leq \lambda_{\tiny max}(P)~\mbox{\rm tr}\left(Q\right).
\end{equation}
The above inequality is also valid when $Q$ is positive semi-definite and $P$ is symmetric. We check this claim
using an orthogonal diagonalization of $P$ and recalling that $Q$ remains positive semi-definite (thus with non negative 
diagonal entries).
When both matrices $P$ and $Q $ are  
negative semi-definite the r.h.s. inequality remains valid if we replace $\left(\lambda_{\tiny min}(P),\lambda_{\tiny max}(P)\right)$ by  $\left(\lambda_{\tiny max}(P),\lambda_{\tiny min}(P)\right)$.

\section{Description of the models}\label{desc-sec-intro}
\subsection{The Kalman-Bucy filter}
Consider a time homogeneous linear-Gaussian filtering model of the following form
\begin{equation}\label{lin-Gaussian-diffusion-filtering}
\left\{
\begin{array}{rcl}
dX_t&=&\left(A~X_t+a\right)~dt~+~R^{1/2}_{1}~dW_t\\
dY_t&=&\left(C~X_t+c\right)~dt~+~R^{1/2}_{2}~dV_{t}.
\end{array}
\right.
\end{equation}
In the above display, $(W_t,V_t)$ is an $(r_1+r_2)$-dimensional Brownian motion, $X_0$ is a $r_1$-valued
Gaussian random vector with mean and covariance matrix $(\EE(X_0),P_0)$ (independent of $(W_t,V_t)$), the symmetric matrices $R^{1/2}_{1}$ and $R^{1/2}_{2}$ are invertible, $A$ is a square  $(r_1\times r_1)$-matrix, $C$ is an  $(r_2\times r_1)$-matrix, $a$ is a given $r_1$-dimensional column vector
and $c$ is an $r_2$-dimensional column vector, and $Y_0=0$. 
We also let $\Fa_t=\sigma\left(Y_s,~s\leq t\right)$ be the filtration generated by the observation process.

It is well-known that
the conditional distribution $\eta_t$ of the signal state $X_t$ given $\Fa_t$ is a $r_1$-dimensional Gaussian  distribution with a 
a mean and covariance matrix 
$$
\widehat{X}_t:=\EE(X_t~|~\Fa_t)\quad\mbox{\rm and}\quad
P_t:=\EE\left(\left(X_t-\EE(X_t~|~\Fa_t)\right)\left(X_t-\EE(X_t)\right)^{\prime}\right)
$$ 
given by the Kalman-Bucy filter
\begin{eqnarray}
d\widehat{X}_t&=&\left(A~\widehat{X}_t+a\right)~dt+P_{t}~C^{\prime}R^{-1}_{2}~\left(dY_t-\left(C\widehat{X}_t+c\right)dt\right)\label{nonlinear-KB-mean}
\end{eqnarray}
and the Riccati equation
\begin{eqnarray}
\partial_tP_t=\ricc(P_t) .\label{nonlinear-KB-Riccati}
\end{eqnarray}
defined in terms of the quadratic drift function
$$
\ricc~:~Q\in\SS_{r_1}~\mapsto~
\ricc(Q)=AQ+QA^{\prime}-QSQ+R\in\SS_{r_1}
$$ 
with $R=R_1$ and $ S:=C^{\prime}R^{-1}_{2}C  $.
When the dimension of the state is too large, as in most ocean and atmosphere stochastic models,
the solving of the Riccati matrix  evolution equation is untractable. Besides the problem of storing high dimensional matrices, we often need to resort to spectral technique and change of vector basis to solve analytically the Riccati equation.
For high dimensional problems these spectral techniques cannot be applied and another level of approximation need to be added. The idea of the EnKF is to replace the covariance matrices by sample covariance matrices associated with a 
well chosen mean-field particle model. These probabilistic models are defined in more details in the next section.

\subsection{A nonlinear Kalman-Bucy diffusion}

We consider the conditional nonlinear McKean-Vlasov type diffusion process
\begin{equation}\label{Kalman-Bucy-filter-nonlinear-ref}
d\overline{X}_t=\left(A~\overline{X}_t+a\right)~dt~+~R^{1/2}_{1}~d\overline{W}_t+\Pa_{\eta_t}C^{\prime}R^{-1}_{2}~\left[dY_t-\left((C\overline{X}_t+c)dt+R^{1/2}_{2}~d\overline{V}_{t}\right)\right]
\end{equation}
where  $(\overline{W}_t,\overline{V}_t,\overline{X}_0)$ are independent copies of $(W_t,V_t,X_0)$ (thus independent of
 the signal and the observation path).  In the above displayed formula $\Pa_{\eta_t}$ stands for the covariance matrix
$$
\Pa_{\eta_t}=\eta_t\left[(e-\eta_t(e))(e-\eta_t(e))^{\prime}\right]
\quad\mbox{\rm with}\quad \eta_t:=\mbox{\rm Law}(\overline{X}_t~|~\Fa_t)\quad\mbox{\rm and}\quad
e(x):=x.
$$
We shall call this probabilistic model the Kalman-Bucy (nonlinear) diffusion process. 

In contrast to conventional nonlinear diffusions the interaction does not take place only on the drift part but also on the diffusion matrix functional. In addition the nonlinearity does not depend on the distribution of the random states $\pi_t=\mbox{\rm Law}(\overline{X}_t)$ but on their conditional distributions $\eta_t:=\mbox{\rm Law}(\overline{X}_t~|~\Fa_t)$. 

Section~\ref{KB-diffusions-sec} discusses in some details the mathematical foundations of this conditional
nonlinear diffusion process.

We will also check that the conditional expectations of the random states $\overline{X}_t$
and their conditional covariance matrices $\Pa_{\eta_t}$ w.r.t. $\Fa_t$ satisfy the Kalman-Bucy and the Riccati Equations
\eqref{nonlinear-KB-mean} and \eqref{nonlinear-KB-Riccati}, {\em even when the initial variable is not Gaussian}; that is we have that
 \begin{equation}\label{ce-intro}
\EE\left(\overline{X}_t~|~\Fa_t\right)=\widehat{X}_t\qquad\mbox{\rm and}\qquad
\Pa_{\eta_t}=P_t.
\end{equation} 

In other words the flow of matrices $\Pa_{\eta_t}$ only depends on the covariance matrix of the initial state $\overline{X}_0$. This property comes from the structure of the nonlinear process equation  which ensures that
the mean and the covariance matrices satisfy the Kalman-Bucy filter and the Riccati equation. 
This property simplifies the stability analysis of this process. 

Given $\Pa_{\eta_0}$ the Kalman-Bucy Diffusion 
\eqref{Kalman-Bucy-filter-nonlinear-ref} can be interpreted as a non homogeneous Ornstein-Uhlenbeck type diffusion with a conditional covariance matrix $P_t=\Pa_{\eta_t}$ that satisfies the Riccati Equation \eqref{nonlinear-KB-Riccati} starting from $P_0=\Pa_{\eta_0}$. In this interpretation, the nonlinearity of the process is encapsulated in the  Riccati Equation \eqref{nonlinear-KB-Riccati}.

A more detailed description of the nonlinear semigroup of \eqref{Kalman-Bucy-filter-nonlinear-ref} is provided
 in Section~\ref{stab-section-p} dedicated to the stability properties of the Kalman-Bucy diffusions (see for instance 
 Lemma~\ref{lem-well-founded}).
If in addition $\overline{X}_0\stackrel{\tiny law}{=}X_0$ (which is Gaussian) then given $\Fa_t$ the random states $\overline{X}_t$
of the nonlinear Diffusion \eqref{Kalman-Bucy-filter-nonlinear-ref} are $r_1$-valued Gaussian random variables with
mean $\widehat{X}_t$ and covariance matrix $P_t$. Notice that deterministic initial states $\overline{X}_0$ can also be seen as Gaussian will a null covariance matrix.

By \eqref{ce-intro} the stability properties of the Kalman-Bucy filter
resume to the stability of the conditional expectations of the Kalman-Bucy diffusion, the reverse is clearly not true.
These questions are developed in some details in Section~\ref{global-contraction-section} and Section~\ref{local-contraction}.

\subsection{The Ensemble Kalman-Bucy filter}
The Ensemble Kalman-Bucy filter coincides with the mean-field particle interpretation
of the nonlinear diffusion process \eqref{Kalman-Bucy-filter-nonlinear-ref}. To be more precise we let
 $(\overline{W}^i_t,\overline{V}^i_t,\xi^i_0)_{1\leq i\leq N}$  be 
 $N$ independent copies of $(\overline{W}_t,\overline{V}_t,\overline{X}_0)$. 
In this notation, the EnKF is given by the Mckean-Vlasov type interacting diffusion process
\begin{equation}\label{fv1-3}
\left\{
\begin{array}{rcl}
d\xi^i_t&=&\left(A~\xi^i_t+a\right)dt+R^{1/2}_{1}d\overline{W}_t^i+p_tC^{\prime}R^{-1}_{2}\left[dY_t-\left((C\xi^i_t+c)dt+R^{1/2}_{2}~d\overline{V}^i_{t}\right)\right]\\
i&=&1,\ldots,N
\end{array}
\right.
\end{equation}
with the rescaled particle covariance matrices
\begin{equation}\label{fv1-3-2}
p_t:=\left(1-\frac{1}{N}\right)^{-1}~\Pa_{\eta^{N}_t}=\frac{1}{N-1}\sum_{1\leq i\leq N}\left(\xi^i_t-m_t\right)\left(\xi^i_t-m_t\right)^{\prime}
\end{equation}
 and the empirical measures
$$
\eta^{N}_t:=\frac{1}{N}\sum_{1\leq i\leq N}\delta_{\xi^i_t}
\quad\mbox{\rm and the sample mean}\quad
 m_t:=\frac{1}{N}\sum_{1\leq i\leq N}\xi_t^i.
$$

We also consider the $N$-particle model $\zeta_t=\left(\zeta^i_t\right)_{1\leq i\leq N}$ defined as $\xi_t=\left(\xi^i_t\right)_{1\leq i\leq N}$
by replacing the sample covariance matrix  $p_t$ by the true covariance matrix $P_t$ (in particular we have $\xi_0=\zeta_0$).

We end this section with some comments on these particle/ensemble filtering processes.

When $C=0$  the EnKF reduce to $N$ independent copies of the Ornstein-Uhlenbeck diffusive signal.
 In the same vein, for a single particle the covariance matrix is null so that the EnKF reduce to a single independent copy of
 the signal. In the case $r_1=1$ we have
 \begin{equation}\label{ref-divergence-1particle}
 \EE\left(\Vert m_t-X_t\Vert^2\right)=2~\mbox{\rm Var}(X_t)\quad\mbox{\rm and}\quad
  \EE\left(\Vert m_t-\widehat{X}_t\Vert^2\right)= \EE\left(\Vert m_t-X_t\Vert^2\right)+P_t
 \end{equation}
 
 In these rather elementary situations, {\em the stability property of the signal drift matrix $A$ is crucial to
 design some useful uniform estimates w.r.t. the time parameter.} The stability
 of the signal is a necessary condition to derive uniform estimates for any type of particle filters~\cite{dm-springer2004,dm-crc2013}
 w.r.t. the time parameter. 
 
 On the other hand by the rank nullity theorem, when $N<r_1$ the sample covariance matrix  $p_t$ is the sample mean of $N$ matrices
 of unit rank so that it has null eigenvalues. As a result, in some principal directions the EnKF is only 
 driven  by the signal diffusion. For unstable drift matrices the EnKF experiences divergence as it is not corrected by the innovation 
 process.
 
It should be clear from the above discussion that the stability of the signal is a necessary and sufficient condition to 
design useful uniform estimates w.r.t. the time horizon.

 We shall return to this question in section~\ref{some-uniform-moments-sec}.

We also recall that evolution of the conditional distributions of a nonlinear signal given the observations satisfy 
a complex 
nonlinear and stochastic measure valued equation. In probability and statistics literature this equation is called the 
Kushner-Stratonovitch nonlinear filtering equation~\cite{kushner,strato}. For linear-Gaussian models the conditional distributions of the states are Gaussian distributions. The evolution of the conditional expectations and covariance error matrices resume to the Kalman-Bucy filter~\cite{bucy}.
When the signal process \eqref{lin-Gaussian-diffusion-filtering} is nonlinear or when its perturbations are 
non Gaussian its is tempting to design a mean-field approximation of a nonlinear diffusion defined as in 
\eqref{Kalman-Bucy-filter-nonlinear-ref} by replacing the linear drift  $x\mapsto Ax+a$ by a nonlinear one, say $x\mapsto A(x)$. 
Of course, under some appropriate regularity conditions 
the sample means  converge to the first moment of the nonlinear process at hand.

Unfortunately, it is well-known that
this nonlinear diffusion process cannot capture any of the statistics of the optimal filter.  

We easily check this assertion by showing that the nonlinear Fokker-Planck equation associated with 
  the nonlinear process differs from the Kushner-Stratonovitch optimal filter equation. 
  
Up to
  the best of our knowledge there does not exist a single result that quantifies the error between the conditional distributions nor some statistics  of the
  random states of the particle model with the ones of the optimal filter. 
 
 Besides these drawbacks and these open important questions, the EnKF is of current use in nonlinear settings. This type of stochastic model can be thought as an extended EnKF approximation scheme. 
 
 The first step to understand these probabilistic models is to ensure the
 convergence of the mean-field approximation to the distribution of the nonlinear diffusion.
 We also refer the reader to the pioneering article~\cite{legland} for a discussion on the convergence of these non optimal mean-field models
 in the discrete time case. The uniform convergence of these nonlinear EnKF is still an open important questions
 for continuous as well as for discrete time filtering problems. We plan to investigate these questions in a forthcoming article.

\section{Statement of the main results}\label{statement-sec-intro}

\subsection{A stochastic perturbation theorem}
Our first main result shows that the stochastic processes $\left(m_t,p_t\right)$ satisfy the same equation as
$\left(\widehat{X}_t,P_t\right)$  up to some local fluctuation
orthogonal martingales with angle brackets that only depends on the sample covariance matrix $p_t$.

\begin{theo}[Perturbation theorem]\label{TH1}

The stochastic processes $\left(m_t,p_t\right)$ defined in \eqref{fv1-3-2} satisfy the diffusion equations
\begin{equation}\label{f21}
dm_t=\left(A~m_t+a\right)dt+p_t~C^{\prime}R^{-1}_{2}~\left(dY_t-\left(Cm_t+c\right)~dt\right)+\frac{1}{\sqrt{N}}~d\overline{M}_t
\end{equation}
with the vector-valued martingale $\overline{M}_t=\left(\overline{M}_t(k)\right)_{1\leq k\leq r_1}$ and the  angle-brackets
\begin{equation}\label{f22}
\displaystyle
\partial_t\langle \overline{M}_t(k),\overline{M}_t(k^{\prime})\rangle_t=R(k,k^{\prime})+
\left(p_tSp_t\right)(k,k^{\prime}).
\end{equation}
We also have the matrix-valued diffusion
\begin{equation}\label{f23}
dp_t=\left(Ap_t+p_tA^{\prime}-p_tSp_t+R\right)~dt+\frac{1}{\sqrt{N-1}}~dM_t
\end{equation}
with a symmetric matrix-valued martingale $M_t=\left(M_t(k,l)\right)_{1\leq k,l\leq r_1}$
and the angle brackets 
\begin{equation}\label{f24}
\begin{array}{rcl}
\displaystyle\partial_t\left\langle  M(k,l), M(k^{\prime},l^{\prime})\right\rangle_t&
=&\left(R+p_tSp_t\right)(k,k^{\prime})~ p_t(l,l^{\prime})+
\left(R+p_tSp_t\right)(l,l^{\prime})~ p_t(k,k^{\prime})\\
&&\\
&&\displaystyle+
\left(R+p_tSp_t\right)(l^{\prime},k) ~p_t(k^{\prime},l)+
\left(R+p_tSp_t\right)(l,k^{\prime}) ~p_t(k,l^{\prime}).
\end{array}
\end{equation}
In addition we have the orthogonality property
$$
\forall 1\leq k,l,l^{\prime}\leq r_1\qquad  
\left\langle  M(k,l), \overline{M}(l^{\prime})\right\rangle_t=0.
$$

\end{theo}

This fluctuation type theorem shows that the EnKF and the interacting sample covariance matrices satisfy a 
Kalman-Bucy recursion and a stochastic type Riccati equation.  
The extra level of randomness comes from the fact the particles mimic the random evolution of the nonlinear
McKean-Vlasov type diffusion (\ref{Kalman-Bucy-filter-nonlinear-ref}). The amplitude of these fluctuation-martingales as of order $1/\sqrt{N}$,
as any mean field type particle process. 

In contrast with the optimal Kalman-Bucy filter, the fluctuations of the covariance  matrices may corrupt the natural stabilizing effects
of the innovation process defined by the centered observations increments $\left(dY_t-\left(Cm_t+c\right)dt\right)$.

For partially observed signals we cannot expect any stability properties of Kalman-Bucy filter and the EnKF without introducing some structural conditions of observability and controllability on the signal-observation equation (\ref{lin-Gaussian-diffusion-filtering}). When $a=0=c$, observe that the Kalman-Bucy equation (\ref{nonlinear-KB-mean}) implies that
 \begin{equation}\label{Error-vector-def}
d(\widehat{X}_t-X_t)=(A-P_tS)~(\widehat{X}_t-X_t)~dt+P_{t}~C^{\prime}R^{-1/2}_2~dV_t+R^{1/2}_{1}~dW_t
 \end{equation}
 The EnKF evolution equation (\ref{f21}) stated in theorem~\ref{TH1} shows that the evolution of 
the error-vector $(m_t-X_t)$ has the same form as above by replacing $P_t$ by $p_t$, up to some fluctuation
martingale due to the interacting sample covariance matrices $p_t$ and the internal perturbations of the particles.

 This equation shows that the stability properties of these processes depends on the nature of the eigenvalues of the 
 matrices $A-P_tS$ and $A-p_tS$.
 
 \subsection{Stability of Kalman-Bucy nonlinear diffusions}\label{stab-nl-kb-dffusions}
 
 We further assume that $(A,R_1^{1/2})$ is a  controllable pair and $(A,C)$ is observable, that is  the matrices
 \begin{equation}\label{def-contr-obs}
\left[R_1^{1/2},A(R_1^{1/2})\ldots, A^{r_1-1}R_1^{1/2}\right]\quad
\mbox{\rm and}\quad
\left[\begin{array}{c}
C\\
CA\\
\vdots\\
CA^{r_1-1}
\end{array}
\right]
 \end{equation}
have rank $r_1$. By Bucy's theorem (cf. theorem~\ref{theo-Bucy} on p.\pageref{theo-Bucy}), under these conditions $P_t$  converge exponentially fast to $P$ as $t\uparrow\infty$. This unique fixed point is called the steady state error covariance matrix and it
satisfies the algebraic Riccati equation  
\begin{equation}\label{steady-state-eq}
\ricc(P)=AP+PA^{\prime}-PSP+R=0.
\end{equation}
In addition, the matrix difference
$A-PS$ is asymptotically stable even when the signal drift matrix $A$ is unstable. 
Under the above condition, there exists some parameters $\upsilon,\varpi_{\pm}>0$ such that
\begin{equation}\label{steady-state-eq-2}
\varpi_-~Id\leq
\int_{0}^{\upsilon}~e^{As}Re^{A^{\prime}s}~ds
\leq \varpi_+~Id
\quad\mbox{\rm
and}
\quad
\varpi_-~Id\leq 
\int_{0}^{\upsilon}~e^{-A^{\prime}s}Se^{-As}~ds \leq \varpi_+~Id
\end{equation}
The parameter $\upsilon$ is often called the interval of observability-controllability.

For a more detailed discussion on these fixed point Riccati equations (a.k.a. algebraic Riccati equation) and the connexions with optimal control theory, we refer the reader to
~\cite{abou,bucy2,poubelle,wonham,zelikin} and the references therein.

In reference to signal processing and control theory literature we adopt the following terminology.
\begin{defi}
 The Kalman-Bucy filter $\widecheck{X}_t$ associated with the initial covariance matrix $\widecheck{P}_0=P$ is called the steady state Kalman-Bucy filter.  The Kalman-Bucy diffusion $\breve{X}_t$ starting from an initial random state  $\breve{X}_0$ with covariance matrix $\widecheck{P}_0=P$
is called the steady state Kalman-Bucy diffusion.
\end{defi}

 It is important to observe that this stationary type
filter doesn't require to solve the Riccati equation (\ref{nonlinear-KB-Riccati}). 
For a more thorough discussion on the stability properties of Kalman-Bucy filters and Riccati equations we refer the reader to~\cite{abou,anderson,lancaster,poubelle,sontag,wonham,zelikin}.  We also refer the reader to section~\ref{local-contraction}
dedicated to contraction estimates of Riccati semigroups and related fundamental matrices.

Our main objective is to extend these stability properties at the level of the Kalman-Bucy nonlinear diffusion.

The conditional distributions $\eta_t$ of the random states $\overline{X}_t$ given $\Fa_t$ is a stochastic measure valued
process. In contrast to more conventional nonlinear Markov diffusions, the nonlinearity depends on the $\Fa_t$-conditional covariance
matrices. Thus, the computation of the distribution $\pi_t=\mbox{\rm Law}(\overline{X}_t)$ of the random states
$\overline{X}_t$ requires to compute the conditional covariance matrices $\Pa_{\eta_t}$.

To study the stability properties of the flow $\pi_t$ it is natural to introduce a copy $\breve{X}_t$ of $\overline{X}_t$
coupled to $\overline{X}_t$  by choosing the same observation process, and by only changing the random perturbations $(\overline{W}_t,\overline{V}_t)$. In other words $(\overline{X}_t,\breve{X}_t)$ are $\Fa_t$-conditionnally independent.

To describe our main
results with some precision we introduce some terminology.

\begin{defi}
Let $(P_t,\widecheck{P}_t)$ be a couple of solutions of the Riccati Equation \eqref{nonlinear-KB-Riccati}
starting at two possibly different values $(P_0,\widecheck{P}_0)$. We also denote by  $(\overline{X}_t,\breve{X}_t)$ be a couple of Kalman-Bucy Diffusions \eqref{Kalman-Bucy-filter-nonlinear-ref} starting from two random states
with covariances matrices $(P_0,\widecheck{P}_0)$. We denote by $\left(\pi_t,\breve{\pi}_t\right)$ and $\left(\eta_t,\breve{\eta}_t\right)$ the distributions and 
 the $\Fa_t$-conditional distributions of
 $\left(\overline{X}_t,\breve{X}_t\right)$, and we set
 
  \begin{equation}\label{f-easy-intro-prop-mean}
\widehat{X}_t:= \EE\left(\overline{X}_t~|~\Fa_t\right)\quad\mbox{and}\quad \widecheck{X}_t= \EE\left(\breve{X}_t~|~\Fa_t\right).
\end{equation}

\end{defi}

Equivalently,
 $(\widehat{X}_t,\widecheck{X}_t)$ are solutions of the Kalman-Bucy Equation \eqref{nonlinear-KB-mean}
running with the covariance matrices $(P_t,\widecheck{P}_t)$.
Since $Y_0=0$, the signal state is not observed at the origin. This clearly implies that
$$
\left(\pi_0,\breve{\pi}_0\right)=\left(\eta_0,\breve{\eta}_0\right)
\quad\mbox{and}\quad
\left(\Pa_{\pi_0},P_{\check{\pi}_0}\right)=\left(\Pa_{\eta_0},\Pa_{\check{\eta}_0}\right)=(P_0,\widehat{P}_0).
$$

Using \eqref{f-easy-intro-prop-mean} it should be clear that most of the stability properties of the Kalman-Bucy diffusions (expressed in terms of some convex criteria) can be used to deduce the ones of their conditional averages
but the reverse is clearly not true.
For instance using \eqref{f-easy-intro-prop-mean}, for any $n\geq 1$ and $t\geq 0$ we readily check that
$$
\WW_{n}\left(
\mbox{\rm Law}(\widehat{X}_t),\mbox{\rm Law}(\widecheck{X}_t)
\right)\leq \WW_{n}\left(
\pi_t,\breve{\pi}_t
\right).
$$

In this context, one of our main result can basically be stated as follows.
\begin{theo}\label{prop-stab-Kalman}
Assume the existence of a positive semi-definite fixed point $P$ of \eqref{steady-state-eq} s.t.
 $\mu(A-PS)<0$.  For any $n\geq 1$,  $\epsilon\in ]0,1/2]$ and $t\geq 0$ we have
 \begin{eqnarray}
\WW_{2n}\left(
\pi_t,\breve{\pi}_t
\right)&\leq& c(\epsilon)
\exp{\left[(1-\epsilon)~\mu(A-PS)~t\right]}~\left(\WW_{2n}\left(
\pi_0,\breve{\pi}_0
\right)+\Vert P_0-\widecheck{P}_0 \Vert_F\right)
\label{cv-steady-df}
\end{eqnarray}
for some constant $c(\epsilon)<\infty$ whose values depend on the parameter $\epsilon$ and on the initial matrices $(P_0,\widecheck{P}_0)$.
 \end{theo}
The proof of this theorem is provided in section~\ref{stab-kalman-sec}.

In the further development of this section we let $\left(\eta_t,\breve{\eta}_t\right)$
 be the $\Fa_t$-conditional probability distributions of a couple of Kalman-Bucy diffusions $\left(\overline{X}_t,\breve{X}_t\right)$ starting from
two Gaussian random variables $(X_0,\breve{X}_0)$ with 
 covariance matrices $\left(P_0,\widecheck{P}_0\right)$.   Also let $(\widehat{X}_t, \widecheck{X}_t)$ be the Kalman-Bucy filters associated with  a couple of Riccati equations starting at $\left(P_0,\widecheck{P}_0\right)$.
\begin{theo}\label{theo-Ent}
 Assume that the algebraic Riccati Equation
 \eqref{steady-state-eq} has a positive definite fixed point $P$ and  $\mu(A-PS)<0$.
 In this situation,  there exists some $t_0$ and some finite constant $c$  that depends on $(P_0,P)$ such that
 for any $t\geq t_0$ and any $\epsilon\in ]0,2]$ we have the quenched almost estimate  
\begin{equation}\label{last-formula}
 \mbox{\rm Ent}\left(\eta_t~|~\breve{\eta}_t\right)\leq  c(\epsilon)~\left[\exp{\left((1-\epsilon)\mu(A-PS)t\right)}~\Vert P_0-\widecheck{P}_0 \Vert_F+
 \Vert \widecheck{X}_t-\widehat{X}_t\Vert^2\right]
\end{equation}
In addition, for any $t\geq t_0$ and any $\epsilon\in ]0,1/2]$ we have the annealed estimate 
\begin{equation}\label{last-formula-bis}
 \mbox{\rm Ent}\left(\pi_t~|~\breve{\pi}_t\right)\leq  c(\epsilon)~\exp{\left((1-\epsilon)\mu(A-PS)t\right)}~
\left(\WW_{2}^2\left(
\pi_0,\breve{\pi}_0
\right)+\Vert P_0-\widecheck{P}_0 \Vert_F\right)
\end{equation}
In the above display, $c(\epsilon)$ stands for some finite constant  whose values depend on the parameter $\epsilon$ and on the initial matrices $(P_0,\widecheck{P}_0)$
\end{theo}

We easily check \eqref{last-formula-bis} using \eqref{cv-steady-df}, \eqref{last-formula} and the convexity property of the relative entropy 
$$
\mbox{\rm Ent}\left(\nu K_1~|~\nu K_2\right)\leq 
\int~\nu(dy)~\mbox{\rm Ent}\left(K_1(y,\point)~|~K_2(y,\point)\right)
$$
which is valid for any Markov transitions $K_i(y,dx)$ from $\RR^{r_2}$ into $\RR^{r_1}$ and any 
probability distribution $\nu$ on $\RR^{r_2}$. In the above display $(\nu K_i)$ stands for the probability measures
on $\RR^{r_1}$ defined by the transport formula
$$
(\nu K_i)(dx):=\int~\nu(dy) K_i(y,dx).
$$

The proof of the almost sure estimate \eqref{last-formula} is provided in section~\ref{sect-last-entropy}.

\subsection{An uniform propagation of chaos theorem}

The EnKF avoid the numerical solving of the Riccati equation (\ref{fv1-3}) or the use of the steady state $P$ 
by using interacting sample covariance matrices. Even when $(\xi^i_0)_{1\leq i\leq N}$ are independent copies of a Gaussian random
variable $\overline{X}_0$ with the steady covariance matrix, the initial sample covariance matrices $p_t$ defined in (\ref{fv1-3-2}) fluctuates
around the limiting value $P$. 

These random fluctuations of the sample covariance matrices $p_t$ eventually corrupt the stability in the EnKF, even if the filtering problem is observable and controllable in the conventional sense.
 For instance, in practical situations the empirical covariance matrices may not be invertible for small sample sizes. This simple observation shows that the Lyapunov theory based on inverse of covariance matrices developed in~\cite{anderson,bensoussan,luenberger} cannot be applied in the context of EnKF.
  
In practice, it has also been observed that fluctuations of the sample-covariances of the EnKF eventually mislead the natural stabilizing effect of the observation process in the Kalman-Bucy filter evolution. 
These fluctuations induce an underestimation of the true error covariances. As a result the EnKF 
   ignores the important information delivered by the sensors. This lack of observation also leads to the divergence of the filter.
 
 Last but not least, from the numerical viewpoint the Kalman-Bucy filter and the EnKF are also know to be non robust, in the sense that arithmetic errors may accumulate even if the exact filter is stable.

All of  these instability properties of the EnKF are well known. They are often referred as the catastrophic filter divergence in data assimilation literature. For a more thorough discussion on these issues we refer the eader to the articles~\cite{gottwald-majda,kelly-majda-tong,majda-harlim}, and the references therein. As mentioned by the authors in~\cite{kelly-majda-tong}, "catastrophic filter divergence is a well-documented but mechanistically mysterious phenomenon whereby ensemble-state estimates explode to machine infinity despite the true state remaining in a bounded region".  The second main result of this article is to provide uniform propagations of chaos properties w.r.t. the time horizon.
To this end, we assume that the following observability condition is satisfied:
\begin{equation}\label{new-condition}
\hskip-5cm(\mbox{\rm S})\hskip3cm S=\rho(S)~Id\quad\mbox{\rm for some}\quad \rho(S)>0.
\end{equation}
Under this condition, we have
 \begin{equation}\label{triangle-ref}
\mu(A)-\rho(S)~\lambda_{\tiny max}(P)\leq 
\mu(A-PS)\leq \mu(A)-\rho(S)\lambda_{\tiny min}(P)\leq \mu(A)
 \end{equation}
 This shows that $\mu(A-PS)<0$ as soon as $\mu(A)<0$.
 When $\mu(A)>0$, it is also met as soon as $\mu(A)\leq \rho(S)\lambda_{\tiny min}(P)$.

Several examples of sensor models satisfying this condition are discussed in section~\ref{observability-section}.

\begin{theo}\label{theo-intro}
Assume that the observability condition $(\mbox{\rm S})$ is met and $\mu\left(A\right)<0$. In this situation, for any $n\geq 1$ and any sufficiently large $N$ we have the uniform estimates
\begin{equation}\label{upxi}
\sup_{t\geq 0}{\EE\left[\Vert p_t-P_t\Vert_{F}^{n}\right]^{\frac{1}{n}}}\leq c(n)/\sqrt{N}\quad
\mbox{and}\quad
\sup_{t\geq 0}{\EE\left[\Vert \xi_t^1-\zeta^1_t\Vert^{n}\right]^{\frac{1}{n}}}\leq c^{\prime}(n)/\sqrt{N}.
\end{equation}
 \end{theo}

When condition $(\mbox{\rm S})$ is not necessarily met but the signal-drift matrix $A$ is stable we have
 the uniform estimate
 $$
 \mu\left(A\right)<0\Longrightarrow\forall 1\leq n<1+\frac{N-1}{2r_1}~\frac{\lambda_{\tiny min}(S)}{\lambda_{\tiny max}(S)}\qquad
 \sup_{t\geq 0}
\EE\left(\left[\mbox{\rm tr}(p_t)\right]^n\right)\leq c(n)
 $$
(cf. proposition~\ref{prop-1}). The proof of the estimates (\ref{upxi})  relies on the observability condition (S).
This condition is satisfied for any filtering problem with $r_1=r_2$ with an invertible matrix $C$, up to a change of basis;
see for instance (\ref{similar-filtering}). 

We conjecture that condition (S) is purely technical and can be relaxed for stable signal-drift matrices. 
For unstable drift matrices the form of the matrix $S$ may also corrupt the regularity of the sample
covariance matrices. We already know that for $r_1>N$ the matrices $p_t$ have necessarily at least one null eigenvalue,
so that the number of null eigenvalues of $p_tS$ is more likely to be larger when $S$ is not of full rank.  

The following corollary is a consequence of  \eqref{upxi}.

\begin{cor}\label{cor-intro-ref}
Under the assumptions stated in theorem~\ref{theo-intro}, for any $n\geq 1$ we have the uniform estimates
\begin{equation}\label{upxi-Sharp-bis}
 \sup_{t\geq 0}{\WW_n\left(\mbox{\rm Law}( \xi_t^1),\eta_t\right)}\leq c(n)/\sqrt{N}
 \quad
\mbox{and}\quad
\EE\left(\left\vert\eta^N_t(f)-\eta_t(f)\right\vert^n\right)^{\frac{1}{n}}\leq c^{\prime}(n)/\sqrt{N}
\end{equation}
for any $1$-Lipschitz function $f$ on $\RR^{r_1}$. In particular, this implies that
\begin{equation}\label{upxi-Sharp}
\sup_{t\geq 0}{\EE\left[\Vert m_t-\widehat{X}_t\Vert^{n}\right]^{\frac{1}{n}}}\leq c(n)/\sqrt{N}.
\end{equation}
Inversely the uniform estimates \eqref{upxi-Sharp} implies that $\mu(A)<0$ as soon as 
$R\not=0=C$. 

\end{cor}

The detailed proof of Theorem~\ref{TH1} is postponed to Section~\ref{sec-EnFK-ss1}. The proof of Theorem~\ref{theo-intro}
is based on uniform moments estimates developed in Section~\ref{unif-moment-EnKF-ss2}. The detailed proof of 
the uniform estimates \eqref{upxi} is presented in full details in Section~\ref{quantitative-section}. 

It is important to observe that all the $\LL_n$-mean error estimates between the sample covariance $p_0$ and $P_0$, as well as the ones between the sample mean $m_0$ and $\widehat{X}_0$, 
are immediate for $t=0$ (as a direct application of the law of large numbers for independent random sequences).

A discussion on the stability condition $\mu\left(A\right)<0$ is provided in Section~\ref{preliminary-sec}. In Section~\ref{some-uniform-moments-sec} we will show that this condition cannot be relaxed to derive uniform
estimates \eqref{upxi-Sharp} as soon as $R\not=0$ and $C=0$. The strong observability condition (S) is discussed in some details in section~\ref{section-regularity}. Section~\ref{sous-section-section-regularity} is dedicated to the stability of stochastic observer processes defined as  the Kalman-Bucy diffusion  (\ref{nonlinear-KB-mean}) 
by replacing $\Pa_{\eta_t}$ by a flow of stochastic covariance matrices.
Several illustrations are presented in section~\ref{observability-section}.

We also emphasize that the conditions (S) and $\mu(A)<0$ are only used to derive uniform estimates
w.r.t. the time horizon. Without these conditions  the statements of theorem~\ref{theo-intro} and corollary~\ref{cor-intro-ref}
remains valid without the supremum operations w.r.t. the time parameter. In this general situation, the constants $c(n)$ and
$c^{\prime}(n)$ are replaced by some constants $c_t(n)$ and
$c^{\prime}_t(n)$ whose values depend on the time horizon.

\section{Some comments on our regularity conditions}\label{section-regularity}

The stability analysis of diffusion processes is always much more documented than the ones on their possible divergence.
For instance, in contrast with conventional Kalman-Bucy filters, the stability properties of the EnKF are not
induced by some kind of observability or controllability condition. The only known results for discrete generation
EnKF is the recent work by X. T. Tong, A. J. Majda and D. Kelly~\cite{tong}. One of the main assumptions of the article is that the sensor-matrix 
is of full rank. The authors also provide a concrete numerical example of filtering problem with sparse observations for which the EnKF experiences a catastrophic divergence.

The stability of linear dynamical systems with time varying drift-matrices $A_t$ is much complex than
the one of time homogeneous models. In our context, the drift-matrix are also random since they 
encapsulates the fluctuations of the stochastic covariance matrices. As shown in~\cite{ilchmann} the fact that  
the real part of the spectrum of $A_t$ is negative  is neither necessary nor sufficient for the exponential stability of the system.
 It may even happen that
the semigroup of the system is unstable even if the real part of the eigenvalues of $A_t$ remain negative for all times. Inversely the 
system may be stable even when one of the eigenvalue is positive for all times.

The analysis of the exponential stability of time varying linear system requires to estimate the variations of 
the matrices $A_t$ w.r.t. the time parameter.
  All type of sufficient conditions
stated in~\cite{brocket,coppel,ilchmann,krause,kreisselmeier,rosenbrock} are restricted to deterministic systems, and thus
cannot
be applied to interacting diffusion processes. They are
  also based on continuity/tightness type properties, as well as on piecewise smoothness properties or on
the uniform boundedness
of the velocity field $\partial_tA_t$. These  regularity properties doesn't hold for stochastic diffusions. 
In our context, the fluctuation matrices are given by
\begin{equation}\label{true-perturbations}
Q_t:=\sqrt{N}~(p_t-P_t)\in\SS_{r_1}\Longleftrightarrow p_t=P_t+\frac{1}{\sqrt{N}}~Q_t\in\SS_{r_1}^+
\end{equation}

Our next objective is to initiate a more refined mathematical analysis to understand 
 these divergence properties in terms of global divergence properties, locally ill-conditioned filtering models, and stochastic type
 observers.

\subsection{Stable and divergence regions}\label{sous-section-section-regularity}

To clarify the presentation, we further assume that
the Kalman-Bucy and the Riccati equation start at the steady state $P_0=P$. 

  Let
  $Z_t$ be some observer type process defined as  (\ref{nonlinear-KB-mean}) 
by replacing $\Pa_{\eta_t}$ by some covariance matrix of the form
$
(P+Q_t)\in\SS_{r_1}^+
$, for some flow of symmetric matrices $Q_t\in\SS_{r_1}$ mimicking the fluctuation matrices (\ref{true-perturbations}). 
In other words, the matrices $Q_t$ reflects the fluctuations of the empirical covariance matrices $p_t$ around the steady state $P$.

In this case we have 
\begin{eqnarray}
dZ_t
&=&\left[A~Z_t-(P+Q_t)S (Z_t-X_t)\right]~dt+(P+Q_t)~C^{\prime}R^{-1/2}_{2} dV_t\label{def-observer}
\end{eqnarray}
When $Q_t=0$, the null matrix, the process $Z_t$ resumes to the steady state Kalman-Bucy filter.
We set
$$
A-(P+Q_t)S=\overline{A}-Q_tS\quad\mbox{\rm with}\quad \overline{A}:=A-PS
\quad\mbox{\rm and}\quad
E_t:=\exp{\left[\oint_0^t (\overline{A}-Q_sS) ds\right]}.
$$ 
The vector error process $\Za_t:=Z_t-X_t$ is given by the stochastic Ornstein-Ulhenbeck
process
$$
d\Za_t=(\overline{A}- Q_tS)~\Za_t~dt~+~\sqrt{2}~d\Wa_t\Rightarrow \displaystyle\Za_t=E_t~\Za_0+\sqrt{2}~\int_0^t~E_tE_s^{-1}~d\Wa_s
$$
In the above display $\Wa_t$ stands for diffusion 
$$
\Wa_t:=~(\overline{W}_t-W_t)/\sqrt{2}+(P+Q_t)C^{\prime}R^{-1/2}_{2}~(V_t-\overline{V}_t)/\sqrt{2}\stackrel{\tiny law}{=}
 W_t+(P+Q_t)C^{\prime}R^{-1/2}_{2}V_t
$$
with covariance matrices
$
\Ra_t=I+(P+Q_t)S(P+Q_t)
$.
When the flow of matrices $Q_t$ enter into the set
\begin{equation}\label{global-div-region-def}
 \Qa_{\mbox{\tiny div}}=\left\{Q\in \SS_{r_1}~:~\varsigma(\overline{A}-QS)>0
\right\}
\end{equation}
the observer experiences a divergence in at least one of the principal directions. Notice that
\begin{equation}\label{key-sg-estimate}
\Vert E_t\Vert_2\leq \exp{\left[\int_0^t\mu(\overline{A}-Q_sS)~ds\right]}
\end{equation}
This semigroup estimate allows to quantify the stability of the process $\Za_t$ as soon as $\mu(\overline{A}-Q_sS)<-\delta$
for some $\delta>0$ for sufficiently large time horizons. 

One natural strategy is to analyze the contraction properties
of the stochastic flow $E_t$ generated by the stochastic matrices $\overline{A}-Q_tS$ and their logarithmic norms $\mu(\overline{A}-Q_tS)$.
More precisely, under the strong observability condition (S) stated in (\ref{new-condition}) we have
\begin{equation}\label{our-conditions}
 \mu(\overline{A}-Q_sS)\leq \mu(A)+\rho(S)~\mu(-(P+Q_s))\leq \mu(A)<0
\end{equation}
as soon as $\mu(A)<0$, for {\em any} possible symmetric fluctuations $Q_s$ s.t. $P+Q_s\geq 0$.

This shows that for stable signal-drift matrices $\mu(A)<0$ the condition (S) ensures that  the stochastic observer is both theoretically and numerically stable for any type of fluctuations $Q_s$. The same reasoning will be used to show that the stability of the signal is transferred to the EnKF filter. 

Without condition $(S)$ is easy to work out several examples of $2$-dimensional filtering problems with a stable-drift matrix $\mu(A)<0$
and such that  $\mu(\overline{A}-Q_sS)>0>\lambda_{\tiny min}((\overline{A}-Q_sS)_{\tiny sym})$ for some flow of symmetric matrices $Q_s$ s.t. $P+Q_s\geq 0$. In this context, even if the EnKF is numerically stable it is difficult to analyze theoretically this class of locally ill conditioned
models using spectral and semigroup techniques.

  In the reverse angle, in practical situations the EnKF generally experiences
 severe divergence when $\mu(A)>0>\mu(\overline{A})$. In this situation, we already know from (\ref{ref-divergence-1particle}) that we cannot expect to
 have uniform propagation of chaos estimates for any fluctuation matrices.  Also observe that
 $$
\left(\mbox{\rm (S)}\quad \mbox{\rm and}\quad
\mu(\overline{A})<0\right)\Longrightarrow \mu(\overline{A}-Q_sS)<\mu(\overline{A})+\rho(S)~\mu(-Q_s)<\mu(\overline{A})<0
$$ for any {\em positive} semidefinite 
fluctuations $Q_s\geq 0$ around the 
steady state $P$. Unfortunately, we cannot ensure that the fluctuations of the sample covariance matrices are always
 positive.

 As mentioned above, the pivotal semigroup estimate (\ref{key-sg-estimate}) requires to estimate the logarithmic norm of the 
 stochastic flow of matrices $\overline{A}-Q_tS$.    Several technical difficulties arise:

  The first one comes from the fact that
 $\varsigma(\overline{A})<0$ doesn't implies that $\mu(\overline{A})<0$, since $\overline{A}\not=\overline{A}^{\prime}$. 
 When $\mu(\overline{A})>0$ the  steady state  Kalman-Bucy is locally ill-conditioned,
 in the sense that the worst
 fluctuation around the true signal behave like $\exp{\left[\mu(\overline{A})\Delta t\right]}$ is a short transient time $\Delta t$.
 This local divergence property may occur even when all the eigenvalues
 of $\overline{A}$ or even the ones of $A$ are negative. This indicates that it is hopeless to analyze the stability of the Kalman-Bucy filters estimates based on the semigroup inequality (\ref{key-sg-estimate}) for such ill-conditioned systems. We are faced to the same issues if we try to quantify the propagations of the 
fluctuation $Q_t$ in the system.  Some illustrations are discussed in the appendix, on page~\pageref{some-illustrations-2d}

This discussion indicates that the stability property  $\mu(A)<0$ and the observability condition (\ref{new-condition}) seem to be essential to control the fluctuations of the EnKF
sample covariance matrices {\em for any number of samples}. These conditions also ensure the semigroup contraction properties needed to derive uniform 
$\LL_n$-mean error estimates of the EnKF particle filter.

\subsection{Full observation sensors}\label{observability-section}

When the observation variables are the same as the ones of the signal; the signal observation has the same dimension as the signal and resumes to
 some equation of the form
\begin{equation}\label{fully-observed-sensor}
dY_t=b~ X_t~dt+\sigma_2~dV_t
\end{equation}
for some parameters $b\in\RR$ and $\sigma_2> 0$.  These sensors are used in data grid-type assimilation
problems when measurements can be evaluated at each cell. These fully observed models are discussed in section 4 in~\cite{harlim-hunt} in the context of the Lorentz-96 filtering problems.  These observation processes are also used to the article~\cite{berry-harlim} for application to nonlinear and multi-scale filtering problem. In this context, the observed variables represents the slow components of the signal. When the fast components are represented by a some Brownian motion with a prescribed covariance matrix, the filtering of the slow components with full observations take the form (\ref{fully-observed-sensor}).

The sensor model discussed in (\ref{fully-observed-sensor}) clearly satisfies
condition (\ref{new-condition})  with the parameter $\rho(S)=(b/\sigma_2)^2$.   
This rather strong condition (\ref{new-condition}) ensures that these fluctuations doesn't propagate w.r.t. the time parameter, regardless of the initial data. Under this condition we shall prove that the particle EnKF has uniformly bounded  $\LL_n$-moments for any $n\geq 1$ (see for proposition~\ref{prop-1}).

We emphasize  that the observability condition (\ref{new-condition})
 is satisfied when the filtering problem are similar to a fully observed sensor model; that is, up to a change of basis functions.
More precisely, any filtering problem (\ref{lin-Gaussian-diffusion-filtering}) with $r_1=r_2$ and s.t.  $\Ca:=(R_2^{-1/2}C)$ is invertible can be turned into a filtering problem equipped with an identity sensor matrix; even when the original matrix $S=C^{\prime}R_2^{-1}C=\Ca^{\prime}\Ca$ doesn't satisfies (\ref{new-condition}). 
To check this claim we observe that
\begin{equation}\label{similar-filtering}
\Ya_t:=R^{-1/2}_2Y_t\quad\mbox{\rm and}\quad \Xa_t:=\Ca X_t
\Longrightarrow \left\{\begin{array}{rcl}
d\Xa_t&=&\Aa~\Xa_t~dt+\Ra_1^{1/2}~dW_t\\
d\Ya_t&=&\Xa_t~dt+dV_t
\end{array}\right.
\end{equation}
with the signal drift matrix
$
\Aa:=\Ca A\Ca^{-1}$ and the diffusion covariance matrix
$\Ra_1:=\Ca R_1\Ca^{\prime}$.
In this situation the filtering model $(\Xa_t,\Ya_t)$ satisfies (\ref{new-condition})  with $\Sa=Id\Rightarrow\rho(\Sa)=1$.  The link between the logarithmic
norm of $\Aa$ and the original signal drift matrix $A$ is given by the formula
$$
\mu(\Aa)=\frac{1}{2}~\lambda_{\tiny max}\left(\Ca A\Ca^{-1}+\left(\Ca A\Ca^{-1}\right)^{\prime}\right)
$$
 For orthogonal matrices $\Ca$ we have $\mu(\Aa)=\mu(\Ca A\Ca^{\prime})=\mu(A)$. Otherwise, the condition $\mu(\Aa)<0$ depends on
 the triplet of matrices $(A,C,R_2)$ associated with the original filtering problem. For instance, when $r_1=r_2=2$, $R_2=Id$ and
 a symmetric negative definite drift matrix $A$, the condition $\mu(A)<0$ is equivalent to the fact that
 \begin{equation}\label{ref-A-def-negative}
 A_{1,1}<0\quad \mbox{\rm and}\quad A_{1,1}A_{2,2}>A_{1,2}^2
 \end{equation}
 For sensor matrices of the form $C=\left(
 \begin{array}{cc}
1&0\\
 0&\beta
 \end{array}
 \right)$, for some~$\beta\not=0$, condition $\mu(\Aa)<0$ takes the strongest  form
 $$
  A_{1,1}<0\quad \mbox{\rm and}\quad A_{1,1}A_{2,2}> \left(\frac{\beta+\beta^{-1}}{2}\right)^2~A_{1,2}^2
 $$

To better understand the importance of the matrix $S$ introduced in (\ref{nonlinear-KB-Riccati}) observe that
\begin{eqnarray}
d\widehat{X}_t
&=&\left(A~\widehat{X}_t+a\right)~dt+P_t~(d\overline{Y}_t-S\widehat{X}_tdt)\label{ref-observers}
\end{eqnarray}
with the $r_1$-dimensional observation process $\overline{Y}_t$ given by
$$
\overline{Y}_t:=C^{\prime}R^{-1}_2~Y_t\Longrightarrow
d\overline{Y}_t=S~X_t~dt+C^{\prime}R^{-1/2}_2~dV_t
$$

 Let us assume that $r_1=1<r_2=2$, and the $(2\times 1)$-sensor matrix is given by $$C=[1,\alpha]\quad \mbox{\rm for some} \quad \alpha\in \RR.
 $$ 
 For unit diffusion covariance matrices the partially observed filtering systems
 $(X_t,Y_t)$ and $(X_t,\overline{Y}_t)$  observable as soon as 
 \begin{equation}\label{observability-condition-alpha}
 A_{1,2}+\alpha A_{2,2}\not=\alpha~(A_{1,1}+\alpha A_{2,1}).
 \end{equation}
 Nevertheless, in this situation the diffusion matrix $S=\left[
 \begin{array}{cc}
 1&\alpha\\
 \alpha&\alpha^2
 \end{array}
 \right]$ of this new $2$-dimensional $\overline{Y}_t$ is no more invertible. This shows that these $1$-dimensional partial observations
models cannot be turned into regular $2$-dimensional sensors.
These $2$-dimensional filtering problems equipped with a $1$-dimensional sensor  are  one of the simplest 
examples of controllable and observable filtering problems that doesn't satisfy the observability condition (\ref{new-condition}) even if the signal drift is stable.

\section{Stability properties of Kalman-Bucy diffusions}\label{stab-section-p}

\subsection{Kalman-Bucy diffusions}\label{KB-diffusions-sec}

As noticed in the introduction, the Kalman-Bucy Diffusion \eqref{Kalman-Bucy-filter-nonlinear-ref} strongly differs from conventional nonlinear diffusion processes. The evolution of this new class of probabilistic models
depend on the $\Fa_t$-conditional distribution of the random states.
 
This section provides a more detailed discussion on this new class of nonlinear McKean-Vlasov type diffusions
with $\Fa_t$-conditional distribution interactions.

 \begin{defi}
Let $\varphi_{s,t}$ be the dynamical semigroup of the Riccati Equation
\eqref{nonlinear-KB-Riccati} given for any $s\leq t$ by
$$
\varphi_{s,t}(P_s)=P_t\quad\mbox{\rm and set}\quad \varphi_{0,t}=\varphi_{t}.
$$
\end{defi}

Next lemma shows that the Kalman-Bucy Diffusion \eqref{Kalman-Bucy-filter-nonlinear-ref} is well-posed.

\begin{lem}\label{lem-well-founded}
Consider the non homogeneous diffusion given by 
$$
d\overline{X}_t=\left(A~\overline{X}_t+a\right)~dt~+~R^{1/2}_{1}~d\overline{W}_t+\varphi_t(P_0)C^{\prime}R^{-1}_{2}~\left(dY_t-\left((C\overline{X}_t+c)dt+R^{1/2}_{2}~d\overline{V}_{t}\right)\right)
$$
where  $(\overline{W}_t,\overline{V}_t,\overline{X}_0)$ are independent copies of $(W_t,V_t,X_0)$ and
 $P_0=\Pa_{\eta_0}$. In this situation, we have
$$
\eta_t:=\mbox{\rm Law}(\overline{X}_t~|~\Fa_t)~\Longrightarrow~
\Pa_{\eta_t}=\varphi_t(P_0)\quad\mbox{and}\quad
 \EE\left(\overline{X}_t~|~\Fa_t\right)=\widehat{X}_t
$$
where $\widehat{X}_t$ stands for the solution of the Kalman-filter \eqref{nonlinear-KB-mean} driven by the solution
$P_t=\Pa_{\eta_t}$ of
Riccati
Equation \eqref{nonlinear-KB-Riccati} starting at $P_0=\Pa_{\eta_0}$. 
\end{lem}

\proof
By construction, we have
\begin{equation}\label{nonlinear-KB-mean-non-mauvais-label}
d\EE\left(\overline{X}_t~|~\Fa_t\right)=\left(A~\EE\left(\overline{X}_t~|~\Fa_t\right)+a\right)~dt+\varphi_t(P_0)~C^{\prime}R^{-1}_{2}~\left(dY_t-\left(C_t\EE\left(\overline{X}_t~|~\Fa_t\right)+c\right)dt\right).
\end{equation}
We set $\widetilde{X}_t:=\overline{X}_t-\EE\left(\overline{X}_t~|~\Fa_t\right)$. In this notation we have
$$
d\widetilde{X}_t=\left[A~-\varphi_t(P_0)C^{\prime}R^{-1}_{2}C\right]~\widetilde{X}_t~dt~+~R^{1/2}_{1}~d\overline{W}_t-
\varphi_t(P_0)~C^{\prime}R^{-1/2}_{2}d\overline{V}_{t}.
$$
This implies that
\begin{eqnarray*}
d\left(\widetilde{X}_t\widetilde{X}_t^{\prime}\right)&
=&\left[A-\varphi_t(P_0)C^{\prime}R^{-1}_{2}C\right]~\widetilde{X}_t\widetilde{X}_t^{\prime}~dt~+\left[~R^{1/2}_{1}~d\overline{W}_t-
\varphi_t(P_0)~C^{\prime}R^{-1/2}_{2}d\overline{V}_{t}\right]\widetilde{X}_t^{\prime}\\
&&+\widetilde{X}_t\widetilde{X}_t^{\prime}\left[A-\varphi_t(P_0)C^{\prime}R^{-1}_{2}C\right]^{\prime}~dt~+~\widetilde{X}_t~\left[R^{1/2}_{1}~d\overline{W}_t-
\varphi_t(P_0)~C^{\prime}R^{-1/2}_{2}d\overline{V}_{t}\right]^{\prime}\\
&&\hskip6.5cm+\left[R_{1}+\varphi_t(P_0)~C^{\prime}R^{-1}_{2}C\varphi_t(P_0)\right]~dt.
\end{eqnarray*}
This shows that the covariance matrix $$Q_t:=\Pa_{\eta_t}:=\EE\left(\widetilde{X}_t\widetilde{X}_t^{\prime}~|~\Fa_t\right)=\EE\left(\widetilde{X}_t\widetilde{X}_t^{\prime}\right)$$ does not depend on the observation process. In addition, taking the expectations in the above displayed formula
\begin{eqnarray*}
\partial_tQ_t&=&AQ_t+Q_tA^{\prime}-\varphi_t(P_0)SQ_t-Q_tS\varphi_t(P_0)+\varphi_t(P_0)S\varphi_t(P_0)+R.
\end{eqnarray*}
We set
$$
U_t:=
Q_t-\varphi_t(P_0).
$$
In this notation, we find that
\begin{eqnarray*}
\partial_tU_t&=&AU_t+U_tA^{\prime}-\varphi_t(P_0)SU_t-U_tS\varphi_t(P_0)=\left(A-\varphi_t(P_0)S\right)U_t+U_t\left(A-\varphi_t(P_0)S\right)^{\prime}.
\end{eqnarray*}
The solution is given by
$$
U_t=\exp{\left(\oint_0^t\left(A-\Pa_{\varphi_s(\eta_0)}S\right)ds\right)}~U_0~\left(
\exp{\left(\oint_0^t\left(A-\Pa_{\varphi_s(\eta_0)}S\right)ds\right)}\right)^{\prime}.
$$
This shows that
$$
Q_0=\Pa_{\eta_0}~\Longleftrightarrow~\forall t\geq 0\quad
Q_t=\varphi_t(P_0).
$$
This ends the proof of the lemma.\cqfd

\subsection{Stable signal processes}\label{global-contraction-section}

This short section provides some rather elementary contraction inequalities when the drif-matrix of the signal process 
is stable w.r.t. the log-norm. Next proposition presents some {\em global} Lipschitz property.

\begin{prop}\label{prop-stability-1}
For any time horizon $t\geq 0$ we have the Lipschitz properties
 \begin{equation}\label{f-easy-intro}
(\mbox{\rm S})\quad\mbox{and}\quad \mu(A)<0~\Longrightarrow~
 \Vert P_t-\widecheck{P}_t\Vert_{F}\leq \exp{\left(2\mu(A)t\right)}~
\Vert P_0-\widecheck{P}_0\Vert_{F}.
\end{equation}
In addition for any $n\geq 1$ we have
\begin{eqnarray}
\WW_{2n}\left(
\pi_t,\breve{\pi}_t
\right)&\leq&
\exp{\left(\mu(A)t/2\right)}~\left[\WW_{2n}\left(\pi_0,\breve{\pi}_0\right)+c~ \Vert P_0-\widecheck{P}_0\Vert_{F}\right]
\label{prop-stability-1-2}
\end{eqnarray}
for some finite constant $c$.
\end{prop}

The detailed proofs of \eqref{f-easy-intro} and \eqref{prop-stability-1-2} are provided in Section~\ref{Lip-global-sec}.

Rewritten in terms of the Riccati semigroup, by \eqref{f-easy-intro-prop-mean} we have
$$
\eqref{f-easy-intro}\Leftrightarrow
\vert P_t-\widecheck{P}_t\vert=\vert\Pa_{\eta_t}-\Pa_{\breve{\eta}_t}\vert=\vert\varphi_{0,t}(\Pa_{\pi_0})-\varphi_{0,t}(\Pa_{\breve{\pi}_0})\vert\leq e^{2\mu(A)t}~
\Vert \Pa_{\pi_0}-\Pa_{\breve{\pi}_0}\Vert_{F}.
$$

Of course there exist many distributions with a prescribed covariance matrix. Next lemma provides some Lipschitz 
properties of the trace and the Frobenius norm  w.r.t. the Wasserstein metric. These properties allow to quantify
the continuity property of the covariation matrices w.r.t. a given distribution.

\begin{lem}\label{lemma-Lip-Cov}
For any probability distributions $(\pi_0,\check{\pi}_0)$ on $\RR^{r_1}$ we have the 
regularity property
$$
\left(4^{-1}\left\vert\mbox{\rm tr}(\Pa_{\pi_0}-\Pa_{\check{\pi}_0})\right\vert\right)\vee
\Vert \Pa_{\pi_0}-\Pa_{\check{\pi}_0}\Vert_F\leq~\WW_2(\pi_0,\check{\pi}_0)
~\Vert\check{\pi}_0(e_2)\Vert^{1/2}+~2^{-1/2}~\WW_2(\pi_0,\check{\pi}_0)^2
$$
with the function $x=(x_i)_{1\leq i\leq r_1}\in\RR^{r_1}\mapsto e_2(x):=(x_i^2)_{1\leq i\leq r_1}\RR^{r_1}$.
\end{lem}
The proof of this lemma is rather technical and lengthy, thus it its housed in the appendix on page~\pageref{proof-lemma-Lip-Cov}. 

Lemma~\ref{lemma-Lip-Cov} can be used to deduce several functional contraction inequalities w.r.t the Wasserstein distance between the initial distributions
of the Kalman-Bucy diffusion.
For instance, combining Lemma~\ref{lemma-Lip-Cov} with \eqref{prop-stability-1-2} we readily obtain the following proposition
\begin{prop}
Assume that $ \mu(A)<0$ and $(\mbox{\rm S})$ is satisfied. In this case,
for any $t\geq 0$, the following nonlinear functional inequality holds :
$$
\WW_{2}\left(
\pi_t,\breve{\pi}_t
\right)\leq c~e^{\mu(A)t/2}~\WW_2(\pi_0,\check{\pi}_0)~\left[
1+
\Vert\check{\pi}_0(e_2)\Vert^{1/2}+
\WW_2(\pi_0,\check{\pi}_0)~\right]
$$ 
for some finite constant $c$.
\end{prop}

\subsection{Unstable signal processes}\label{local-contraction}

The two main theorems stated in section~\ref{stab-nl-kb-dffusions} (theorems~\ref{prop-stab-Kalman} and~\ref{theo-Ent}) 
show that the nonlinear Kalman-Bucy diffusions can be stable even when the drift-matrix of the signal is unstable. 
The proof of these stability properties rely on Bucy's analysis of the Riccati equation. 

The following theorem is a direct consequence of the Lyapunov inequalities and the uniform spectral estimates stated in 
lemma 4 and 5 and theorem 4,  in the pioneering article by R.S. Bucy~\cite{bucy}.

\begin{theo}[Bucy~\cite{bucy}]\label{theo-Bucy}
When the filtering problem is uniformly observable and controllable,  for any $ t\geq s\geq \upsilon$ we have
the uniform estimates
$$
\sup_{P_0\in\SS^+_{r_1}}\left\Vert  \exp{\left[\oint_s^t(A-P_uS)du\right]}\right\Vert_2
\leq \alpha_{\upsilon}
~\exp{\left\{-\beta_{\upsilon}(t-s)\right\}}
$$
for some parameters $\alpha_{\upsilon}<\infty$ and $\beta_{\upsilon}>0$. In addition,   for any $t\geq 0$ we have
\begin{equation}\label{bucy-estimate}
\Vert P_t-\widecheck{P}_t\Vert_2
\leq \alpha_{\upsilon}(P_0,\widecheck{P}_0)
~\exp{\left\{-2\beta_{\upsilon}t\right\}}~\Vert P_0-\widecheck{P}_0\Vert_2
\end{equation}
for some constant $\alpha_{\upsilon}(P_0,\widecheck{P}_0)$ whose values only depend on $(\upsilon,P_0,\widecheck{P}_0))$.
\end{theo}
These important contributions were published in 
1967 by R.S. Bucy in~\cite{bucy2}. 

 Combining theorem~\ref{theo-Bucy} with the perturbation lemma~\ref{perturbation-lemma-intro} 
we find the following corollary.
\begin{cor}\label{cor-1-bucy}
Under the assumptions of theorem~\ref{theo-Bucy} for any $\epsilon\in ]0,1]$, any $P_0\in\SS_{r_1}^+$, and any $s\leq t$ we have the exponential semigroup estimates
\begin{equation}\label{with-varsigma}
 \left\Vert \exp{\left[\oint_s^t(A-P_uS)du\right]}\right\Vert_2\leq~\overline{\kappa}_{\epsilon,\varsigma}(P_0,\upsilon)~\exp{\left((1-\epsilon)\varsigma(A-PS)~(t-s)\right)}
\end{equation}
and
\begin{equation}\label{with-mu}
 \left\Vert \exp{\left[\oint_s^t(A-P_uS)du\right]}\right\Vert_2\leq~\overline{\kappa}_{\mu}(P_0,\upsilon)~\exp{\left(\mu(A-PS)~(t-s)\right)}
\end{equation}
In the above displayed formulae the finite constants $\overline{\kappa}_{\mu}(P_0,\upsilon)$ and $\overline{\kappa}_{\epsilon,\varsigma}(P_0,\upsilon)$ defined by
$$
\log{\overline{\kappa}_{\mu}(P_0,\upsilon)}=
 \kappa(\epsilon)^{-1}\log{\left[\overline{\kappa}_{\epsilon,\varsigma}(P_0,\upsilon)/\kappa(\epsilon)\right]}=\Vert P_0-P\Vert_2\Vert S\Vert_2~ \alpha_{\upsilon}(P_0,P)/(2\beta_{\upsilon})
$$
with the parameters $(\kappa(\epsilon), \alpha_{\upsilon}(P_0,P),\beta_{\upsilon})$ presented in  (\ref{common-estimate-Jordan-Schur}) and (\ref{bucy-estimate}).

\end{cor}
We also have
 \begin{equation}\label{implicit-formulation}
 \begin{array}{l}
\partial_t(P_t-\widecheck{P}_t)=(A-\widecheck{P}_tS)~(P_t-\widecheck{P}_t)+(P_t-\widecheck{P}_t)(A-P_tS)^{\prime}\\
\\
\displaystyle\Rightarrow(P_t-\widecheck{P}_t)=\exp{\left(\oint_s^t (A-\widecheck{P}_uS)du\right)}~(P_s-\widecheck{P}_s)~\left[\exp{\left(\oint_s^t (A-P_uS)du\right)}\right]^{\prime}
\end{array}
\end{equation}
 This readily implies the following result.
 \begin{cor}\label{cor-bucy-entended}
Under the assumptions of theorem~\ref{theo-Bucy} for any $\epsilon\in ]0,1]$ and any $t\geq 0$ we have the exponential semigroup estimates
\begin{equation}\label{with-varsigma}
 \left\Vert P_t-\widecheck{P}_t\right\Vert_2\leq~\overline{\kappa}_{\epsilon,\varsigma}(P_0,\upsilon)~\overline{\kappa}_{\epsilon,\varsigma}(\widecheck{P}_0,\upsilon)~\exp{\left(2(1-\epsilon)\varsigma(A-PS)t\right)}~ \left\Vert P_0-\widecheck{P}_0\right\Vert_2
\end{equation}
as well as
\begin{equation}\label{with-mu}
 \left\Vert P_t-\widecheck{P}_t\right\Vert_2\leq~\overline{\kappa}_{\mu}(P_0,\upsilon)
 \overline{\kappa}_{\mu}(\widecheck{P}_0,\upsilon)~\exp{\left(2\mu(A-PS)t\right)}~ \left\Vert P_0-\widecheck{P}_0\right\Vert_2
\end{equation}
with functions $Q\mapsto\overline{\kappa}_{\mu}(Q,\upsilon)$ and $\overline{\kappa}_{\varsigma}(Q,\upsilon)$ defined in corollary~\ref{cor-1-bucy}.
\end{cor}

The article~\cite{ocone-pardoux} also provides a similar exponential decay when $\widecheck{P}_0=P$, without the Lipschitz property 
w.r.t.
the initial covariance matrix, and with half of the order of the rate of decays to equilibrium stated above.

For completeness and to better connect our work with existing literature on Riccati differential matrix equations we end this section with some comments
on the contraction theory of Riccati flows w.r.t the Thompson metric. 

We recall that the Thompson's metric
(a.k.a. part metric) on the space of definite positive matrices $P_1,P_2$ is defined by
$$
d_{T}(P_1,P_2)=\log{\max{\left(M(P_1/P_2),M(P_2/P_1)\right)}}$$
with
$$
M(P_1/P_2):=\inf{\left\{u\geq~0~:~P_1\leq u~P_2\right\}}
$$
By a recent article by D.A. Snyder~\cite{snyder} we have
 \begin{eqnarray*}
\Vert P_1-P_2\Vert_F&\leq& \left(\exp{\left[d_{T}(P_1,P_2)\right]}-1\right)~\sqrt{\frac{\Vert P_1\Vert_F^2+\Vert P_2\Vert_F^2}{1+\exp{\left(2d_{T}(P_1,P_2)\right)}}}\\
&\leq & e~d_{T}(P_1,P_2)~\sqrt{\Vert P_1\Vert_F^2+\Vert P_2\Vert_F^2}
 \end{eqnarray*}
The last assertion is valid as soon as $d_{T}(P_1,P_2)\leq 1$.

Let $(P_t,\widecheck{P}_t)$ be two solutions of the Riccati equation starting at some possibly different states
$(P_0,\widecheck{P}_0)$ such that
$$
\alpha_1^{-1}~\widecheck{P}_0~\leq P_0\leq \alpha_2~\widecheck{P}_0\Longrightarrow M(P_0/\widecheck{P}_0)\leq \alpha_2\quad
\mbox{\rm and}\quad M(\widecheck{P}_0/P_0)\leq \alpha_1
$$
for some $\alpha_1,\alpha_2>0$.
By  theorem 8.5 in~\cite{lawson}, for any $\beta>0$ we have the contraction inequality
 \begin{eqnarray}
  R\geq \beta~S&\Longrightarrow&
 d_{T}(P_t,\widecheck{P}_t)\leq e~\exp{\left(-2\sqrt{\beta}~ t\right)}~ d_{T}(P_0,\widecheck{P}_0)\leq e~
 \exp{\left(-2\sqrt{\beta}~ t\right)}~ \left(\alpha_1\vee\alpha_2\right)\nonumber\\
 &\Longrightarrow&\Vert  P_t-\widecheck{P}_t\Vert_F\leq e~\exp{\left(-2\sqrt{\beta}~ t\right)~ \left(\alpha_1\vee\alpha_2\right)}~\sqrt{\Vert P_t\Vert_F^2+\Vert \widecheck{P}_t\Vert_F^2}\label{PtP-ref}
 \end{eqnarray}
 as soon as $t\geq ~\left(2\sqrt{\beta}\right)^{-1}\log{ \left(\alpha_1\vee\alpha_2\right)}$. Choosing $~\widecheck{P}_0=P$ and
 $\alpha_1=1\leq \alpha:=\alpha_2$  we conclude that
\begin{equation}\label{PtP-ref-2}
 \begin{array}{l}
 P\leq P_0\leq \alpha~P\quad \mbox{\rm and}\quad
   R\geq \beta~S\\
   \\
   \Longrightarrow
   \Vert  P_t-P\Vert_F\leq \alpha e~e^{-2\sqrt{\beta}~ t}~\left(2\Vert P\Vert_F+e^{2\mu(A-PS)t}~\Vert P_0-P\Vert_F\right)
\end{array} 
\end{equation}
 
The estimates (\ref{PtP-ref})  and (\ref{PtP-ref-2}) are useful as soon $\Vert P_t\Vert_F$ is uniformly bounded. This property is ensured when the filtering problem is uniformly observable and controllable. In this situation, the exponential rate to equilibrium given in (\ref{PtP-ref-2}) is related to a signal to noise ratio associated with the pair of matrices $(R,S)$.

We have derived a series of quantitive estimates for Kalman-Bucy diffusions. These estimates can be used to analyze the stability properties of Kalman-Bucy filters. For a more thorough discussion and a more recent account
on the stability of discrete generation Kalman filters we refer to~\cite{costa,bougerol,sua} and the references therein.
See also the pioneering article of Anderson~\cite{anderson}, the one by  Ocone and 
Pardoux~\cite{ocone-pardoux} on the stability of continuous time Kalman-Bucy filters, as well as the book by H. Kwakernaak, R. Sivan~\cite{kwaker-sivan}.

\section{A brief review on Ornstein-Ulhenbeck processes and Riccati equations}\label{preliminary-sec}

\subsection{Some uniform moment estimates}\label{some-uniform-moments-sec}

Our analysis on the convergence of the EnKF
requires that $\mu(A)<0$. This is not really surprising.
The EnKF is designed in terms of interacting covariance matrices and interacting Monte Carlo samples
based on the signal evolution. When $\mu(A)\geq 0$ the signal contains an unstable component. 
In this case the fluctuations induced by the Monte Carlo samples may increase dramatically the 
global error variances. 
To analyze these interacting filters based on an extra level of randomness
we need to strengthen the usual condition $\mu(A-PS)<0$ discussed in Section~\ref{stab-section-p} to ensure that the signal itself is stable. 

Next we show that the condition $\mu(A)<0$ cannot be relaxed.
In the one dimensional case when $C=0$ the EnKF resumes to $N$ independent copies of the signal
and $\widehat{X}_t=\EE(X_t)$. In this case we have
$$
\mu(A)\geq 0~\Longleftrightarrow~
N~\EE\left[\left(m_t-\widehat{X}_t\right)^2\right]=R(2\mu(A))^{-1}\left(e^{2\mu(A)t}-1\right)~\longrightarrow_{t\uparrow\infty}~\infty.
$$
When $\mu(A)=0$ we use the convention $(2\mu(A))^{-1}\left(e^{2\mu(A)t}-1\right)=t$.
This shows that, even for one dimensional Brownian signal motions it is hopeless to try to find some uniform estimates for the sample mean.

Next we provide a brief discussion on the stability of multi-dimensional  signal processes.
For multidimensional filtering problems, the solution of the signal
stochastic differential equation is given by the Ornstein-Uhlenbeck formula
$$
X_t=e^{At}~X_0-\left({\rm Id}-e^{At}\right)A^{-1}a+\int_0^t~e^{A(t-s)}R^{1/2}~dW_s.
$$
The mean vector and the covariance matrix 
are given by
$$
\left(\EE(X_t)+A^{-1}a\right)=e^{tA}~\left(\EE(X_0)+A^{-1}a\right)\quad\mbox{\rm and}
\quad P^X_t=e^{tA}P_0~e^{A^{\prime}t}+\int_0^t~e^{As}Re^{A^{\prime}s}~ds.
$$
In signal processing and control theory, the integral in the r.h.s. term is called the controllability Grammian. 
Recalling that $\left\Vert  \exp{\left(tA\right)}\right\Vert\leq \exp{\left(\mu(A)t\right)}$ we find that
$$
\left\Vert\EE(X_t)+A^{-1}a\right\Vert\leq \exp{\left(\mu(A)t\right)}~\left\Vert\EE(X_0)+A^{-1}a\right\Vert
$$
and
$$
\left\Vert \int_t^{\infty}~e^{As}Re^{A^{\prime}s}~ds\right\Vert \leq\exp{\left(2t\mu(A)\right)}~\frac{\left\Vert R\right\Vert}{2\left\vert\mu(A)\right\vert}.
$$
Recall that $-A^{-1}=\int_0^{\infty}~e^{tA}dt\Rightarrow\Vert A^{-1}\Vert\leq 1/\vert\mu(A)\vert$.
Thus,  the condition $\mu\left(A\right)<0$ (and  of course course $\Vert a\Vert<\infty$) ensures that 
$$
\lim_{t\rightarrow\infty}
\left\Vert\EE(X_t)+A^{-1}a\right\Vert=0=
\lim_{t\rightarrow\infty}\left\Vert P^X_t-\int_0^{\infty}e^{As}Re^{A^{\prime}s}~ds\right\Vert.
$$
It also yields the uniform moment estimates
\begin{equation}\label{fv1-4}
\sup_{t\geq 0}{\EE\left(\Vert \widehat{X}_t\Vert^{n}\right)}<\infty \quad\mbox{\rm and}
\quad
\sup_{t\geq 0}{\EE\left(\left\Vert X_t\right\Vert^{n}\right)}<\infty
\quad\mbox{\rm as well as}
\quad
\sup_{t\geq 0}\left\Vert P^X_t\right\Vert_F<\infty
\end{equation}
for any $n\geq 1$.
The last assertion is easily checked using Bernstein inequality 
$$
\mbox{\rm tr}\left(e^{A}~e^{A^{\prime}}\right)\leq \mbox{\rm tr}\left(e^{A+A^{\prime}}\right)\quad\left(\leq r_1~e^{2\mu(A)}\right).
$$
A proof of this inequality result can be found in~\cite{bernstein}, see also~\cite{yang,yang2} for a more thorough discussion on trace inequalities.

\subsection{Riccati equations}\label{riccati-section}
The Riccati Equation \eqref{nonlinear-KB-Riccati} can be solved analytically when for non observed or
noise free signals (i.e. $C=0$ or $R=0$). The situation $C=0$ has already been discussed above.
In this case,  $P_t=P^X_t$ resumes to the covariance matrix $P^X_t$ of the signal process.

When $R=0$ solution is given by
$$
P_t=e^{tA}\left(P_0^{-1}+\int_0^te^{sA^{\prime}}S~e^{sA}~ds\right)^{-1}e^{tA^{\prime}}
\quad\left(\Rightarrow
P_t=\left(P_0^{-1}+St\right)^{-1}\quad\mbox{\rm when $A=0$}\right).
$$
In more general situations we need to resort to some numerical scheme or to some additional algebraic development such as
the Bernoulli substitution approach to reduce the problem to an ordinary linear differential equation in $2r_1$-dimensions. For one dimensional signal processes ($r_1=1$) $P_t$ coincides with the 
variance between $\widehat{X}_t$ and $X_t$. 
When $S\not=0\not=R$ the Riccati equation
takes the form
$$
\partial_tP_t=-S~\left(P_t-z_1\right)~\left(P_t-z_2\right)
$$
with the couple of roots
$$
z_1=\frac{A-\sqrt{A^2+SR}}{S}<0<z_2=\frac{A+\sqrt{A^2+S R}}{S}.
$$
The solution
 is given by the formula
\begin{equation}\label{fv1-1}
P_t-z_2=(P_0-z_2)~\frac{(z_2-z_1)~e^{-2t\sqrt{A^2+S R}}}{(z_2-P_0)~e^{-2t\sqrt{A^2+S R}}+(P_0-z_1)}~\longrightarrow_{t\rightarrow\infty}~0.
\end{equation}
The above formula underline the fact that the Riccati equation is stable even for unstable signals, that is when $A>0$. Also observe that
\begin{equation}\label{fv1-1-2}
0\leq P_t\leq z_2+\left(P_0-z_2\right)_+~e^{-2t\sqrt{A^2+S R}}
\end{equation}
as soon as $P_0\geq 0$, where $a_+:=\max{(a,0)}$, for any $a\in\RR$. We check this claim using the decomposition
$$
P_t-z_2=(P_0-z_2)~\frac{(z_2-z_1)~e^{-2t\sqrt{A^2+S R}}}{(P_0-z_2)~(1-e^{-2t\sqrt{A^2+S R}})+(z_2-z_1)}
$$
Assume that $S\not=0$ and let $P=z_2$ the unique positive fixed point. 
In this case, we notice that $-\sqrt{A^2+S R}=A-PS<0$ iff $A^2\wedge R>0$ ; and
 $A-PS/2=z_1S/2<0$ iff $R>0$. Our regularity assumption $A-PS/2<0$ cannot capture the case $R=0$ and $A\not=0$. 
 
Another direct consequence of this result is that  the minimum variance function $t\mapsto P_t$ is uniformly bounded w.r.t. the time parameter. For matrix valued Riccati equation we can use the following comparison lemma.
\begin{lem}\label{Lem-v1-1}
We assume that $\mu(A)<0$. In this situation, we have
$
\mbox{\rm tr}(P_t)\leq g_t
$
where $g_t$ stands for the solution of the Riccati equation
$$
 \partial_t g_t=2~\alpha~g_t- \beta~g_t^2+r\quad\mbox{starting at $g_0=\mbox{\rm tr}(P_0)$}
$$
with the parameters $\left(\alpha,\beta,r\right)=\left(\mu(A),r_1^{-1}\lambda_{\tiny min}(S),\mbox{\rm tr}(R)\right)$.

\end{lem}
\proof

The key idea is to use the commutation inequality
\begin{equation}\label{f32}
\begin{array}{l}
\mbox{\rm tr}\left(AP_t+P_tA^{\prime}-P_tSP_t+R\right)\\
\\
=2~\mbox{\rm tr}\left(A_{\tiny sym}P_t\right)-\mbox{\rm tr}\left(SP^2_t\right)+\mbox{\rm tr}\left(R\right)
\leq 
2~\alpha~\mbox{\rm tr}(P_t)- \beta~\left(\mbox{\rm tr}(P_t)\right)^2+r.
\end{array}
\end{equation}
In the last display we have used \eqref{f30}.
This yields the Riccati differential inequality
$$
\partial_t\mbox{\rm tr}(P_t)\leq 2~\alpha~\mbox{\rm tr}(P_t)- \beta~\left(\mbox{\rm tr}(P_t)\right)^2+r
$$
from which we conclude that
$$
\mbox{\rm tr}(P_t)\leq g_t+e^{2\alpha t}~\left(\mbox{\rm tr}(P_0)-g_0\right)\quad\mbox{\rm with}\quad
 \partial g_t=2~\alpha~g_t- \beta~g_t^2+r.
$$
This ends the proof of the lemma.\cqfd

Using theorem~\ref{prop-stab-Kalman} or Lemma~\ref{Lem-v1-1} we readily deduce the following uniform estimates:
\begin{equation}\label{fv1-2}
\mu(A)<0\Longrightarrow~\sup_{t\geq 0}{\mbox{\rm tr}(P_t)}<\infty\qquad\mbox{or equivalently}\qquad \sup_{t\geq 0}{\left\Vert P_t\right\Vert_F}<\infty.
\end{equation}
The l.h.s. inequality in \eqref{fv1-2} is proven using the same analysis as the one of the scalar Riccati Equation \eqref{fv1-1}.
The equivalence property in  \eqref{fv1-2} is a direct consequence of the fact that $\left\Vert P_t\right\Vert_F=\mbox{\rm tr}(P_t^2)\leq \left(\mbox{\rm tr}(P_t)\right)^2\leq r_1~\left\Vert P_t\right\Vert_F$.  The estimate (\ref{fv1-1-2}) also implies
that
\begin{eqnarray*}
0\leq \mbox{\rm tr}(P_t)&\leq&  r_1\mu(A)~\frac{1+\delta}{\lambda_{\tiny min}(S)}+\left(\mbox{\rm tr}(P_0)- r_1\mu(A)~\frac{1+\delta}{\lambda_{\tiny min}(S)}\right)_+~\exp{\left[2t\mu(A)\delta\right]}
\end{eqnarray*}
with the parameter
$
\delta:=\sqrt{1+r_1^{-1}\left[\mbox{\rm tr}(R)\lambda_{\tiny min}(S)\right]/\mu(A)^2}
$.
These estimates are useful as soon as $P_0\geq P$. When $P_0\leq P$ we clearly have $\mbox{\rm tr}(P_t)\leq \mbox{\rm tr}(P)$.

\section{The Ensemble Kalman-Bucy filter equations}\label{sec-EnFK}
\subsection{Sample mean and Covariance diffusions}\label{sec-EnFK-ss1}
This section is mainly concerned with the proof of the stochastic differential Equations \eqref{f21} and \eqref{f22}.

The stochastic diffusion equation of the EnKF sample mean \eqref{f21} is easily checked using \eqref{fv1-3}
with the $r_1$-multidimensional martingale defined by
$$
d\overline{M}_t=\frac{1}{\sqrt{N}}\sum_{1\leq i\leq N}R^{1/2}~d\overline{W}_t^i-\frac{1}{\sqrt{N}}\sum_{1\leq i\leq N}~p_t~C^{\prime}~R^{-1/2}_{2}~d\overline{V}^i_{t}.
$$
This clearly implies \eqref{f22}. Using  \eqref{fv1-3} we also readily check that
$$
d\left(\xi^i_t-m_t\right)=\left(A-p_tS\right)\left(\xi^i_t-m_t\right)~dt+dM^i_t
$$
with the $r_1$-dimensional martingale
$$
dM^i_t:=R^{1/2}d\widetilde{W}^i_t-p_tC^{\prime}
R^{-1/2}_2d\widetilde{V}^i_t
$$
defined in terms of the diffusion processes
$$
\widetilde{W}^i_t=\overline{W}^i_t-N^{-1}\sum_{1\leq j\leq N}\overline{W}^j_t\quad\mbox{\rm and}
\quad
\widetilde{V}^i_t=\overline{V}^i_t-N^{-1}\sum_{1\leq j\leq N}\overline{V}^j_t.
$$
The angle-brackets $\langle M^i(k),M^j(k^{\prime})\rangle_t$ associated with the collection 
of vector valued martingales $M^i_t$ are given by the formulae 
\begin{equation}\label{details-ab-1}
\partial_t\langle M^i(k),M^i(k^{\prime})\rangle_t=\left(1-\frac{1}{N}\right)~\left(R+p_tSp_t\right)(k,k^{\prime})
\end{equation}
and for $i\not=j$
\begin{equation}\label{details-ab-2}
\partial_t\langle M^i(k),M^j(k^{\prime})\rangle_t=-\frac{1}{N}~\left(R+p_tSp_t\right)(k,k^{\prime}).
\end{equation}
To check this claim we observe that
\begin{eqnarray*}
dM^i_t&:=&\left(1-\frac{1}{N}\right)~\left[R^{1/2}d\overline{W}^i_t-p_tC^{\prime}
R^{-1/2}_2d\overline{V}^i_t\right]-\frac{1}{N}\sum_{j\not =i}\left[R^{1/2}d\overline{W}^j_t-p_tC^{\prime}
R^{-1/2}_2d\overline{V}^j_t\right]\\
\\
&=&\sum_{j}\epsilon^i_j~\left[R^{1/2}d\overline{W}^j_t-p_tC^{\prime}
R^{-1/2}_2d\overline{V}^j_t\right]
\end{eqnarray*}
with
$$
\epsilon^i_j=1_{i=j}\left(1-\frac{1}{N}\right)-\frac{1}{N}~1_{i\not= j}
$$
Therefore
\begin{eqnarray*}
dM^i_t(k)
&=&\sum_{j,l}\epsilon^i_jR^{1/2}(k,l)d\overline{W}^j_t(l)-\sum_{j,l}\epsilon^i_j (p_tC^{\prime}
R^{-1/2}_2)(k,l)d\overline{V}^j_t(l)
\end{eqnarray*}
from which we conclude that
\begin{eqnarray*}
\partial_t\langle M^i(k),M^i(k^{\prime})\rangle_t
&=&\sum_{j,l}\sum_{j^{\prime},l^{\prime}}\epsilon^i_j\epsilon^i_{j^{\prime}}
R^{1/2}(k,l)R^{1/2}(k^{\prime},l^{\prime})~1_{j=j^{\prime}}~1_{l=l^{\prime}}\\
&&
\hskip1cm+\sum_{j,l}\sum_{j^{\prime},l^{\prime}}\epsilon^i_j\epsilon^i_{j^{\prime}}
(p_tC^{\prime}
R^{-1/2}_2)(k,l)(p_tC^{\prime}
R^{-1/2}_2)
(k^{\prime},l^{\prime})~1_{j=j^{\prime}}~1_{l=l^{\prime}}\\
&=&\sum_{j}(\epsilon^i_j)^2~\left[\sum_{l}R^{1/2}(k,l)R^{1/2}(k^{\prime},l)\right.\\
&&\hskip4cm\left.+\sum_{j,l}
(p_tC^{\prime}
R^{-1/2}_2)(k,l)(p_tC^{\prime}
R^{-1/2}_2)
(k^{\prime},l)\right]\\
&=&\sum_{j}(\epsilon^i_j)^2~(R-p_tSp_t)(k,k^{\prime})
\end{eqnarray*}
In the last assertion we have used the symmetry of the matrices $R^{1/2}$ and
$$
(p_tC^{\prime}
R^{-1/2}_2)
(k^{\prime},l)=(p_tC^{\prime}
R^{-1/2}_2)^{\prime}
(l,k^{\prime})=(R^{-1/2}_2Cp_t)(l,k^{\prime})
$$
Observe that
$$
\sum_{j}(\epsilon^i_j)^2=\left(1-\frac{1}{N}\right)^2+\sum_{j\not=i}\frac{1}{N^2}=
\left(1-\frac{1}{N}\right)\left(\left(1-\frac{1}{N}\right)+\frac{1}{N}\right)=\left(1-\frac{1}{N}\right)
$$
When $i\not=i^{\prime}$ we have
\begin{eqnarray*}
\partial_t\langle M^i(k),M^{i^{\prime}}(k^{\prime})\rangle_t
&=&\sum_{j}\epsilon^i_j\epsilon^{i^{\prime}}_j~(R-p_tSp_t)(k,k^{\prime})
\end{eqnarray*}
and
\begin{eqnarray*}
\sum_{j}\epsilon^i_j\epsilon^{i^{\prime}}_j&=&\epsilon^i_i\epsilon^{i^{\prime}}_i+\epsilon^{i}_{i^{\prime}}\epsilon^{i^{\prime}}_{i^{\prime}}+\sum_{j\not\in{\{i,i^{\prime}\}}}\epsilon^i_j\epsilon^{i^{\prime}}_j\\
&=&-2~\frac{1}{N}\left(1-\frac{1}{N}\right)~+(N-2)~\frac{1}{N^2}=
\frac{1}{N}\left(\left(\frac{2}{N}-2\right)+\left(1-\frac{2}{N}\right)\right)=-\frac{1}{N}
\end{eqnarray*}
This ends the proof of the angle-bracket Formulae \eqref{details-ab-1} and \eqref{details-ab-2}.

This implies that
$$
\begin{array}{l}
d\left[\left(\xi^i_t-m_t\right)\left(\xi^i_t-m_t\right)^{\prime}\right]\\
\\
=\left(A-p_tS\right)\left(\xi^i_t-m_t\right)\left(\xi^i_t-m_t\right)^{\prime}~dt
+\left(\xi^i_t-m_t\right)\left(\xi^i_t-m_t\right)^{\prime}\left(A^{\prime}-Sp_t\right)~dt\\
\\
+\left(1-N^{-1}\right)~\left(R+p_tSp_t\right)~dt+\left(\xi^i_t-m_t\right)\left(dM^i_t\right)^{\prime}+dM^i_t\left(\xi^i_t-m_t\right)^{\prime}.
\end{array}
$$
Summing the indices we find that 
\begin{eqnarray*}
dp_t
&=&\left[\left(A-p_tS\right)p_t+p_t\left(A^{\prime}-Sp_t\right)+\left(R+p_tSp_t\right)\right]~dt+\frac{1}{\sqrt{N-1}}~dM_t
\end{eqnarray*}
with
$$
dM_t:=
\frac{1}{\sqrt{N-1}}\sum_{1\leq i\leq N}\left[
\left(\xi^i_t-m_t\right)\left(dM^i_t\right)^{\prime}+dM^i_t\left(\xi^i_t-m_t\right)^{\prime}\right].
$$
This ends the proof of \eqref{f23}. To check \eqref{f24} we set
$
\epsilon^i_t:=\xi^i_t-m_t
$.
In this notation we have
$$
dM_t(k,l):=
\frac{1}{\sqrt{N-1}}\sum_{1\leq i\leq N}
\left[
\epsilon^i_t(k)~dM^i_t(l)+dM^i_t(k)\epsilon^i_t(l)\right].
$$
This implies that
$$
\begin{array}{l}
(N-1)\partial_t\langle M(k,l),M(k^{\prime},l^{\prime})\rangle\\
\\
=
\displaystyle\sum_{1\leq i,i^{\prime}\leq N}
\left[
\epsilon^i_t(k)\epsilon^{i^{\prime}}_t(k^{\prime})~\partial_t\langle M^i(l),M^{i^{\prime}}(l^{\prime})\rangle_t+\epsilon^i_t(l)\epsilon^{i^{\prime}}_t(l^{\prime})~\partial_t\langle M^i(k), M^{i^{\prime}}(k^{\prime})
\rangle_t
\right]\\
\\ 
+\displaystyle\sum_{1\leq i,i^{\prime}\leq N}\left[
\displaystyle\epsilon^i_t(k)\epsilon^{i^{\prime}}_t(l^{\prime})
~\partial_t \langle M^i(l) , M^{i^{\prime}}(k^{\prime})
\rangle_t
+
\epsilon^i_t(l)\epsilon^{i^{\prime}}_t(k^{\prime})~\partial_t \langle M^i(k), M^{i^{\prime}}(l^{\prime})\rangle_t
\right]
\end{array}
$$
By (\ref{details-ab-1}) and (\ref{details-ab-2}) we have
$$
\begin{array}{l}
(N-1)\partial_t\langle M(k,l),M(k^{\prime},l^{\prime})\rangle\\
\\
=
\displaystyle\left(1-\frac{1}{N}\right)\sum_{1\leq i\leq N}
\left[\epsilon^i_t(k)\epsilon^{i}_t(k^{\prime})~\left(R+p_tSp_t\right)(l,l^{\prime})+\epsilon^i_t(l)\epsilon^{i}_t(l^{\prime})~\left(R+p_tSp_t\right)(k,k^{\prime})
\right]\\
\\ 
\hskip1cm+\displaystyle \left(1-\frac{1}{N}\right)\sum_{1\leq i\leq N}\left[
\displaystyle\epsilon^i_t(k)\epsilon^{i}_t(l^{\prime})
\left(R+p_tSp_t\right)(l,k^{\prime})
+
\epsilon^i_t(l)\epsilon^{i}_t(k^{\prime})\left(R+p_tSp_t\right)(k,l^{\prime})
\right]\\
\\
\hskip2cm~\displaystyle-\frac{1}{N}\sum_{1\leq i\not=i^{\prime}\leq N}
\left[
\epsilon^i_t(k)\epsilon^{i^{\prime}}_t(k^{\prime})~~\left(R+p_tSp_t\right)(l,l^{\prime})+\epsilon^i_t(l)\epsilon^{i^{\prime}}_t(l^{\prime})~\left(R+p_tSp_t\right)(k,k^{\prime})
\rangle_t
\right]\\
\\ 
\hskip3cm\displaystyle-\frac{1}{N}~\sum_{1\leq i\not=i^{\prime}\leq N}\left[
\displaystyle\epsilon^i_t(k)\epsilon^{i^{\prime}}_t(l^{\prime})
~\left(R+p_tSp_t\right)(l,k^{\prime})
+
\epsilon^i_t(l)\epsilon^{i^{\prime}}_t(k^{\prime})~\left(R+p_tSp_t\right)(k,l^{\prime})
\right].
\end{array}
$$
Recalling that
$$
\frac{1}{N-1}\sum_{1\leq i\leq N}\epsilon^i_t(k)\epsilon^i_t(k^{\prime})=p_t(k,k^{\prime})
$$
and
$$
\begin{array}{l}
\displaystyle\sum_{1\leq i\leq N}\epsilon^i_t(k)\sum_{1\leq i^{\prime}\leq N}\epsilon^i_t(k^{\prime})=0\\
\\
\Rightarrow
\displaystyle\frac{1}{N-1}\sum_{1\leq i\not=i^{\prime}\leq N}\epsilon^i_t(k)\epsilon^{\prime}_t(k^{\prime})=-\frac{1}{N-1}
\sum_{1\leq i\leq N}\epsilon^i_t(k)\epsilon^i_t(k^{\prime})=p_t(k,k^{\prime})
\end{array}
$$
we find that
$$
\begin{array}{l}
\partial_t\langle M(k,l),M(k^{\prime},l^{\prime})\rangle\\
\\
=
\displaystyle\left(1-\frac{1}{N}\right)
\left[p_t(k,k^{\prime})~\left(R+p_tSp_t\right)(l,l^{\prime})+p_t(l,l^{\prime})~\left(R+p_tSp_t\right)(k,k^{\prime})
\right]\\
\\ 
\hskip1cm+\displaystyle \left(1-\frac{1}{N}\right)\left[
\displaystyle p_t(k,l^{\prime})
\left(R+p_tSp_t\right)(l,k^{\prime})
+
p_t(l,k^{\prime})\left(R+p_tSp_t\right)(k,l^{\prime})
\right]\\
\\
\hskip2cm~\displaystyle+\frac{1}{N}
\left[
p_t(k,k^{\prime})~~\left(R+p_tSp_t\right)(l,l^{\prime})+p_t(l,l^{\prime})~\left(R+p_tSp_t\right)(k,k^{\prime})
\rangle_t
\right]\\
\\ 
\hskip3cm\displaystyle+\frac{1}{N}~\left[
\displaystyle p(k,l^{\prime})
~\left(R+p_tSp_t\right)(l,k^{\prime})
+
p_t(l,k^{\prime})~\left(R+p_tSp_t\right)(k,l^{\prime})
\right].
\end{array}
$$
This ends the proof of  \eqref{f24}.
The last assertion can be checked easily using the fact that
$\langle
\overline{M}(n),M^{i}
\rangle_t
$ does not depend on the index $i$.

This ends the proof of the EnKF differential Equations \eqref{f21} and \eqref{f22}.\cqfd
\subsection{Uniform moments estimates}\label{unif-moment-EnKF-ss2}

The following technical lemma combines a Foster-Lyapunov approach with martingale techniques to control the moments of Riccati type
stochastic differential equations uniformly w.r.t. the time horizon. 
\begin{lem}\label{tech-1}
Let $Z_t$ be some stochastic processes adapted to some filtration $\Fa_{t}$ and taking values in some measurable state 
space $(E,\Ea)$. 
Let $H$ be some non negative measurable function on $(E,\Ea)$ such that
\begin{equation}\label{f39}
dH\left(Z_t\right)= \La_{t}(H)\left(Z_t\right)~dt~+~d\Ma_{t}(H)
\end{equation}
with
 an $\Fa_t$-martingale $\Ma_{t}(H)$ and some $\Fa_t$-adapted process $\La_{t}(H)\left(Z_t\right)$. 
 
 \begin{itemize}
 \item Assume that
 \begin{eqnarray*}
 \La_{t}(H)\left(Z_t\right)&\leq& {2\gamma~\sqrt{H(Z_t)}}+3\alpha~H\left(Z_t\right)-\beta~H\left(Z_t\right)^2+r
\\
\partial_t\left\langle
\Ma(H)
\right\rangle_t&\leq&  ~H\left(Z_t\right)~\left(\tau_0+\tau_1~H\left(Z_t\right)+\tau_2~H\left(Z_t\right)^2\right)
\end{eqnarray*}
for some parameters $\alpha<0$ and ${\gamma},\beta,r,\tau_0,\tau_1,\tau_2\geq 0$. In this situation 
we have the uniform moment estimates
$$
1\leq n<1+2\min{\left(\beta/\tau_2,\vert\alpha\vert/\tau_1\right)}~\Longrightarrow~
\sup_{t\geq 0}
\EE\left(H\left(Z_t\right)^n\right)<\infty
$$
{with the convention $\beta/0=\infty=\vert\alpha\vert/0$ when $\tau_2=0$ or when $\tau_1=0$.}
\item Assume that
 \begin{eqnarray*}
 \La_{t}(H)\left(Z_t\right)&\leq& {2\tau_t(Z_t)~\sqrt{H\left(Z_t\right)}}+2{\alpha}~H\left(Z_t\right)+\beta_{t}\left(Z_t\right)
\\
\partial_t\left\langle
\Ma(H)
\right\rangle_t&\leq& ~H\left(Z_t\right)~\gamma_{t}(Z_t)
\end{eqnarray*}
for some $\alpha<0$ and some non negative functions $(\tau_t,\beta_t,\gamma_t)$ s.t.
 \begin{eqnarray*}
\delta_{\tau,t}(n)&:=&\EE\left(\tau_{t}\left(Z_t\right)^n\right)^{\frac{1}{n}}<\infty \\
\delta_{\beta,t}(n)&:=&\EE\left(\beta_{t}\left(Z_t\right)^n\right)^{\frac{1}{n}}<\infty\quad 
 \mbox{and}\quad
\delta_{\gamma,t}(n):=\EE\left(\gamma_{t}\left(Z_t\right)^n\right)^{\frac{1}{n}}<\infty
\end{eqnarray*}
for any $n\geq 1$.
In this situation, we have the estimate
$$
\begin{array}{l}
\displaystyle\EE\left(H\left(Z_t\right)^n\right)^{\frac{1}{n}}\\
\\
\displaystyle\leq e^{\alpha t}~
\EE\left(H\left(Z_0\right)^n\right)^{\frac{1}{n}} +\int_0^t~e^{\alpha (t-s)}~\left[\left(\delta_{\tau,s}(2n)^2/\vert\alpha\vert+\delta_{\beta,s}(n)+(n-1)\delta_{\gamma,s}(n)/2 \right)\right]~ds.
\end{array}
$$
\end{itemize}
\end{lem}

\proof

{Firstly we observe that
$$
\forall \epsilon>0\,\,:\qquad 2~\gamma~\sqrt{H(Z_t)}\leq \epsilon~H(Z_t)+\frac{1}{\epsilon}~\gamma^2.
$$
Choosing $\epsilon=\vert\alpha\vert>0$ we find that
$$
 \La_{t}(H)\left(Z_t\right)\leq 2\alpha~H\left(Z_t\right)-\beta~H\left(Z_t\right)^2+r+\gamma^2/\vert\alpha\vert.
$$}
For any $n\geq 1$ we have
\begin{eqnarray*}
dH\left(Z_t\right)^n
&=&\La_{n,t}(H)\left(Z_t\right)~dt+d\Ma_{n,t}(H)
\end{eqnarray*}
with the martingale  $d\Ma_{n,t}(H):=n~H\left(Z_t\right)^{n-1}~d\Ma_{t}(H)$ and the drift
\begin{eqnarray*}
\La_{n,t}(H)\left(Z_t\right)&=&n~\left[H\left(Z_t\right)^{n-1}\La_{t}(H)\left(Z_t\right)+\frac{(n-1)}{2}~H\left(Z_t\right)^{n-2}\partial_t\left\langle
\Ma(H)
\right\rangle_t\right]\\
&\leq &2n\alpha_n
H\left(Z_t\right)^{n}
-n\beta_n
H\left(Z_t\right)^{n+1}
+n\rho_n
H\left(Z_t\right)^{n-1}
\end{eqnarray*}
with
$$
\alpha_n=\alpha+\tau_1\frac{(n-1)}{2}<0\,,\quad
\beta_n=\left(\beta-\frac{(n-1)}{2}~\tau_2\right)>0
~~\mbox{\rm and}~~
\rho_n=\left(r+\frac{(n-1)}{2}~\tau_0\right)>0
$$
as soon as
$
1\leq n<1+2\beta\tau_2^{-1}
$. This implies that
\begin{eqnarray*}
\partial_t\EE(H\left(Z_t\right)^n)&\leq& 2\alpha_n
\EE\left(H\left(Z_t\right)^{n}\right)
-\beta_n
\EE\left(H\left(Z_t\right)^{n+1}\right)
+\rho_n
\EE\left(H\left(Z_t\right)^{n-1}\right)\\
&\leq &2\alpha_n~
\EE\left(H\left(Z_t\right)^{n}\right)
-\beta_n~
\EE\left(H\left(Z_t\right)^{n}\right)^{1+\frac{1}{n}}
+\rho_n~
\EE\left(H\left(Z_t\right)^{n}\right)^{1-\frac{1}{n}}
\end{eqnarray*}
and therefore
\begin{eqnarray*}
\partial_t\EE(H\left(Z_t\right)^n)^{\frac{1}{n}}&=& \frac{1}{n}~\EE(H\left(Z_t\right)^n)^{\frac{1}{n}-1}\partial_t\EE(H\left(Z_t\right)^n)\\
&\leq  &2\alpha_n~
\EE\left(H\left(Z_t\right)^{n}\right)^{\frac{1}{n}}
-\beta_n~
\EE\left(H\left(Z_t\right)^{n}\right)^{\frac{2}{n}}
+\rho_n.
\end{eqnarray*}
This shows that 
$$
\EE(H\left(Z_t\right)^n)^{\frac{1}{n}}\leq g_{n,t}+e^{2\alpha_n}~\left(\EE(H_0\left(Z_0\right)^n)^{\frac{1}{n}}-g_{n,0}\right)
$$
with
$$
\partial_tg_{n,t}=2\alpha_n~
g_{n,t}
-\beta_n~
g_{n,t}^2
+\rho_n.
$$
The end of the proof of the first assertion follows the same lines of arguments as the ones of Lemma~\ref{Lem-v1-1}. 

Now we come to the proof of the second assertion.

Arguing as above we have

$$
\forall \epsilon>0\,\,:\quad 2\tau_t(Z_t)~\sqrt{H\left(Z_t\right)}\leq \epsilon~H\left(Z_t\right)+\tau_t(Z_t)^2/\epsilon.
$$
Choosing $\epsilon=\vert\alpha\vert$ we find that
$$
 \La_{t}(H)\left(Z_t\right)\leq \alpha~H\left(Z_t\right)+\beta_{t}\left(Z_t\right)+\tau_t(Z_t)^2/\vert\alpha\vert.
$$
Therefore, there is no loss of generality to assume that $\tau_t(Z_t)=0$ by changing $2\alpha$ by $\alpha$ and $\beta_{t}\left(Z_t\right)$
by $\beta_{t}\left(Z_t\right)+\frac{\tau_t(Z_t)^2}{\vert\alpha\vert}$.
By \eqref{f39} we have
\begin{eqnarray*}
dH\left(Z_t\right)^n&=&nH\left(Z_t\right)^{n-1} \La_t(H)(Z_t)~dt\\
&&+\frac{n(n-1)}{2}~H\left(Z_t\right)^{n-2} 
\partial_t\left\langle
\Ma(H)
\right\rangle_t~dt+n H\left(Z_t\right)^{n-1} ~d\Ma_{t}(H).
\end{eqnarray*}
This implies that
\begin{eqnarray*}
\displaystyle\partial_t\EE\left(H\left(Z_t\right)^n\right)&
\displaystyle\leq& n~\alpha~\EE\left(H\left(Z_t\right)^{n}\right)\\
&&\displaystyle+n~\left[\EE\left(\beta_t\left(Z_t\right)H\left(Z_t\right)^{n-1}\right)+\frac{(n-1)}{2}~\EE\left(
\gamma_{t}(Z_t)H\left(Z_t\right)^{n-1} \right)\right].
\end{eqnarray*}
Using H\"older inequality we have
\begin{eqnarray*}
\EE\left(\beta_t\left(Z_t\right)H\left(Z_t\right)^{n-1}\right)&\leq &
\EE\left(H\left(Z_t\right)^{n}\right)^{1-\frac{1}{n}}
\EE\left(\beta_t^{n}\left(Z_t\right)\right)^{\frac{1}{n}}\leq \delta_{\beta,t}(n)~\EE\left(H\left(Z_t\right)^{n}\right)^{1-\frac{1}{n}}.
\end{eqnarray*}
In much the same way, we have
\begin{eqnarray*}
\EE\left(\gamma_t\left(Z_t\right)
H\left(Z_t\right)^{n-1}\right)
&\leq& \delta_{\gamma,t}(n)~\EE\left(H\left(Z_t\right)^{n}\right)^{1-\frac{1}{n}}.
\end{eqnarray*}
This yields the estimate
$$
\begin{array}{l}
\displaystyle\partial\EE\left(H\left(Z_t\right)^n\right)
\displaystyle\leq n~\alpha~\EE\left(H\left(Z_t\right)^{n}\right)+n~\delta_t(n)~\EE\left(H\left(Z_t\right)^{n}\right)^{1-\frac{1}{n}}
\end{array}
$$
with
$
\delta_t(n)=\delta_{\beta,t}(n)~
+(n-1)~\delta_{\gamma,t}(n)/2
$.
We conclude that
\begin{eqnarray*}
\displaystyle\partial_t\EE\left(H\left(Z_t\right)^n\right)^{\frac{1}{n}}&=&
\frac{1}{n}~\EE\left(H\left(Z_t\right)^n\right)^{\frac{1}{n}-1}\partial_t\EE\left(H\left(Z_t\right)^n\right)\leq  \alpha~\EE\left(H\left(Z_t\right)^{n}\right)^{\frac{1}{n}}+~\delta_t(n).
\end{eqnarray*}
The last assertion is a direct application of Gr\"onwall inequality.

The proof of the lemma
is now completed
\cqfd

\begin{prop}\label{prop-1}
Assume that $\mu(A)<0$. In this situation we have the uniform 
trace moment estimates
$$
1\leq n<1+\frac{N-1}{2r_1}~\frac{\lambda_{\tiny min}(S)}{\lambda_{\tiny max}(S)}~
\Longrightarrow~
\sup_{t\geq 0}
\EE\left(\left[\mbox{\rm tr}(p_t)\right]^n\right)\leq c(n)
$$
with the convention $\frac{\lambda_{\tiny min}(S)}{\lambda_{\tiny max}(S)}=0$ when $S=0$.
In addition, when condition (S) is met we have
$$
\sup_{t\geq 0}\EE\left(\left\Vert \zeta^1_t\right\Vert^{2n}\right)\leq c(n)
\quad\mbox{and}\quad
\sup_{t\geq 0}\EE\left(\left\Vert \xi^1_t\right\Vert^{2n}\right)\leq c(n).
$$
The l.h.s. estimates are valid for any $n\geq 1$, while
the r.h.s. ones are valid for $ 3n<1+\frac{N-1}{2r_1}\frac{\lambda_{\tiny min}(S)}{\lambda_{\tiny max}(S)}$.
\end{prop}
\proof

We set
$
H(p_t):=\mbox{\rm tr}(p_t)
$ the trace function of the random sample covariance matrices $p_t$. 
Using \eqref{f23} we prove evolution equation
$$
dH(p_t)=L(H)(p_t)~dt+\frac{1}{\sqrt{N-1}}~d\Ma_t
\quad\mbox{\rm with}\quad\Ma_t:=\mbox{\rm tr}(M_t)
$$
and the drift 
$$
L(H)(p_t):=2~\mbox{\rm tr}\left(A_{\tiny sym}p_t\right)-\mbox{\rm tr}\left(Sp^2_t\right)+\mbox{\rm tr}\left(R\right).
$$
Following the proof of Lemma~\ref{Lem-v1-1} we also have the estimates
\begin{eqnarray*}
0\leq \partial_t\left\langle \Ma\right\rangle_t&=&4~\mbox{\rm tr}\left((R+p_tSp_t)p_t\right)\leq
4 ~H(p_t)~\left(\mu(R)+\mu(S)~H(p_t)^2\right)
\end{eqnarray*}
and
$$
L(H)(p_t)~\leq 2~\alpha~H(p_t)-\beta~H(p_t)^2+r
$$
with 
$$
\alpha:=2\mu(A)<0\,,\qquad\beta:=r_1^{-1}\lambda_{\tiny min}(S)\quad\mbox{\rm and}\quad
r:=\mbox{\rm tr}(R).
$$
The moment estimates of the trace of the sample covariance matrices are now easily checked using Lemma~\ref{tech-1}. 

Now we come to the proof of the moments estimates of the norm of the samples. 
Notice that
$$
d\xi^i_t=\left((A-p_tS)~\xi^i_t+a+p_tSX_t\right)dt+dM^{\xi,i}_t
$$
with $r_1$-dimensional martingale
 $$
dM^{\xi,i}_t:=R^{1/2}d\overline{W}_t^i+p_tC^{\prime}R^{-1/2}_{2}d\left(V_t-\overline{V}^i_{t}\right).
$$
We set
$
H(X_t,\xi_t):=\left\Vert\xi^1_t\right\Vert^2
$. In this notation we have
$$
dH(X_t,\xi_t)
=\La_t(H)(X_t,\xi_t)~dt+d\Ma_t(H)\quad\mbox{\rm with}\quad
d\Ma_t(H)=2\langle\xi^1_t,dM^{\xi,1}_t\rangle
$$
and
$$
\La_t(H)(X_t,\xi_t)=\langle \xi^1_t,(A+A^{\prime}-(p_tS+Sp_t))~\xi^1_t\rangle+2\langle \xi^1_t,a+p_tSX_t\rangle+\mbox{\rm tr}\left(R+2~p_tSCp_t\right).
$$
On the other hand
$$
(S)\Longrightarrow
\langle \xi^1_t,(A+A^{\prime}-(p_tS+Sp_t))~\xi^1_t\rangle\leq 2\mu(A)\left\Vert\xi^1_t\right\Vert^2.
$$
Using the inequality
\begin{equation}\label{f42}
2\langle x,y\rangle\leq \frac{1}{\epsilon}~\Vert x\Vert^2+{\epsilon}~\Vert y\Vert^2
\end{equation}
with $0<\epsilon=-\mu(A)<-2\mu(A)$ (recall that $\mu(A)<0$) $y=\xi^i_t$ and $x=a+p_tSX_t$, we find
that
$$
\La_t(H)(X_t,\xi_t)\leq 
\mu(A)~\left\Vert\xi^i_t\right\Vert^2+\beta_t(X_t,\xi_t)=\mu(A)~H(X_t,\xi_t)+\beta_t(X_t,\xi_t)$$
with
$$
\beta_t(X_t,\xi_t)={\left\vert \mu(A)\right\vert^{-1}}~\Vert a+p_tSX_t\Vert^2+\mbox{\rm tr}\left(R+2~p_tSCp_t\right).
$$
On the other hand we have
$$
\Vert a+p_tSX_t\Vert\leq \Vert a\Vert+\Vert S\Vert~\Vert X_t\Vert ~\Vert p_t\Vert\quad \mbox{\rm with}\quad
\Vert p_t\Vert\leq \Vert p_t\Vert_F\leq \mbox{\rm tr}(p_t)
$$
and
$$
\mbox{\rm tr}\left((R+p_tSp_t)p_t\right)\leq
\mbox{\rm tr}(p_t)~\left(\mu(R)+\mu(S)~\left(\mbox{\rm tr}(p_t)\right)^2\right).
$$
Combining these estimates with the trace estimates we have just proven and the signal state uniform moment 
estimates stated in \eqref{fv1-4} a direct
application of the Lemma~\ref{tech-1} yields 
$$\sup_{t\geq 0}\EE\left(\Vert \xi^1_t\Vert^{2n}\right)\leq c(n).$$ Notice that the control of the $n$-th
moment of $\beta_t(X_t,\xi_t)$ involves the control   of the $(3n)$-th
moment of the trace of $p_t$. 
The same analysis applies to $\zeta^1_t$ with $p_t$ replaced by $P_t$.

This ends the proof of the proposition.\cqfd

\section{Quantitative properties}\label{quantitative-section}

This section is mainly concerned with the proof of the uniform estimates presented in 
Theorem~\ref{theo-intro}. 
The first step is to control and to estimate the fluctuations of the particle covariance matrices
involved in the EnKF filter uniformly w.r.t. the time horizon. In Section~\ref{sample-cov-sec}
we present a key uniform control of the Frobenius norm between the particle covariance matrices 
and their limiting values. These estimates are used in Section~\ref{upc-sec} to 
derive the uniform propagation of chaos properties of the EnKF.

\subsection{Particle covariance matrices}\label{sample-cov-sec}
Next theorem is pivotal.
It describes the evolution of the Frobenius norm of the ``centered'' sample covariance matrices
in terms of a nonlinear diffusion and provides some key uniform convergence results. 

\begin{theo}\label{theo3}
The Frobenius norm of the  sample covariance matrix fluctuations satisfies the diffusion equation
\begin{equation}\label{f38}
\displaystyle d\Vert p_t-P_t\Vert_{F}^2
\displaystyle=\left(\alpha_t(p_t)+\frac{1}{N-1}~\beta_t(p_t)\right)~dt +
\frac{2}{\sqrt{N-1}}~d\Ma_t
\end{equation}
with the drift functions
\begin{eqnarray*}
\alpha_t(p_t)&:=&2~\mbox{\rm tr}\left(\left[(A+A^{\prime})-\frac{1}{2}\left\{(p_t+P_t)S+S(p_t+P_t)\right\}\right](p_t-P_t)^2\right)\\
\beta_t(p_t)&:=&2\left[\gamma_t(p_t)+\mbox{\rm tr}\left(R+p_tSp_t\right)\mbox{\rm tr}\left(p_t\right)\right]\qquad
\gamma_t(p_t):=\mbox{\rm tr}\left((R+p_tSp_t)p_t\right)
\end{eqnarray*}
and a martingale $\Ma_t$ with angle bracket 
$$
\partial_t\left\langle \Ma\right\rangle_t=4\mbox{\rm tr}(p_t(p_t-P_t)(R+p_tSp_t)(p_t-P_t)).
$$
In addition, when $\mu(A)<0$  and  (S) is satisfied we have the uniform mean error estimates
$$
2n<1+\frac{N-1}{2r_1}~\Longrightarrow~
\sup_{t\geq 0}{\EE\left(\Vert p_t-P_t\Vert_{F}^{n}\right)^{\frac{1}{n}}}\leq c(n)/\sqrt{N}
$$

\end{theo}

\proof

By \eqref{nonlinear-KB-Riccati} and \eqref{f23} we have

$$
d(p_t-P_t)=\left[A(p_t-P_t)+(p_t-P_t)A^{\prime}-\left(p_tSp_t-P_tSP_t\right)
\right]dt+\frac{1}{\sqrt{N-1}}~dM_t
$$
with the martingale $M_t$
with  angle brackets defined in \eqref{f24}.
Using the decomposition
$$
p_tSp_t-P_tSP_t=(p_t-P_t)S(p_t-P_t)+P_tS(p_t-P_t)+(p_t-P_t)SP_t,
$$
we readily check that
$$
\begin{array}{l}
A(p_t-P_t)+(p_t-P_t)A^{\prime}-\left(p_tSp_t-P_tSP_t\right)\\
\\
=\left[A-\frac{1}{2}(p_t+P_t)S\right](p_t-P_t)+(p_t-P_t)\left[A^{\prime}-\frac{1}{2}S(p_t+P_t)\right].
\end{array}
$$
This implies that
$$
\begin{array}{l}
\displaystyle d\left((p_t-P_t)^2\right)\\
\\
\displaystyle=2~(p_t-P_t)~\left\{\left[A-\frac{1}{2}(p_t+P_t)S\right](p_t-P_t)+(p_t-P_t)\left[A^{\prime}-\frac{1}{2}S(p_t+P_t)\right]\right\}~dt\\
\\
\displaystyle+\frac{1}{N-1}~\left[(R+p_tSp_t)p_t+p_t(R+p_tSp_t)+\mbox{\rm tr}(p_t)~(R+p_tSp_t)+
\left(\mbox{\rm tr}(R+p_tSp_t)\right)p_t\right]~dt\\
\\
\displaystyle\hskip6cm +
\frac{2}{\sqrt{N-1}}~(p_t-P_t)~dM_t.
\end{array}
$$
Taking the trace we find that
$$
\begin{array}{l}
\displaystyle d\Vert (p_t-P_t)\Vert_{F}^2\\
\\
\displaystyle=2~\mbox{\rm tr}\left(\left[(A+A^{\prime})-\frac{1}{2}\left\{(p_t+P_t)S+S(p_t+P_t)\right\}\right](p_t-P_t)^2\right)  dt\\
\\
\displaystyle+\frac{2}{N-1}~\left[\mbox{\rm tr}\left((R+p_tSp_t)p\right)+\mbox{\rm tr}\left((R+p_tSp_t)\right)\mbox{\rm tr}\left(p_t\right)\right]~dt +
\frac{2}{\sqrt{N-1}}~d\Ma_t
\end{array}
$$
with the martingale
\begin{eqnarray*}
d\Ma_t&=&\mbox{\rm tr}\left((p_t-P_t)~dM_t\right)=\sum_{1\leq k,l\leq r_1}(p_t-P_t)(l,k)~dM_t(k,l).
\end{eqnarray*}
The angle bracket of $\Ma_t$ is computed using \eqref{f24}. More precisely we have
$$
\partial_t\left\langle \Ma\right\rangle_t=4\mbox{\rm tr}\left((R+p_tSp_t)(p_t-P_t)p_t(p_t-P_t)\right).
$$
We check this claim using the decomposition
$$
\partial_t\left\langle \Ma\right\rangle_t=\sum_{1\leq k,l,k^{\prime},l^{\prime}\leq r_1}(p_t-P_t)(l,k)(p_t-P_t)(l^{\prime},k^{\prime})~\partial_t\langle M(k,l),M(k^{\prime},l^{\prime})\rangle_t.
$$
Recalling that
$$
\begin{array}{rcl}
\displaystyle\partial_t\left\langle  M(k,l), M(k^{\prime},l^{\prime})\right\rangle_t&
=&\left(R+p_tSp_t\right)(k,k^{\prime})~ p_t(l,l^{\prime})+
\left(R+p_tSp_t\right)(l,l^{\prime})~ p_t(k,k^{\prime})\\
&&\\
&&\displaystyle+
\left(R+p_tSp_t\right)(l^{\prime},k) ~p_t(k^{\prime},l)+
\left(R+p_tSp_t\right)(l,k^{\prime}) ~p_t(k,l^{\prime}).
\end{array}
$$
and using the symmetry of the matrices $(p_t,P_t)$ and $\left(R+p_tSp_t\right)$ we find that
\begin{eqnarray*}
\partial_t\left\langle \Ma\right\rangle_t&
\displaystyle=&\sum_{1\leq k,l,k^{\prime},l^{\prime}\leq r_1}
p_t(l^{\prime},l)~(p_t-P_t)(l,k)~\left(R+p_tSp_t\right)(k,k^{\prime})~(p_t-P_t)(k^{\prime},l^{\prime}) \\
&&\displaystyle+\sum_{1\leq k,l,k^{\prime},l^{\prime}\leq r_1}
(p_t-P_t)(k,l)~\left(R+p_tSp_t\right)(l,l^{\prime})~(p_t-P_t)(l^{\prime},k^{\prime})~p_t(k^{\prime},k)~\\
&&
\displaystyle+\sum_{1\leq k,l,k^{\prime},l^{\prime}\leq r_1}
p_t(k^{\prime},l)(p_t-P_t)(l,k)~\left(R+p_tSp_t\right)(k,l^{\prime}) ~(p_t-P_t)(l^{\prime},k^{\prime})\\
&&
\displaystyle+\sum_{1\leq k,l,k^{\prime},l^{\prime}\leq r_1}
(p_t-P_t)(l,k)p_t(k,l^{\prime})~(p_t-P_t)(l^{\prime},k^{\prime})~\left(R+p_tSp_t\right)(k^{\prime},l) 
\end{eqnarray*}
This shows that
\begin{eqnarray*}
\partial_t\left\langle \Ma\right\rangle_t&
\displaystyle=&4\mbox{\rm tr}(p_t(p_t-P_t)(R+p_tSp_t)(p_t-P_t))
\end{eqnarray*}

This ends the proof of the first assertion.

We set $H(p_t,P_t)=\Vert p_t-P_t\Vert_{F}^2$. In this notation, Equation \eqref{f38} takes the form
$$
\displaystyle\partial_tH(p_t,P_t)
\displaystyle=\La_t(H)(p_t,P_t)~dt + d\Ma_t(H)
\quad\mbox{\rm with}\quad \Ma_t(H):=\frac{2}{\sqrt{N-1}}~\Ma_t
$$
and
$$
\La_t(H)(p_t,P_t)=\alpha_t(p_t)+\frac{1}{N-1}~\beta_t(p_t).
$$

Under condition (S) we have
$$
\begin{array}{l}
\mbox{\rm tr}\left(\left[(A+A^{\prime})-\frac{1}{2}\left\{(p_t+P_t)S+S(p_t+P_t)\right\}\right](p_t-P_t)^2\right)\\
\\
=\mbox{\rm tr}\left((A+A^{\prime})(p_t-P_t)^2\right)-\rho(S)~\mbox{\rm tr}\left((p_t+P_t)(p_t-P_t)^2\right)\
\end{array}
$$
Using \eqref{f30} this implies that
\begin{eqnarray*}
\alpha_t(p_t)&\leq& 2
\lambda_{\tiny max}\left(A+A^{\prime}
\right)~\Vert p_t-P_t\Vert_{F}^2\leq 4\mu(A)~H(p_t,P_t)
\end{eqnarray*}
and by \eqref{form-pq-ref} we have
$$
\partial_t\left\langle \Ma(H)\right\rangle_t\leq {4^2}{(N-1)^{-1}}~\mbox{\rm tr}\left(R+p_tSp_t\right)\mbox{\rm tr}\left(p_t\right)~H(p_t,P_t).
$$
This implies that
$$
\La_t(H)(p_t,P_t)\leq  4\mu(A)~H(p_t,P_t)+{(N-1)^{-1}}~\beta_t(p_t).
$$
Arguing as in the end of the proof of Proposition~\ref{prop-1} we also have
$$
{\sup_{t\geq 0}{\EE\left(\beta_t(p_t)^n\right)}}\leq c(n)
\quad\mbox{and}\quad
{\sup_{t\geq 0}{\EE\left(\gamma_t(p_t)^n\right)}}\leq c^{\prime}(n)
$$
for any $n$ s.t. $ 3n<1+\frac{N-1}{2r_1}$.
Using Lemma~\ref{tech-1} we conclude that
$$
\EE\left(\Vert p_t-P_t\Vert_{F}^{2n}\right)^{\frac{1}{n}}\leq c(n)/N.
$$
For odd numbers, we use H\"older inequality
$$
\EE\left(\Vert p_t-P_t\Vert_{F}^{2n+1}\right)\leq \EE\left(\Vert p_t-P_t\Vert_{F}^{2n}\right)^{\frac{1}{2}}
\EE\left(\Vert p_t-P_t\Vert_{F}^{2(n+1)}\right)^{\frac{1}{2}}.
$$
The end of the proof of the uniform estimates is now easily completed.
The proof of the theorem is now completed.
\cqfd

\subsection{Uniform propagation of chaos}\label{upc-sec}

This section is mainly concerned with the proof of the uniform estimates in \eqref{upxi}.
We set 
$$
Z_t=(P_t,X_t,\zeta_t,\xi_t)\qquad
\chi_t:=\xi^1_t-\zeta^1_t\qquad H(Z_t):=\left\Vert\chi_t\right\Vert^2\quad \mbox{\rm and}\quad q_t:=p_t-P_t.
$$
We have
$$
d\chi_t=\left(A\chi_t-\left(p_tS~\xi^1_t-P_tS\zeta^1_t\right)
+q_tSX_t\right)dt+q_tC^{\prime}R^{-1/2}_{2}d\left(V_t-\overline{V}^1_{t}\right).
$$
Observe that
$$
p_tS\xi^1_t-P_tS\zeta^1_t=p_tS\chi_t+q_tS\zeta^1_t.
$$
This implies that
\begin{equation}\label{fv1-5}
d\chi_t=\left((A-p_tS)\chi_t-q_tS\left(\zeta^1_t-X_t\right)\right)dt+dM_t
\end{equation}
with
$$
dM_t=q_tC^{\prime}R^{-1/2}_{2}d\left(V_t-\overline{V}^1_{t}\right).
$$
The angle bracket matrix $\langle M\rangle_t=\left(\langle M(k),M(l)\rangle_t\right)_{1\leq k,l\leq r_1}$ is given by
the formula
\begin{eqnarray}\label{ab-precision}
\partial_t\langle M\rangle_t=2q_tSq_t.
\end{eqnarray}
To check this claim observe that
$$
dM_t(k)=\sum_{l}\left(q_tC^{\prime}R^{-1/2}_{2}\right)(k,l)~d\left(V_t-\overline{V}^1_{t}\right)(l)
$$
This implies that
\begin{eqnarray*}
\partial _t\langle M(k),M(k^{\prime})\rangle_t&=&2\sum_{l}\left(q_tC^{\prime}R^{-1/2}_{2}\right)(k,l)
\left(q_tC^{\prime}R^{-1/2}_{2}\right)(k^{\prime},l)\\
&=&2\sum_{l}\left(q_tC^{\prime}R^{-1/2}_{2}\right)(k,l)
\left(R^{-1/2}_{2}Cq_t\right)(l,k^{\prime})=2(q_tSq_t)(k,k^{\prime}).
\end{eqnarray*}
This ends the proof of (\ref{ab-precision}).
On the other hand we have
\begin{eqnarray*}
d\left\Vert\chi_t\right\Vert^2&=&2\langle \chi_t,d\chi_t\rangle+\mbox{\rm tr}(\langle M\rangle_t)~dt\\
&=&2\left[\langle\chi_t,(A-p_tS)\chi_t\rangle-
\langle \chi_t,q_tS\left(\zeta^1_t-X_t\right)\rangle+\mbox{\rm tr}\left(q_t^2S\right)\right]+2\langle \chi_t,dM_t\rangle
\end{eqnarray*}
This yields 
$$
dH(Z_t)
=\La_t\left(H\right)(Z_t)~dt+d\Ma_t(H)\quad \mbox{\rm with}\quad d\Ma_t(H):=2\langle \chi_t,dM_t\rangle
$$
with
$$
\La_t\left(H\right)(Z_t)=2\left[\langle\chi_t,(A-p_tS)\chi_t\rangle-
\langle \chi_t,q_tS\left(\zeta^1_t-X_t\right)\rangle+\mbox{\rm tr}\left(q_t^2S\right)\right].
$$
We also have
$$
\begin{array}{l}
d\Ma_t(H)=2\sum_{k} \chi_t(k)~dM_t(k)\\
\\
\Rightarrow
\begin{array}[t]{rcl}
\partial_t\langle \Ma(H)\rangle_t&=&4\sum_{k,k^{\prime}} \chi_t(k)\chi_t(k^{\prime})~\partial_t\langle M(k),M(k^{\prime})\rangle_t\\
&&\\
&=&8~\sum_{k,k^{\prime}} \chi_t(k)~(q_tSq_t)(k,k^{\prime})~\chi_t(k^{\prime})
\end{array}
\end{array}
$$
This implies that
$$
\partial_t\langle \Ma(H)\rangle_t:=8~\langle\chi_t,q_tSq_t\chi_t\rangle.
$$

Using \eqref{new-condition}, \eqref{f30} and \eqref{f42} we check that
$$
\La_t\left(H\right)(Z_t)
\displaystyle\leq \mu(A)\left\Vert\chi_t\right\Vert^2+{\left\vert \mu(A)\right\vert^{-1}}~\left\Vert q_tS\left(\zeta^1_t-X_t\right)\right\Vert^2+2\mu(S)~\left\Vert q_t\right\Vert^2_F.
$$
On the other hand, we have
$$
\left\Vert q_tS\left(\zeta^1_t-X_t\right)\right\Vert\leq 
\left\Vert q_t\right\Vert_F\left\Vert S\right\Vert\left\Vert\zeta^1_t-X_t\right\Vert.
$$
This implies that
$$
\La_t\left(H\right)(Z_t)\leq \mu(A)~H(Z_t)+\beta_t(Z_t)
$$
with
$$
\beta_t(Z_t)=\left[{\left\vert \mu(A)\right\vert^{-1}}~\left\Vert S\right\Vert^2\left\Vert\zeta^1_t-X_t\right\Vert^2+2\mu(S)\right]~\left\Vert q_t\right\Vert^2_F.
$$

In much the same way we have
$$
\partial_t\langle \Ma(H)\rangle_t\leq~H(Z_t)~ \gamma_t(Z_t)\quad
\mbox{\rm with}\quad
\gamma_t(Z_t):=8\left\Vert S\right\Vert~
\left\Vert q_t\right\Vert^2_F.
$$
By Theorem~\ref{theo3}, and using the uniform estimates stated in \eqref{fv1-4} and in Proposition~\ref{prop-1}
for any $2n<1+\frac{N-1}{2r_1}$ we have 
$$
\sup_{t\geq 0}{\EE\left(\beta_t(Z_t)^{n}\right)^{\frac{1}{n}}}\leq c(n)/\sqrt{N}
\quad\mbox{\rm and}\quad
\sup_{t\geq 0}{\EE\left(\gamma_t(Z_t)^{n}\right)^{\frac{1}{n}}}\leq c(n)/\sqrt{N}.
$$
The proof of the uniform estimates in \eqref{upxi} is now a direct consequence of Lemma~\ref{tech-1}. This ends the proof of Theorem~\ref{theo-intro}.\cqfd

\section{Proof of proposition~\ref{prop-stability-1}}\label{Lip-global-sec}

The estimate (\ref{f-easy-intro}) is a direct consequence of the perturbation lemma~\ref{perturbation-lemma-intro}
and the triangle inequality (\ref{triangle-ref}).

Now we come to the proof of (\ref{prop-stability-1-2}). We set
$$
\Xa_t=\breve{X}_t-\overline{X}_t
\quad\mbox{\rm and}\quad
Q_t=\widecheck{P}_t-P_t.
$$
Arguing as in \eqref{fv1-5} we have
\begin{equation}\label{fin-eq-1}
d\Xa_t=\left((A-\widecheck{P}_tS)\Xa_t-Q_tS\left(\overline{X}_t-X_t\right)\right)dt+dM_t
\quad\mbox{\rm
with}\quad\partial_t\langle M\rangle_t=Q_tSQ_t.
\end{equation}
This implies that
$$
d\Vert \Xa_t\Vert^2=2\left[\langle\Xa_t,(A-\widecheck{P}_tS)\Xa_t\rangle-
\langle \Xa_t,Q_tS\left(\overline{X}_t-X_t\right)\rangle+\mbox{\rm tr}\left(Q_t^2S\right)/2\right]+d\Ma_t
$$
with a real valued martingale with angle bracket 
$$
\partial_t\langle \Ma\rangle_t={4}~\langle\Xa_t,Q_tSQ_t\Xa_t\rangle\leq {4}\Vert S\Vert~\Vert \Xa_t\Vert^2~\Vert Q_t\Vert^2.
$$
Using the same arguments as in the proof of \eqref{upxi} given in the end of Section~\ref{quantitative-section}
we conclude that
\begin{eqnarray*}
\EE\left(\Vert \Xa_t\Vert^{2n}\right)^{\frac{1}{n}}&\leq&
e^{\mu(A)t}~\EE\left(\Vert \Xa_0\Vert^{2n}\right)^{\frac{1}{n}}+~c~\int_0^t~e^{\mu(A)(t-s)}~\Vert Q_s\Vert_F^2~ds
\end{eqnarray*}
for any $n$ and some finite constant $c$.  Using \eqref{f-easy-intro} we arrive at the estimate 
\begin{eqnarray*}
\EE\left(\Vert \overline{X}_t-\breve{X}_t\Vert^{2n}\right)^{\frac{1}{2n}}&\leq&
e^{\mu(A)t/2}~\EE\left(\Vert \overline{X}_0-\breve{X}_0\Vert^{2n}\right)^{\frac{1}{2n}}\\
&&\hskip3cm+c~\left[
\int_0^t~e^{\mu(A)(t-s)}~ e^{2\mu(A)s}~ds\right]^{1/2}~ \Vert P_0-\widecheck{P}_0\Vert_{F}\\
&\leq &e^{\mu(A)t/2}~\left[\EE\left(\Vert \overline{X}_0-\breve{X}_0\Vert^{2n}\right)^{\frac{1}{2n}}+c^{\prime}~
 \Vert P_0-\widecheck{P}_0\Vert_{F}
\right]
\end{eqnarray*}
for some finite constant $c^{\prime}$.
Also notice that
$$
\begin{array}{l}
\eqref{f-easy-intro-prop-mean}\quad\mbox{\rm and}\quad
\eqref{prop-stability-1-2}~\\
\\
\Rightarrow~
\EE\left(\Vert \widehat{X}_t-\widecheck{X}_t\Vert^{2n}\right)^{\frac{1}{2n}}\leq
e^{\mu(A)t/2}~\left[\EE\left(\Vert \widehat{X}_0-\widecheck{X}_0\Vert^{2n}\right)^{\frac{1}{2n}}+c~
 \Vert P_0-\widecheck{P}_0\Vert_{F}
\right].
\end{array}$$

This ends the proof of proposition~\ref{prop-stability-1}.\cqfd

\section{Proof of theorem~\ref{prop-stab-Kalman}}\label{stab-kalman-sec}

 We have
$$
d(\overline{X}_t-X_t)=(A-P_tS)(\overline{X}_t-X_t)+R_1^{1/2}(d\overline{W}_t-W_t)+P_t~C^{\prime}R_2^{-1/2}(dV_t-\overline{V}_t)
$$
This implies that
$$
d\Vert\overline{X}_t-X_t\Vert^2=\left[2\langle \overline{X}_t-X_t, (A-P_tS)(\overline{X}_t-X_t)\rangle+\mbox{\rm tr}(R_1+P_tSP_t)\right]~dt+
dM_t
$$
with some martingale $M_t$ s.t. 
$$
\partial_t\langle M_t\rangle\leq \Vert\overline{X}_t-X_t\Vert^2~\mbox{\rm tr}(R_1+P_tSP_t)
$$
We set $A_P:=A-PS$.
For any $u\in ]0,1]$  there exists some time horizon
 $\tau_{u}(P_0)\geq 0$ such that for any $t\geq \tau_{u}(P_0) $ 
 $$
 \sup_{t\geq \tau_u(P_0)}{\mu(A-P_tS)}\leq (1-u)~\mu(A_P)
 $$
 This yields
 $$
d\Vert\overline{X}_t-X_t\Vert^2=\left[2(1-u)~\mu(A_P)~\Vert\overline{X}_t-X_t\Vert^2+c(u)\right]~dt+
dM_t
$$
and
$$ \partial_t\langle M_t\rangle\leq \Vert\overline{X}_t-X_t\Vert^2~c(u)
$$
 for some finite constant $c(u)$ whose values only depend on $u$. By lemma~\ref{tech-1} we conclude that
 $$
 \sup_{t\geq 0}\EE\left(\Vert\overline{X}_t-X_t\Vert^n\right)<\infty
 $$
for any $n\geq 1$.

When $\breve{X}_t$ is the steady state Kalman-Bucy diffusion we have $\widecheck{P}_t=P$, for any $t\geq 0$.
In this situation (\ref{fin-eq-1}) takes the form
$$
d\Xa_t=\left(A_P\Xa_t-Q_tS\left(\overline{X}_t-X_t\right)\right)~dt+dM_t
$$
with $Q_t=P-P_t$ and a martingale $M_t$ with
$\partial_t\langle M\rangle_t=Q_tSQ_t
$.
This implies that
$$
d\Vert \Xa_t\Vert^2=2\left[\langle\Xa_t,A_P\Xa_t\rangle-
\langle \Xa_t,Q_tS\left(\overline{X}_t-X_t\right)\rangle+\mbox{\rm tr}\left(Q_t^2S\right)/2\right]+d\Ma_t
$$
with a real valued martingale with angle bracket 
$$
\partial_t\langle \Ma\rangle_t={4}~\langle\Xa_t,Q_tSQ_t\Xa_t\rangle\leq {4}\Vert S\Vert~\Vert \Xa_t\Vert^2~\Vert Q_t\Vert^2.
$$
Using the log-norm Lipschitz estimate (\ref{with-mu}) in  corollary~\ref{cor-bucy-entended} we have
\begin{eqnarray*}
\Vert Q_t S\Vert_F\leq \Vert S\Vert_F~\Vert Q_t \Vert_F
&\leq& c(u)~\exp{\left(2(1-u)\mu(A_P)t\right)}~\Vert S\Vert_F~\Vert P_0-P \Vert_F
\\
\mbox{\rm tr}\left(Q_t^2S\right)\vee
\Vert Q_t SQ_t\Vert_F&\leq& c(u)^2~\exp{\left(4(1-u)\mu(A_P)t\right)}~\Vert S\Vert_F~\Vert P_0-P \Vert_F^2
\end{eqnarray*}
 for any $u\in ]0,1]$, and $t\geq 0$, and  for some constant $c(u)$ that depends on $u$. 
 
 We set
 $$
 Z_t=(\Xa_t,P_t)\qquad
 H(Z_t)=\Vert \Xa_t\Vert^2\quad\mbox{\rm and}\quad M_t(H):=M_t
 $$
and we consider the parameter $ \alpha:=\mu(A_P) $ and the function 
 \begin{eqnarray*}
 \tau_t(Z_t)&:=&c(u)~\exp{\left(2(1-u) \alpha t\right)}~\Vert S\Vert_F\Vert P_0-P \Vert_F
 \Vert\overline{X}_t-X_t\Vert\\
  \beta_t(Z_t)&:=&c(u)^2~\exp{\left(4(1-u) \alpha t\right)}~\Vert S\Vert_F\Vert P_0-P \Vert_F^2\\
  \gamma_t(Z_t)&:=&({4}\Vert S\Vert)~c(u)^2~\exp{\left(4(1-u) \alpha t\right)}\Vert P_0-P \Vert_F^2
 \end{eqnarray*}
 In this notation we have
 $$
 dH(Z_t)\leq {\cal L}_t(H)(Z_t)~dt+M_t(H)
 $$
 with
 $$
 {\cal L}_t(H)(Z_t)\leq 2\alpha~H(Z_t)+ 2~ \tau_t(Z_t)~\sqrt{H(Z_t)}+\beta_t(Z_t)\quad\mbox{\rm and}\quad
  \partial_t\langle M(H)\rangle_t\leq  ~H(Z_t)~ \gamma_t(Z_t).
 $$
We also have
 \begin{eqnarray*}
\EE\left(\tau_{t}\left(Z_t\right)^n\right)^{\frac{1}{n}}&\leq& c_1(u)~\exp{\left(2(1-u)\alpha t\right)}\Vert P_0-P \Vert_F\\
\EE\left(\beta_{t}\left(Z_t\right)^n\right)^{\frac{1}{n}}\vee\EE\left(\gamma_{t}\left(Z_t\right)^n\right)^{\frac{1}{n}}&\leq& c_2(u)~\exp{\left(4(1-u)\alpha t\right)}\Vert P_0-P \Vert_F^2
 \end{eqnarray*}
for any $n\geq 1$, for some constants $c_i(u)$ that depends on $u$, with $i\in\{1,2\}$. 
 Using lemma~\ref{tech-1} we check that for any $u\in ]0,1]$ and any $t\geq 0$ we have
 $$
\displaystyle\EE\left(\Vert \Xa_t\Vert^{2n}\right)^{\frac{1}{n}}
\displaystyle\leq ~e^{\alpha t}~
\EE\left(\Vert \Xa_0\Vert^{2n}\right)^{\frac{1}{n}} + c(u)~e^{(1-u)\alpha t}\Vert P_0-P \Vert_F^2.
$$
 for some constants $c(u)$ that depends on $u$. Using \eqref{f-easy-intro-prop-mean}, we also check that
$$
\begin{array}{l}
\EE\left[\Vert \widehat{X}_t-\widecheck{X}_t\Vert^{2n}\right]^{\frac{1}{2n}}\leq c(u)~e^{(1-u)\alpha t/2}\left(\EE\left(\Vert
\widehat{X}_0-\widecheck{X}_0
\Vert^{2n}\right)^{\frac{1}{2n}}+\Vert P_0-P \Vert_F\right)
\end{array}
$$
where $\widecheck{X}_t$ is the steady state Kalman-Bucy filter. This ends the proof of the proposition.
\cqfd

\section{Proof of Theorem~\ref{theo-Ent}}\label{sect-last-entropy}
To simplify the presentation we assume that $\widecheck{P}_0=P$.

We further assume that the algebraic Riccati Equation
 \eqref{steady-state-eq} has a positive definite fixed point $P$ (so that $P$ is invertible). We also assume
 that  $\mu(A-PS)<0$. 
 
We let
 $\left(\overline{X}_t,\breve{X}_t\right)$ be a couple of Kalman-Bucy Diffusions \eqref{Kalman-Bucy-filter-nonlinear-ref} starting from
two possibly different {\em Gaussian} random variables with
 covariance matrices $\left(P_0,P\right)$. We recall that
 $$
\widehat{X}_t= \EE\left(\overline{X}_t~|~\Fa_t\right)\quad\mbox{\rm and}\quad \widecheck{X}_t= \EE\left(\breve{X}_t~|~\Fa_t\right)
 $$
satisfy the Kalman-Bucy Recursion \eqref{nonlinear-KB-mean}
 associated with the covariance matrices $(P_t,P)$, with $P_t$ given by the Riccati Equation \eqref{nonlinear-KB-Riccati}.
 
 Let $\left(\eta_t,\breve{\eta}_t\right)$ be the (Gaussian) conditional distributions of  $\left(\overline{X}_t,\breve{X}_t\right)$ 
 given the $\sigma$-field  $\Fa_t$ generated by the observation process.  The conditional Boltzmann-Kullback Liebler relative entropy of $\eta_t$ w.r.t. $\breve{\eta}_t$ is given by the formula
 $$
 \mbox{\rm Ent}\left(\eta_t~|~\breve{\eta}_t\right)=-\frac{1}{2}\left(
 \mbox{\rm tr}\left(I-P^{-1}P_t\right)+\log{\mbox{det}\left(P_tP^{-1}\right)}-\left\langle\left(\widecheck{X}_t-\widehat{X}_t\right),P^{-1}\left(\widecheck{X}_t-\widehat{X}_t\right)\right\rangle\right).
 $$
 To estimate the logarithm of the determinant of the matrices $P_tP^{-1}$ as $t\uparrow\infty$ we use the following technical lemma.
\begin{lem}\label{lem-tech-2}
For any $(r\times r)$-matrices $(A,B)$  we have
$$
\Vert A \Vert_F~\Vert B\Vert_F< \frac{1}{2}\Longrightarrow
\left\vert\log{\mbox{\rm det}\left(I-A\right)}\right\vert 
\leq \frac{3}{2}~\Vert A \Vert_F~\Vert B\Vert_F
$$

\end{lem}

\proof
For any $n\geq 1$ we have
$$
\vert\mbox{\rm tr}(A^nB^n)\vert\leq \Vert A \Vert^n_F~\Vert B\Vert^n_F
$$
Using the well-known trace formulae $$
\log{\mbox{\rm det}(I-AB)}=\mbox{\rm tr}(\log{(I-AB)})=-\sum_{n\geq 1}~n^{-1}~\mbox{\rm tr}(A^nB^n)
$$
 we conclude that
 $$
 \vert \log{\mbox{\rm det}(I-AB)}\vert\leq -\log{\left( 1-\Vert A \Vert_F~\Vert B\Vert_F\right)}
 $$
The last assertion comes from the inequality
$$
0\leq -\log{(1-u)}\leq u+\frac{1}{2}~\frac{u^2}{1-u}=u\left(1+\frac{1}{2}~\frac{u}{1-u}\right)\leq 3u/2
$$
which is valid for any $u\in [0,1/2[$.

This ends the proof of the lemma.\cqfd

For any  $u\in ]0,1]$, and $t\geq 0$, 
 there exists some $c(u)$ that depends on $u$ s.t.
 $$
\Vert I-P^{-1}P_t\Vert_F\leq \Vert P^{-1}\Vert_F~\Vert P-P_t\Vert_F\leq c(u)~\exp{\left(2(1-u)~\mu(A-PS)~t\right)}~\Vert P_0-P \Vert_F
 $$

Applying Lemma~\ref{lem-tech-2} to $A=(P-P_t)P^{-1}$ there exists some $t_0$ that depends on $(P_0,P)$ and some finite constant $c$ such that
  \begin{eqnarray*}
\left\vert\log{\mbox{det}\left(P_tP^{-1}\right)}\right\vert
&\leq &\frac{3}{2}~\Vert P-P_t \Vert_F~\Vert P^{-1}\Vert_F
\end{eqnarray*}
for any $t\geq t_0$.
On the other hand, using the monotocity properties of $P_t$ we have
  \begin{eqnarray*}
\left\vert\mbox{\rm tr}\left(I-P^{-1}P_t\right)\right\vert&=&\left\vert\mbox{\rm tr}\left(P^{-1}(P-P_t)\right)\right\vert\\
&\leq & \left\vert\mbox{\rm tr}\left(P^{-1}\right)\right\vert\left\vert\mbox{\rm tr}\left(P-P_t\right)\right\vert\leq  \left\vert\mbox{\rm tr}\left(P^{-1}\right)\right\vert~\left\Vert P-P_t\right\Vert_F.
\end{eqnarray*}
Finally, we notice that
$$
\left\vert
\left\langle\left(\widecheck{X}_t-\widehat{X}_t\right),P^{-1}\left(\widecheck{X}_t-\widehat{X}_t\right)\right\rangle
\right\vert\leq \left\Vert P^{-1}\right\Vert~ \Vert \widecheck{X}_t-\widehat{X}_t\Vert^2.
$$
We conclude that
$$
0\leq 
 \mbox{\rm Ent}\left(\eta_t~|~\breve{\eta}_t\right)\leq  \frac{3}{2}~\left\Vert P^{-1}\right\Vert_F\left( 
 \Vert \widecheck{X}_t-\widehat{X}_t\Vert^2
+\left\Vert P-P_t\right\Vert_F\right)
$$
for any $t\geq t_0$. The end of the proof of (\ref{last-formula}) is now clear. 
This ends the proof of theorem~\ref{theo-Ent}.

\cqfd

\section{Conclusion}\label{conclusion}

We have designed and analyzed a new class of conditional nonlinear diffusion processes arising in 
filtering theory. In contrast with conventional nonlinear Markov models,
these Kalman-Bucy diffusion type models depends on the conditional covariance matrices 
of the internal random states.  To
analyze the stability properties of these models, a series of functional contraction inequalities have been developed
w.r.t. the Wasserstein distance, 
Frobenius norms on random matrices and relative entropy criteria. 

In this framework, the traditional Kalman-Bucy filter resumes to the time evolution of the conditional averages of 
these nonlinear diffusions. The stability properties of the filter are now deduced directly from the ones of the 
nonlinear model.

The second important contribution of the article concerns the long-time behaviour and the refined convergence analysis of Ensemble Kalman filters.
The EnKF is interpreted as a natural mean-field particle approximation of nonlinear Kalman-Bucy diffusions.
The performance of the EnKF is measured in terms of uniform $\LL_n$-mean error estimates and uniform propagation of chaos properties w.r.t. the time horizon. 

We end this article of an avenue of open research problems. 

The first project is to extend the analysis to nonlinear diffusions with an interacting function that depends on the covariance
matrices of the random states. A toy model of that form is given by the one dimensional diffusion
$$
dX_t=\mbox{\rm Var}(X_t)~(dW_t-X_tdt)
$$
where $W_t$ stands for a Brownian motion. It is readily check that this nonlinear diffusion is well-posed. In addition, the variance
$P_t=\mbox{\rm Var}(X_t)=P_0/(1+P_0t)$ satisfies the Riccati equation
$$
\partial_t P_t=-P_t^2\quad\mbox{\rm and}\quad X_t={(X_0+P_0W_t)}/{(1+P_0~t)}.
$$
Besides the fact that the convergence rate of $P_t$ towards $0$ is not exponential,  following the stochastic analysis developed in the present article several
uniform propagation of chaos properties can be developed for this toy model is a rather simple way. The extension of these results to
more general multidimensional diffusions with drift remains an open research question.

The second open question is to analyze the long-time behaviour of the extended EnKF commonly used in nonlinear filtering theory, and more particularly in the numerical solving of data assimilation problems arising in ocean-atmosphere
sciences and oil reservoir simulations.

Another important problem is clearly to develop uniform propagations of chaos properties of the EnKF in discrete time settings.
Last, but not least a series of research projects can be developed around the fluctuations and the large deviations of this new class of mean-field type particle models.

\subsubsection*{Acknowledgements}
We would like to thank Adrian N. Bishop and Sahani Pathiraja. Our discussion in UNSW and UTS in Sydney as well as their detailed comments greatly improved the 
presentation of the article.

\section*{Appendix}

\subsection*{Proof of lemma~\ref{perturbation-lemma-intro}}\label{proof-lemma-intro}

The first assertion is a direct consequence of the inequality
\begin{eqnarray*}
\left\Vert   \Ea_{s,t}(A+B)\right\Vert_2&\leq& 
\exp{\left(\int_s^t\mu(A_u)~du+\int_s^t~\Vert B_u\Vert_2~du\right)}
\end{eqnarray*}
The above estimate is a direct consequence of the matrix log-norm inequality
$$
\mu(A_t+B_t) < \mu(A_t)+\mu(B_t)\quad\mbox{\rm and the fact that}\quad
\mu(B_t)\leq \Vert B_t\Vert_2
$$
This ends the proof of the first assertion. To check the second assertion we observe that
$$
\partial_{t} \Ea_{s,t}(A+B) =\left(\partial_{t} \Ea_{t}(A+B)\right) \Ea_{s}(A+B)^{-1}=A_t\Ea_{s,t}(A+B)+B_t\Ea_{s,t}(A+B)
$$
This implies that
 $$
 \Ea_{s,t}(A+B) =
\Ea_{s,t}(A)+\int_s^t~\Ea_{u,t}(A)~B_u\Ea_{s,u}(A+B)~du
 $$
 for any $s\leq t$ from which we prove that
   \begin{eqnarray*}
 e^{\beta_A(t-s)}\Vert\Ea_{s,t}(A+B) \Vert&\leq &
\alpha_A+\alpha_A~\int_s^t~ e^{\beta_A(t-s)}~ e^{-\beta_A(t-u)}~\Vert B_u\Vert~\Vert
\Ea_{s,u}(A+B)\Vert~du\\
&=&\alpha_A+\alpha_A~\int_s^t~ 
\Vert B_u\Vert~ e^{\beta_A(u-s)}\Vert
\Ea_{s,u}(A+B)\Vert~du
 \end{eqnarray*}
  By Gr\"onwall's lemma this implies that
  $$
   e^{\beta_A(t-s)}\Vert\Ea_{s,t}(A+B) \Vert\leq \alpha_A \exp{\left[\int_s^t \alpha_A\Vert B_u\Vert~du\right]}
  $$
   This ends the proof of the lemma.
  \cqfd

\subsection*{Proof of lemma~\ref{lemma-Lip-Cov}}\label{proof-lemma-Lip-Cov}

This section is mainly concerned with the proof of Lemma~\ref{lemma-Lip-Cov}.  Let $P_1$ and $P_2$ be the covariance matrices
of some $r$-valued random variables $Z_1$ and $Z_2$. Let $\overline{Z}_i$ be an independent copie
of $Z_i$, with $i=1,2$, and set
$$
Z_{1,2}:=Z_1-Z_2\quad \mbox{\rm and}\quad
\overline{Z}_{1,2}:=\overline{Z}_1-\overline{Z}_2.
$$

Observe that
$$
Z_1-\overline{Z}_1=\left(Z_2-\overline{Z}_2\right)+\left(Z_{1,2}-\overline{Z}_{1,2}\right)
$$
and for any $i=1,2$ we have
$$
2~P_i=\EE\left(\Za_i\right)\quad
\mbox{\rm with the random matrix}\quad \Za_i:= \left(Z_i-\overline{Z}_i\right)\left(Z_i-\overline{Z}_i\right)^{\prime}.
$$
This yields the decomposition
$$
\begin{array}{l}
\Za_1-\Za_2\\
\\
=\left(\left(Z_2-\overline{Z}_2\right)+\left(Z_{1,2}-\overline{Z}_{1,2}\right)\right)\left(\left(Z_2-\overline{Z}_2\right)+\left(Z_{1,2}-\overline{Z}_{1,2}\right)\right)^{\prime}-\left(Z_2-\overline{Z}_2\right)\left(Z_2-\overline{Z}_2\right)^{\prime}\\
\\
=\left(Z_2-\overline{Z}_2\right)\left(Z_{1,2}-\overline{Z}_{1,2}\right)^{\prime}+\left(Z_{1,2}-\overline{Z}_{1,2}\right)\left(Z_2-\overline{Z}_2\right)^{\prime}+\left(Z_{1,2}-\overline{Z}_{1,2}\right)\left(Z_{1,2}-\overline{Z}_{1,2}\right)^{\prime}.
\end{array}$$
This shows that
\begin{eqnarray*}
2~(P_1-P_2)&=&\EE\left(\left(Z_2-\overline{Z}_2\right)\left(Z_{1,2}-\overline{Z}_{1,2}\right)^{\prime}\right)+
\EE\left(\left(Z_{1,2}-\overline{Z}_{1,2}\right)\left(Z_2-\overline{Z}_2\right)^{\prime}\right)\\
&&\hskip4cm+
\EE\left(\left(Z_{1,2}-\overline{Z}_{1,2}\right)\left(Z_{1,2}-\overline{Z}_{1,2}\right)^{\prime}\right)
\end{eqnarray*}
from which we prove the trace formula
\begin{eqnarray*}
2~\mbox{\rm tr}(P_1-P_2)&=&2~\EE\left(\left\langle Z_{1,2}-\overline{Z}_{1,2},Z_2-\overline{Z}_2\right\rangle\right)+\EE\left(\Vert Z_{1,2}-\overline{Z}_{1,2}\Vert^2\right).
\end{eqnarray*}

By Cauchy-Schwarz inequality, we find that
\begin{eqnarray*}
2~\left\vert\mbox{\rm tr}(P_1-P_2)\right\vert&\leq &2~\EE\left(\Vert Z_{1,2}-\overline{Z}_{1,2}\Vert^2\right)^{1/2}
\EE\left(\Vert Z_{2}-\overline{Z}_{2}\Vert^2\right)^{1/2}+\EE\left(\Vert Z_{1,2}-\overline{Z}_{1,2}\Vert^2\right)\\
&\leq &8~\EE\left(\Vert Z_{1,2}\Vert^2\right)^{1/2}
\EE\left(\Vert Z_{2}\Vert^2\right)^{1/2}+4~\EE\left(\Vert Z_{1,2}\Vert^2\right).
\end{eqnarray*}
This yields
$$
\frac{1}{2}~\left\vert\mbox{\rm tr}(P_1-P_2)\right\vert\leq 2~
\EE\left(\Vert Z_{1}-Z_2\Vert^2\right)^{1/2}
\EE\left(\Vert Z_{2}\Vert^2\right)^{1/2}+\EE\left(\Vert Z_{1}-Z_2\Vert^2\right)
$$
from which we find that
$$
\frac{1}{2}~\left\vert\mbox{\rm tr}(\Pa_{\eta_0})-\mbox{\rm tr}(\Pa_{\check{\eta}_0})\right\vert\leq 2~\WW_2(\eta_0,\check{\eta}_0)
~\Vert\check{\eta}_0(e_2)\Vert^{1/2}+\WW_2(\eta_0,\check{\eta}_0)^2.
$$

In much the same way, we have
$$
\begin{array}{l}
\displaystyle 4~\Vert P_1-P_2\Vert^2_F=\sum_{1\leq i,j\leq r}~\left[\EE\left(\Big(Z_{1,2}-\overline{Z}_{1,2}\Big)(i)\Big(Z_{1,2}-\overline{Z}_{1,2}\Big)(j)\right)\right.\\
\\
\left.\hskip1cm+\displaystyle
\EE\left(\Big(Z_{2}-\overline{Z}_{2}\Big)(i)\Big(Z_{1,2}-\overline{Z}_{1,2}\Big)(j)\right)
+\EE\left(\Big(Z_{2}-\overline{Z}_{2}\Big)(j)\Big(Z_{1,2}-\overline{Z}_{1,2}\Big)(i)\right)\right]^2.
\end{array}
$$
This implies that
$$
\begin{array}{l}
\displaystyle 2~\Vert P_1-P_2\Vert^2_F=\sum_{1\leq i,j\leq r}~\left[\left\{\EE\left(\Big(Z_{1,2}-\overline{Z}_{1,2}\Big)(i)\Big(Z_{1,2}-\overline{Z}_{1,2}\Big)(j)\right)\right\}^2\right.\\
\\
\left.\hskip1cm+\displaystyle
\left\{\EE\left(\Big(Z_{2}-\overline{Z}_{2}\Big)(i)\Big(Z_{1,2}-\overline{Z}_{1,2}\Big)(j)\right)\right\}^2
+\left\{\EE\left(\Big(Z_{2}-\overline{Z}_{2}\Big)(j)\Big(Z_{1,2}-\overline{Z}_{1,2}\Big)(i)\right)\right\}^2\right].
\end{array}
$$
Using Cauchy-Schwarz inequality we prove that
$$
\begin{array}{l}
\displaystyle 2~\Vert P_1-P_2\Vert^2_F
\displaystyle\leq 2~
\sum_{1\leq i\leq r}\EE\left((Z_{1,2}-\overline{Z}_{1,2})(i)^2\right)\sum_{1\leq j\le r}
\EE\left((Z_{2}-\overline{Z}_{2}(j)^2\right)\\
\\
\displaystyle\hskip7cm+\left[
\sum_{1\leq i\leq r}\EE\left((Z_{1,2}-\overline{Z}_{1,2})(i)^2\right)\right]^2.
\end{array}
$$
This implies that
\begin{eqnarray*}
\displaystyle \sqrt{2}~\Vert P_1-P_2\Vert_F
\displaystyle&\leq& \sqrt{2}~ \EE\left(\Vert Z_{1,2}-\overline{Z}_{1,2}\Vert^2\right)^{1/2}
\EE\left(\Vert Z_{2}-\overline{Z}_{2}\Vert^2\right)^{1/2}+\EE\left(\Vert Z_{1,2}-\overline{Z}_{1,2}\Vert^2\right)\\
\displaystyle&\leq& 4\sqrt{2}~ \EE\left(\Vert Z_{1,2}\Vert^2\right)^{1/2}~
\EE\left(\Vert Z_{2}\Vert^2\right)^{1/2}+4~\EE\left(\Vert Z_{1,2}\Vert^2\right).
\end{eqnarray*}
We find that
$$~
\Vert \Pa_{\eta_0}-\Pa_{\check{\eta}_0}\Vert_F\leq\WW_2(\eta_0,\check{\eta}_0)
~\Vert\check{\eta}_0(e_2)\Vert^{1/2}+~ \frac{1}{\sqrt{2}}~\WW_2(\eta_0,\check{\eta}_0)^2.
$$
This ends the proof of the lemma.\cqfd

\section*{Divergence regions - 2d observers}\label{some-illustrations-2d}

We illustrate the spectral analysis discussion given in section~\ref{sous-section-section-regularity} with  $2$-dimensional partially
observed filtering problems associated with the parameters
$$
(r_1,r_2)=(2,1),~
C=[1,0]~,~(R_1,R_2)=(Id,1)$$
and some unstable drift matrix $A$ with a saddle equilibrium, that is
$$
-\mbox{\rm det}(A)=A_2\wedge A_1>0
$$
with the column vectors $A_1=\left(\begin{array}{c}A_{1,1}\\A_{1,2}\end{array}
\right)$ and $A_2=\left(\begin{array}{c}A_{2,1}\\A_{2,2}\end{array}
\right)$. In the above display $A_2\wedge A_1$ stands for the cross product of the vectors $A_2$ and $A_1$.
Whenever $A_{1,2}\not=0$ the system is observable and controllable; thus there exists some unique steady state $P$ and
$\mbox{\rm Re}(\lambda_{\tiny max}(\overline{A}))<0$, or equivalently
$$
 \mbox{\rm tr}(A)<P_{1,1}\quad \mbox{\rm and}\quad 
A_2\wedge P_1>A_2\wedge A_1
$$ 
The set of admissible fluctuation matrices $Q$ (that is s.t.
$P+Q\geq 0$) is defined by
$$
\begin{array}{rcl}
\Qa(P)&:=&\left\{Q\in \SS_{r_1}~:~
  \left[P_1-\frac{1}{2}~Q_1\right]\wedge Q_2+\left[P_2+\frac{1}{2}~Q_2\right]\wedge Q_1 <P_1\wedge P_2
\right\}\\
&&\\
&&\hskip3cm\cap \left\{Q\in \SS_{r_1}~:~Q_{1,1}\in [-P_{1,1},\infty[~~\mbox{\rm and}~~
Q_{2,2}\in [-P_{2,2},\infty[
\right\}
\end{array}$$
 Given some $Q\in \Qa(P)$ several cases can happen. In the most favorable case,  we have
$$
(P_{1,1}-\mbox{\rm tr}(A))>- Q_{1,1}~~\mbox{\rm and}~~\mbox{\rm det}([(A_1-P_1)-Q_1,A_2])>0
\Longleftrightarrow \mbox{\rm Re}(\lambda_{\tiny max}(\overline{A}-QS))<0
$$
The determinant condition is equivalent to
$$
Q_1\wedge A_2<
A_2\wedge P_1-A_2\wedge A_1
$$
To be more precise, we have
 two negative eigenvalues when
$\mbox{\rm tr}(\overline{A}-QS)^2>4\mbox{\rm det}(\overline{A}-QS)$, otherwise we have a spiral phase portrait 
with complex eigenvalues with negative real parts.
In both cases the matrix $\overline{A}-QS$ remains stable.
Skipping the discussion on borderline cases, the other situation that may arise is that 
$$
Q_1\wedge A_2>
A_2\wedge P_1-A_2\wedge A_1
\quad\mbox{\rm or}\quad 0<P_{1,1}-\mbox{\rm tr}(A)< - Q_{1,1}
$$
When the l.h.s. condition is met, the eigenvalues have opposite sign and the stochastic observer experience 
a catastrophic divergence in the direction of the eigenvector associated with the positive one. When the r.h.s. 
condition is met, both eigenvalues are negative real numbers
if $\mbox{\rm tr}(\overline{A}-QS)^2>4\mbox{\rm det}(\overline{A}-QS)$, otherwise they are both 
complex with the same negative real part. In both situations the observer diverges. The divergence set (\ref{global-div-region-def})
is given by
\begin{eqnarray*}
 \Qa_{\mbox{\tiny div}}&=&\left\{Q\in \Qa(P)~:~
 Q_1\wedge A_2>
A_2\wedge P_1-A_2\wedge A_1
 \right\}\\
 &&\hskip2cm\cup \left\{Q\in \Qa(P)~:~Q_{1,1}<-(P_{1,1}-\mbox{\rm tr}(A))
 \right\}
\end{eqnarray*}
For instance for 
$$
A=\left[
\begin{array}{cc}
1&2\\
1&3
\end{array}
\right]\Longrightarrow P\simeq\left[
\begin{array}{cc}
8.7&14.5\\
14.5&30
\end{array}
\right]\quad
\mbox{\rm and}\quad 
 \overline{A}\simeq\left[
\begin{array}{cc}
-7.4&2\\
-13.5&3
\end{array}
\right]
$$
Notice that the system is locally ill-conditioned as  $\lambda_{\tiny max}(\overline{A})=5.4913>0$. 
The set of admissible fluctuations is given by
$$
\begin{array}{rcl}
\Qa(P)&:=&\left\{Q\in \SS_{r_1}~:~
Q_{1,2}Q_{2,1}+14.5 (Q_{1,2}+Q_{2,1})<52.2+8.7Q_{2,2}+30Q_{1,1}+Q_{1,1}Q_{2,2} 
\right\}\\
&&\\
&&\hskip3cm\cap \left\{Q\in \SS_{r_1}~:~Q_{1,1}\in [-8.7,\infty[~~\mbox{\rm and}~~
Q_{2,2}\in [-30,\infty[
\right\}
\end{array}$$

In this situation we have
\begin{eqnarray*}
 \Qa_{\mbox{\tiny div}}&\simeq&\left\{Q\in \Qa(P)~:~Q_{1,2}<1.5 Q_{1,1}-1.8
 \right\}\cup \left\{Q\in \Qa(P)~:~Q_{1,1}\in \left] -8.7, -4.7\right[
 \right\}
\end{eqnarray*}
In this situation the stable subset of $\Qa(P)$ is given by 
$$
Q\in \Qa(P)~:~
Q_{1,2}>1.5 Q_{1,1}-1.8
\quad\mbox{\rm and}\quad Q_{1,1}\in \left]-4.7,\infty\right[
$$

\begin{figure}[h]
\centering
    \vskip-3cm
  \includegraphics[trim=1cm 1cm 3cm 2cm, clip=true,width=0.5\textwidth]{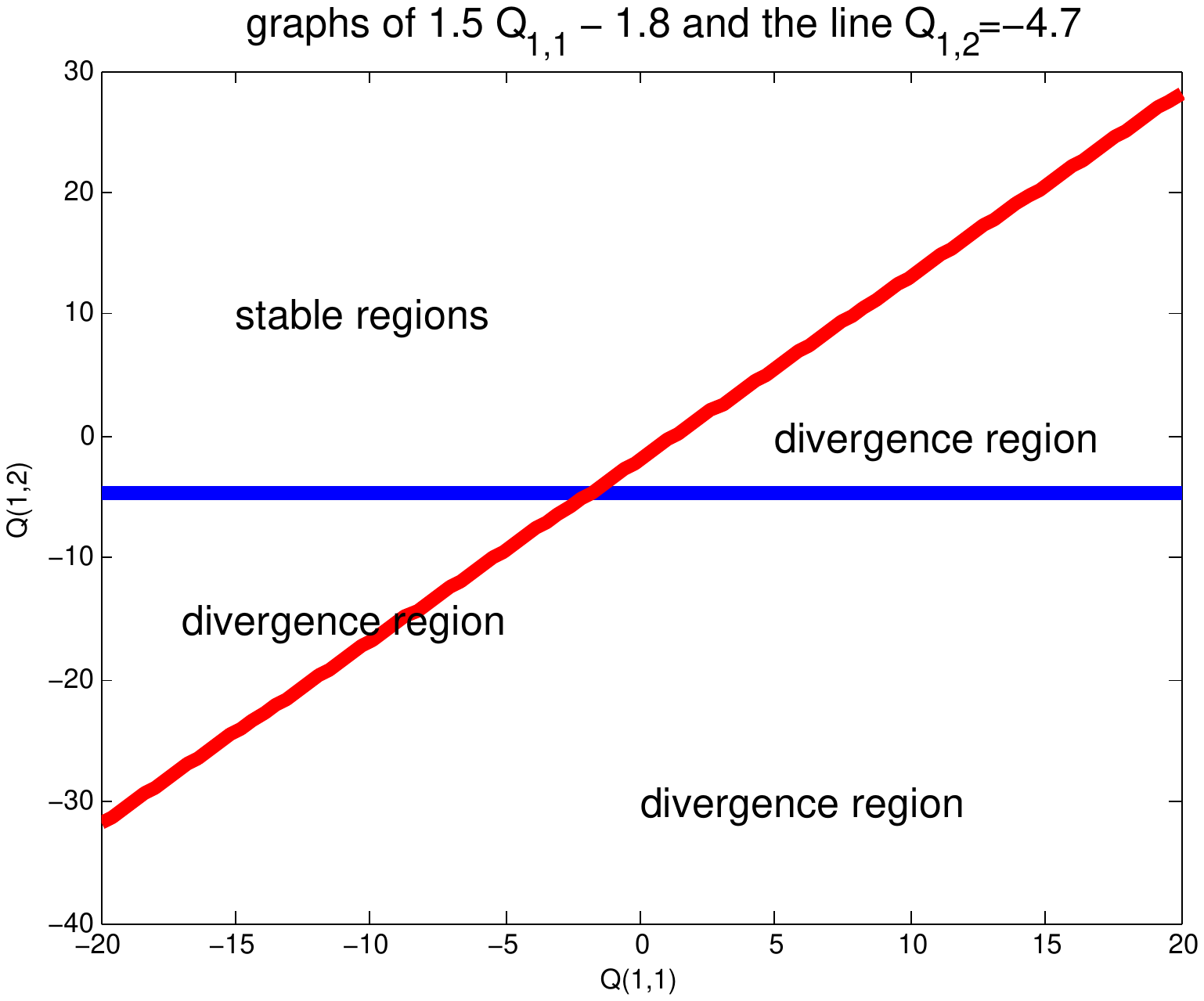}
    \vskip-3cm
      \caption{fluctuation-divergence effects}  
      \label{figure-convergence-divergence}
 \end{figure}

These convergence and divergence sets without the admissible conditions are illustrated in figure~\ref{figure-convergence-divergence}.

\begin{figure}[h]
\centering
  \includegraphics[trim=0cm 0cm 0cm 0cm, clip=true,width=0.6\textwidth]{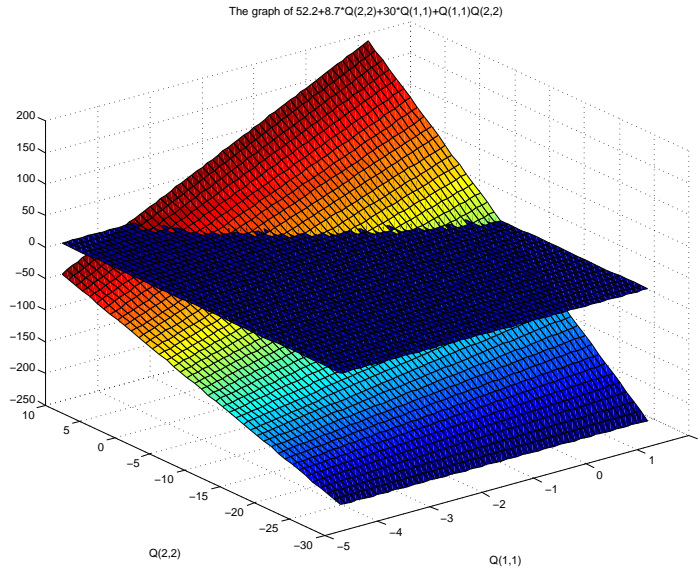}
    \vskip-2cm
      \caption{Trace of the stable domain with diagonal matrices}  
      \label{figure-trace-stable}
 \end{figure}

The trace of the divergence domain with diagonal matrices (i.e. $Q_{1,2}=0=Q_{2,1}$) resume to
diagonal matrices s.t.
$$
 Q_{1,1}\in \left] 1.2, \infty\right[\quad
\mbox{\rm and}\quad
0<52.2+8.7Q_{2,2}+30Q_{1,1}+Q_{1,1}Q_{2,2} 
$$
The trace with the stable domain is 
$$
Q_{1,1}\in \left]-4.7,1.2\right[\quad \mbox{\rm and}\quad0<52.2+8.7Q_{2,2}+30Q_{1,1}+Q_{1,1}Q_{2,2} 
$$
An illustration of this set  is given in figure~\ref{figure-trace-stable}

\begin{figure}[h]
\centering
  \includegraphics[trim=3cm 18cm 3cm 2cm, clip=true,width=0.6\textwidth]{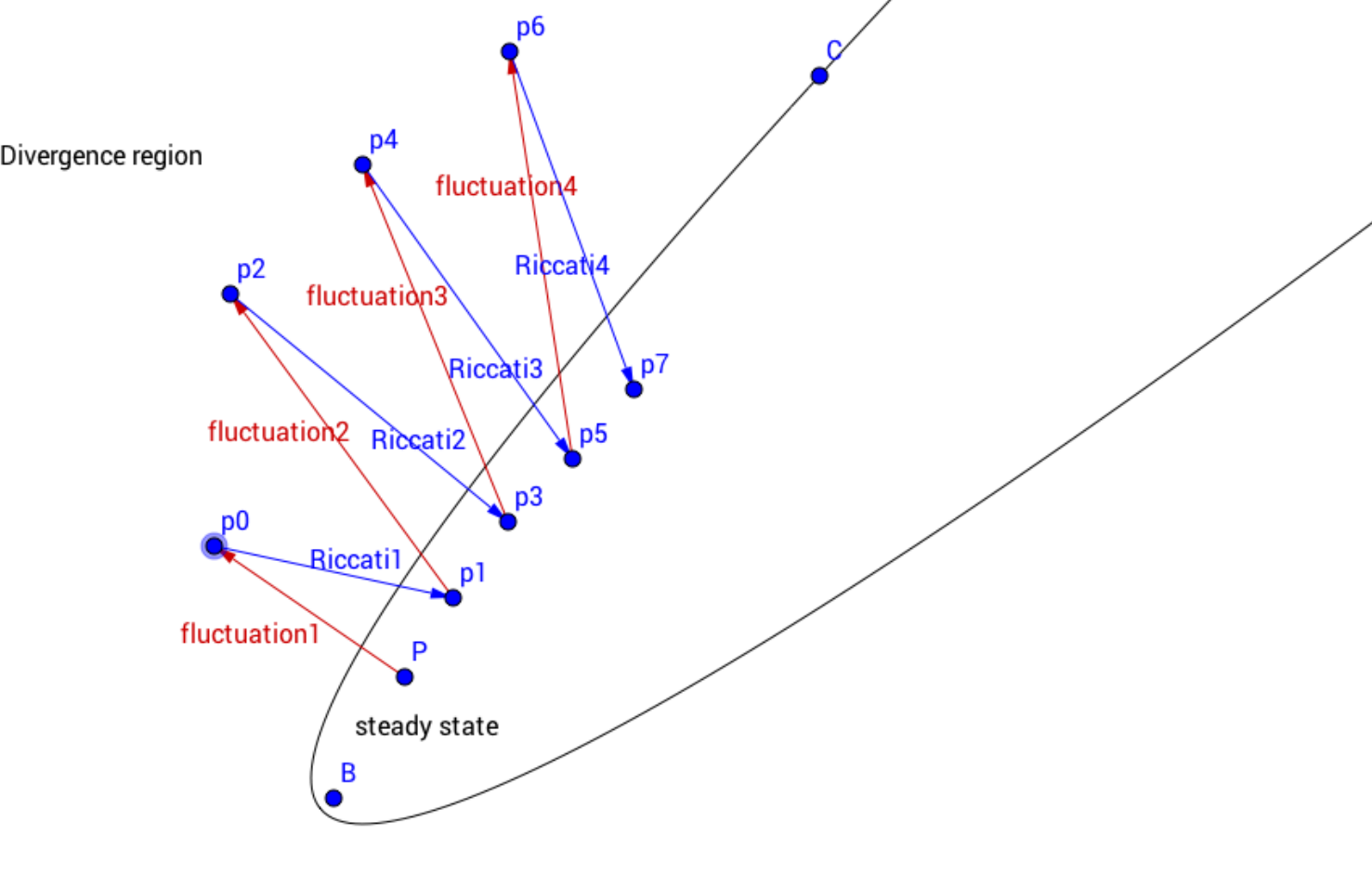}
  \vskip-.3cm
      \caption{fluctuation-divergence effects}  
      \label{figure-divergence}
 \end{figure}

\begin{figure}[h]
  \begin{subfigure}[b]{0.5\linewidth}
    \centering
\includegraphics[trim = 3mm 5mm 4mm 1mm, clip, width=1.1\textwidth]{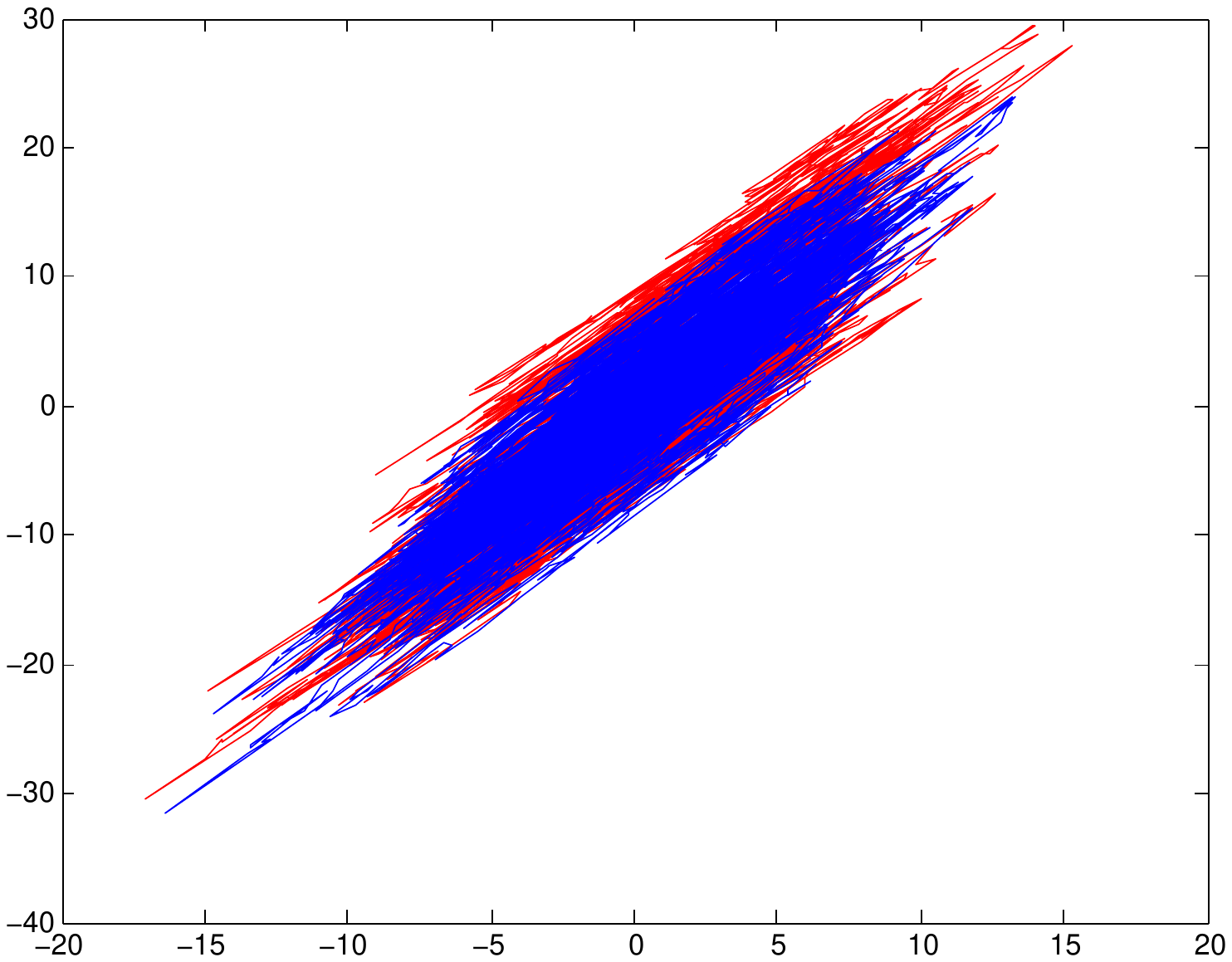} 
\vskip-2.5cm
    \caption{\begin{tabular}{l}Observer [red] /Steady state [blue]\\ $Q_{1,1}=1.1$ with $t=100$\end{tabular}} 
    \label{fig11}
    \vspace{4ex}
  \end{subfigure}
  \begin{subfigure}[b]{0.5\linewidth}
    \centering
    \includegraphics[width=1.1\linewidth]{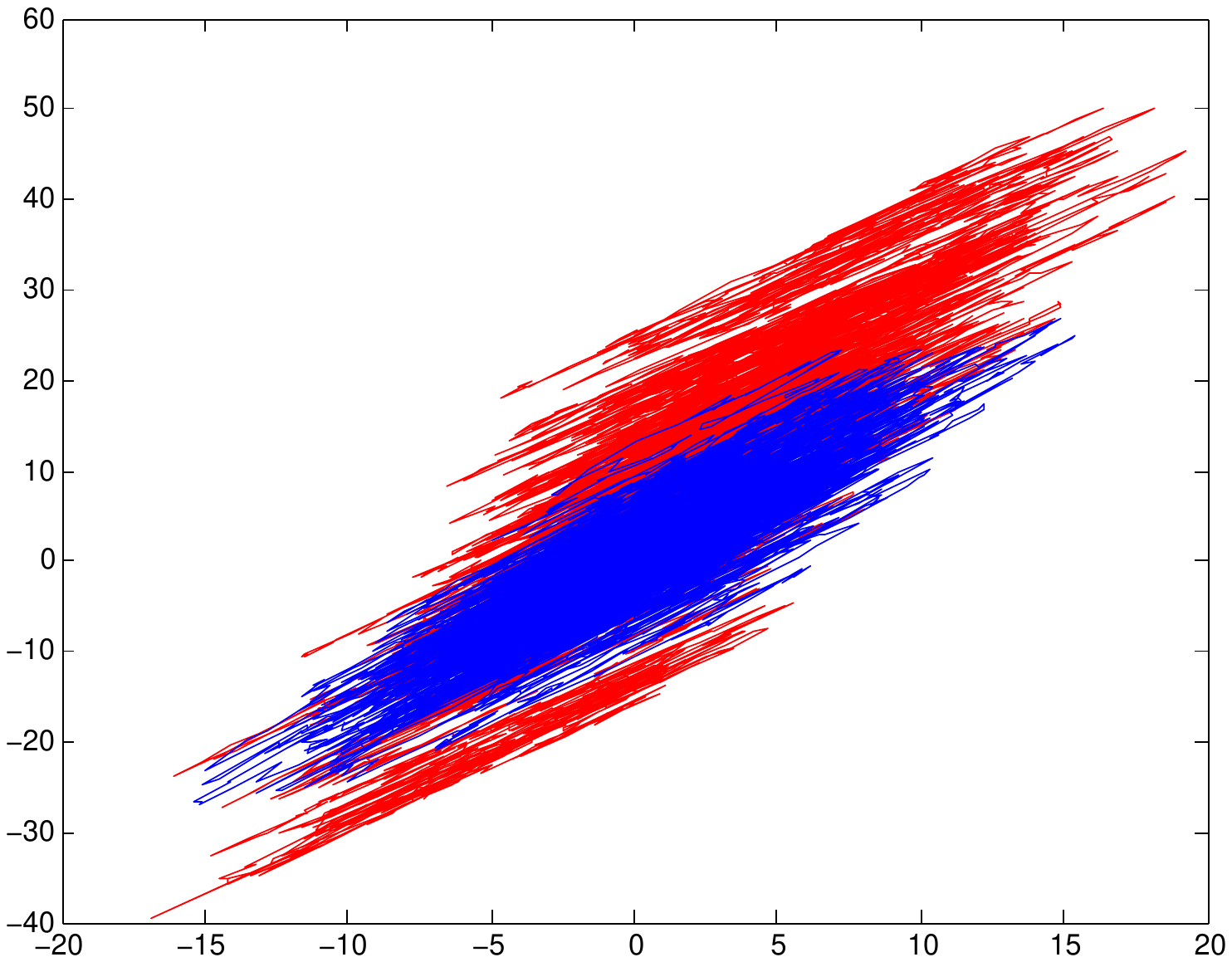} 
    \vskip-2.5cm
    \caption{\begin{tabular}{l}Observer [red] /Steady state [blue] \\
    $Q_{1,1}=1.2$ with $t=100$\end{tabular}} 
    \label{fig12}
    \vspace{4ex}
  \end{subfigure} 
  \begin{subfigure}[b]{0.5\linewidth}
  \vskip-3cm
    \centering
    \includegraphics[width=1.1\linewidth]{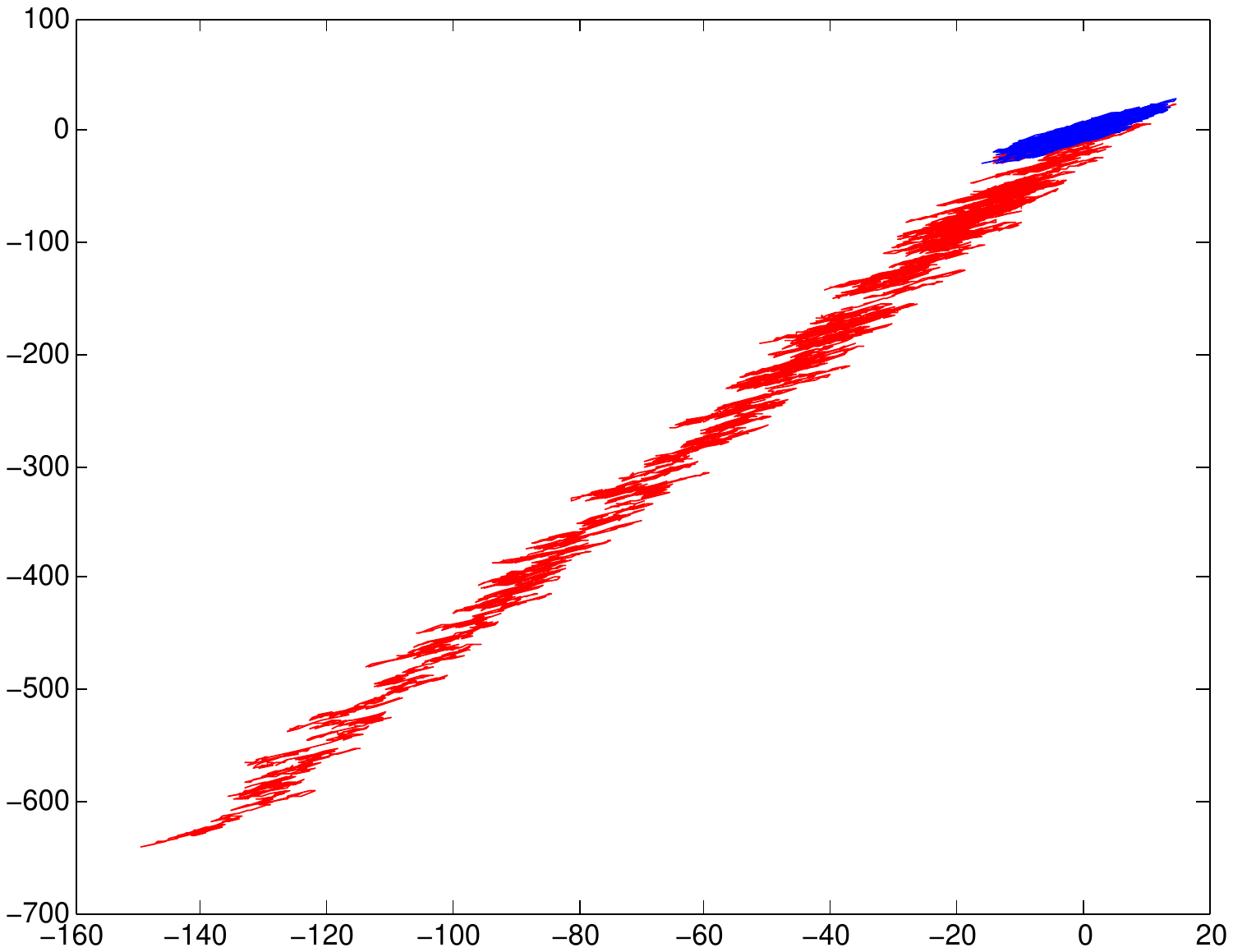} 
    \vskip-3cm
    \caption{\begin{tabular}{l}Observer [red] /Steady state [blue]\\ $Q_{1,1}=1.3$ with $t=125$\end{tabular}} 
    \label{fig13} 
  \end{subfigure}
  \begin{subfigure}[b]{0.5\linewidth}
  \vskip-3cm
    \centering
    \includegraphics[width=1.1\linewidth]{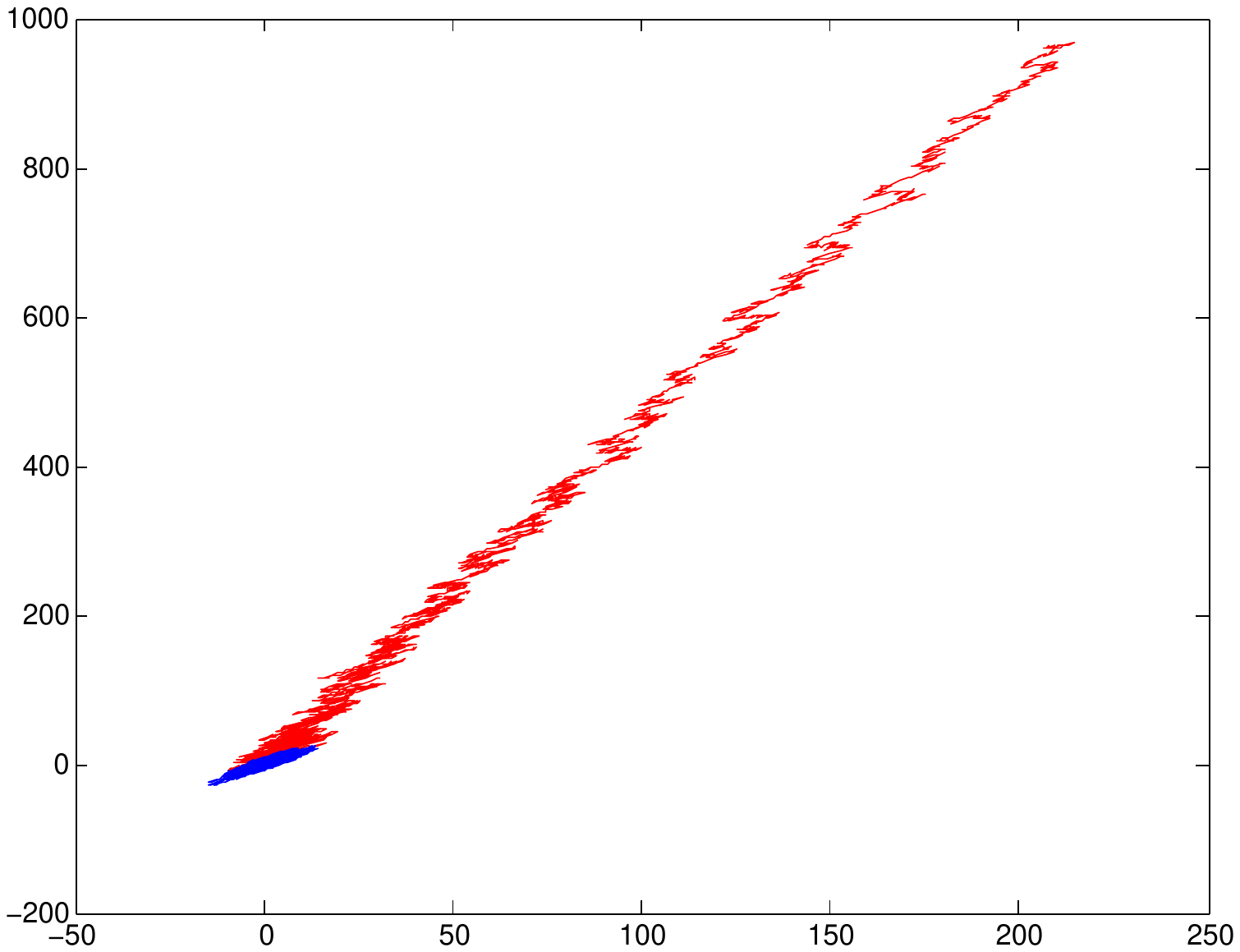} 
    \vskip-3cm
    \caption{\begin{tabular}{l}Observer [red] /Steady state [blue]\\
    $Q_{1,1}=1.4$ with $t=75$\end{tabular}} 
    \label{fig14} 
  \end{subfigure}
\end{figure}

The fluctuation/divergence effects we can expect when the observer is driven by 
fluctuations entering into the divergence domain are illustrated in figure~\ref{figure-divergence} 
 A series of realization of the stochastic observer driven by fluctuation matrices 
 in the stable domain are presented in figures~\ref{fig11},\ref{fig12}; the
 ones driven by fluctuation matrices in the divergence set are presented  in figures~\ref{fig13},\ref{fig14}.
The entries $Q_{2,2}$ are not seen by the observer so we assume that $Q_{2,2}=0$.

\newpage
\section*{Figures}

\end{document}